\documentclass[a4paper,11pt]{article} 
\usepackage[utf8]{inputenc}
\usepackage[T1]{fontenc} 
\usepackage{amsmath,amsthm,amsfonts,amssymb,euscript,stmaryrd,graphicx,enumitem,yhmath,mathrsfs}
\usepackage[toc,page]{appendix}
\usepackage{algorithm,algorithmic}
\usepackage{fullpage}
\usepackage{color,psfrag}
\usepackage{dsfont}
\usepackage[babel=true]{csquotes} 
\usepackage{lmodern}
\usepackage{bibentry}

\usepackage[dvipsnames]{xcolor}

\usepackage[french,english]{babel}
\usepackage{ifthen}
\usepackage{array}
\usepackage{hyperref} 
\usepackage{verbatim}
\usepackage[normalem]{ulem}
 \usepackage{soul}
\usepackage{empheq}

\newcommand{\footnoteremember}[2]{
 \footnote{#2}
 \newcounter{#1}
 \setcounter{#1}{\value{footnote}}
}

\author
{        {Joachim \textsc{Lebovits}\footnote{MAThematic Center of Heidelberg, University of Heidelberg, INF 294,
69120 Heidelberg, Germany.}\footnoteremember{footnoteP13}{Laboratoire Analyse, Géométrie et Applications, C.N.R.S. (UMR 7539), Université Paris 13, Sorbonne Paris Cité, 99 avenue Jean-Baptiste Clément
93430, Villetaneuse, France. Email address: \url{jolebovits@gmail.com}}}
  } 
\newtheorem{theo}{Theorem}[section]
\newtheorem{theodef}{Theorem-Definition}[section]
\newtheorem{prop}[theo]{Proposition}
\newtheorem{coro}[theo]{Corollary}

\newtheorem{rem}{Remark}
\newtheorem{lem}[theo]{Lemma}
\newtheorem{defi}{Definition}

\newtheorem{ex}[theo]{Example}

\newenvironment{pr}{\begin{proof}[\bfseries \textup{Proof.}]}{\end{proof}}

\newenvironment{g}[2]{\begin{array}{cl} %
<\hspace{-0.2cm}<\hspace{-0.1cm} #1, \hspace{-0.3cm} &#2 \hspace{-0.1cm} >\hspace{-0.2cm}>}%
{\end{array}}

\newcommand{\bD}{\mathbb{D}}

\newcommand{\bN}{\mathbb{N}}

\newcommand{\sifbm}{\mathbf{B}}

\newcommand{\E}{\mathbf{E}}

\newcommand{\N}{\mathbf{N}}
\newcommand{\Q}{\mathbf{Q}}
\newcommand{\R}{\mathbf{R}}
\newcommand{\C}{\mathbf{C}}
\newcommand{\zZ}{\mathbf{Z}}

\newcommand{\oS}{${(\cS)}^*$}
\newcommand{\ooS}{{(\cS)}^*}

\newcommand{\Z}{\mathcal{Z}}

\newcommand{\w}{\Gamma_{\hspace{-0.1cm}\scriptscriptstyle R}}

\newcommand{\z}{\mathcal{Z}^{\scriptscriptstyle \hspace{-0.01cm} T}_{\hspace{-0.05cm}\scriptscriptstyle R}}

\newcommand{\di}{d^{\diamond} \hspace{-0.05cm}}

\newcommand{\ds}{\displaystyle}

\newcommand{\cl}{\centerline}

\def\i1{\mathbf{1}}  

\newcommand{\cS}{\mathcal{S}}
\newcommand{\sS}{\mathscr{S}}

\def\Aa{\ref{Aa} }
\def\Aap{\ref{Aa}}

\def\Aai{\ref{Ai} }

\def\Aaii{\ref{Aii} }
\def\Aaiip{\ref{Aii}}

\def\Aaiii{\ref{Aiii} }

\def\Aaiv{\ref{Aiv} }
\def\Aaivp{\ref{Aiv}}

\def\Da{\ref{D} }
\def\Dap{\ref{D}}

\def\Dai{\ref{Di} }
\def\Daip{\ref{Di}}

\def\Daii{\ref{Dii} }

\def\Di{{\mathscr{D}_{\text{(i)}}}}

\def\iIp{{\eqref{i}} }
\def\iIpp{{\eqref{i}}}

\def\iIlp{{\eqref{ipq}} }

\def\iI{{(\mathcal{I})}}

\def\ie{{\textit{i.e. }}}

\def\rR{\mathscr {R}}

\def\cB{{\cal B}}
\def\cC{C}

\def\cE{{\cal E}}
\def\cG{{\mathcal G}}
\def\cH{{\cal H}}
\def\cI{{\cal I}}

\def\cK{{\cal K}}
\def\cL{{\cal L}}

\def\cN{{\cal N}}

\def\cP{{\cal P}}

\def\cR{{\cal R}}
\def\cS{{\cal S}}
\def\cT{{\cal T}}
\def\cU{{\cal U}}

\newcommand{\Dom}{\bD\text{om}}

\newcommand{\cF}{\mathcal{F}}

\newcommand{\bit} {\begin{itemize} }
\newcommand{\eit} {\end{itemize} }
\def\ben{\begin{enumerate}}
\def\een{\end{enumerate}}

\def\ia{\item[H1)]}
\def\iaa{\item[H2)]}
\def\iaaa{\item[H3)]}

\def\iti{\item[(i)]}

\def\ita{\item[a)]}
\def\itap{\item[${\text{a}}_{1})$]}
\def\itapp{\item[${\text{a}}_{2})$]}
\def\itappp{\item[${\text{a}}_{3})$]}

\def\itappq{\item[$a'_{2})$]}

\def\itaps{\item[${\text{b}}_{1})$]}
\def\itapps{\item[${\text{b}}_{2})$]}
\def\itappps{\item[${\text{b}}_{3})$]}
\def\itapppps{\item[${\text{b}}_{4})$]}
\def\itappqs{\item[$b'_{2})$]}

\def\itii{\item[(ii)]}

\def\itb{\item[b)]}

\def\itiii{\item[(iii)]}

\def\ {\hspace{0.1cm}}

\def\keywordname{{\bf Keywords:}}

\newcommand{\keywords}[1]{\par\addvspace\baselineskip\noindent\keywordname\enspace\ignorespaces#1}

{\end{array}}

\makeatletter
\renewcommand\theequation{\thesection.\arabic{equation}}
\@addtoreset{equation}{section}
\makeatother

\graphicspath{{./}{figures/}}

\DeclareFontFamily{U}{mathx}{\hyphenchar\font45}
\DeclareFontShape{U}{mathx}{m}{n}{
      <5> <6> <7> <8> <9> <10>
      <10.95> <12> <14.4> <17.28> <20.74> <24.88>
      mathx10
      }{}
\DeclareSymbolFont{mathx}{U}{mathx}{m}{n}
\DeclareFontSubstitution{U}{mathx}{m}{n}
\DeclareMathAccent{\widecheck}{0}{mathx}{"71}
\DeclareMathAccent{\wideparen}{0}{mathx}{"75}

\title{%
    \begin{minipage}\linewidth
        \centering\bfseries\sffamily
     Stochastic Calculus with respect to Gaussian processes
            \end{minipage}
}

\begin{document}

\maketitle

\frenchspacing

\vspace{-5ex}

\begin{abstract}
\begin{footnotesize}
\par Stochastic integration \textit{wrt} Gaussian processes has raised strong interest in recent years, motivated in particular by its applications in  Internet traffic modeling, biomedicine and finance.  The aim of this work is to define and develop a White Noise Theory-based anticipative stochastic calculus with respect to 
all Gaussian processes that have an integral representation over a real (maybe infinite)
interval. 
Very rich,  
this class  of Gaussian processes contains, among many others, 
Volterra processes (and thus fractional Brownian motion) as well as processes 
 the regularity of which varies along the time (such as multifractional Brownian motion).
A systematic comparison of the stochastic calculus (including Itô formula) we provide here, to the ones given by Malliavin calculus in  \cite{nualart,MV05,NuTa06,KRT07,KrRu10,LN12,SoVi14,LN12}, and by Itô stochastic calculus is also made. Not only our stochastic calculus fully generalizes and extends the ones originally proposed  in 
\cite{MV05} and in \cite{NuTa06} for Gaussian processes, but also the ones proposed in \cite{ell,bosw,ben1} for fractional Brownian motion (\textit{resp.} in \cite{JLJLV1,JL13,LLVH} for multifractional Brownian motion).
\end{footnotesize}
\end{abstract}

\vspace{-3ex}


\smallskip



\vspace{-2ex}

\section{Introduction}
\vspace{-1ex}
The purpose of this paper is to develop an anticipative stochastic calculus with 
respect to Gaussian process  $G:={(G_{t})}_{t\in\rR}$ 
that can be written under the form:
\vspace{-1ex}
\begin{equation}
\label{erofkorekp}
G_{t}=\int_{\R} g_{t}(u) \ dB_{u},
\vspace{-1ex}
\end{equation}
where $\R$ denotes the set of real numbers, $\rR$ denotes a closed 
interval of $\R$ (that may be equal to $\R$), $B:={(B_{u})}_{u
\hspace{-0.01cm}\in\R}$ is Brownian motion on $\R$ and  ${(g_{t})}_{t\in
\rR}$  is a family of a square integrable functions\footnote{\ie such that for all $t$ in $\rR$,  $u\mapsto g_{t}(u)$ is measurable on $\R$ and such that $\int_{\R} {|g_{t}(u)|}^{2} \ du<+\infty$.} on $\R$. 
Denote $\mathscr{G}$ the set of Gaussian processes that can be written under the form \eqref{erofkorekp}. This class  of Gaussian processes contains, among many others, 
Volterra processes (and thus fractional Brownian motion),
Gaussian Fredholm processes as well as processes 
the regularity of which varies along the time (such as multifractional Brownian motion). 
For every positive real $T$, the process ${(V_{t})}_{t\in [0,T]}$ 
 is said to be a Volterra 
 process on $[0,T]$ (\textit{resp.} a Fredholm process), if it can be written under the form: 
\vspace{-1ex}
\begin{equation}
\label{erofkoredefefekdeddedep}  
V_{t}:=\int^{t}_{0} K(t,s) \ dW_{s}; \hspace{0.5cm} \hspace{3ex} (\textit{resp.} \ F_{t}:=\int^{T}_{0} K_{T}(t,s) \ dW_{s}),  \hspace{3ex} \forall \ t\in [0,T].
\vspace{-1ex}
\end{equation}
where  ${(W_{s})}_{s\in [0,T]}$ is a Brownian motion and $K$  
belongs to $L^{2}({[0,T]}^{2},ds)$.
Note moreover that $\mathscr{G}$ also contains the Gaussian processes that can be written under the form:
\begin{equation}
\label{ooijoierfkfpopfokerfpoefoierfu} 
H_{t}:=\int^{t}_{-\infty} K(t,s) \ dB_{s}; \hspace{0.5cm} \forall\ t \in \R.
\vspace{-1ex}
\end{equation}
%
Our main result is an Itô formula, that reads: 
\vspace{-1ex}
\bit
\item for every $T>0$ and every $C^{1,2}([0,T]\times \R)$ function $f$, with sub exponential growth: 
\eit
\begin{equation*}
f(T,G_{T}) = f(0,0) + \int^T_0 \ \tfrac{\partial f}{\partial t}
(t,G_{t}) \ dt + \int^T_0 \ \tfrac{\partial f}{\partial x}(t,G_{t}) 
\ \di \hspace{-0.01cm}G_{t}
				+ \tfrac{1}{2} \ \int^T_0 \ \tfrac{\partial^2 f}{\partial x^2}(t,G_{t}) \  dR_{t},
\end{equation*}

where the equality holds in $L^2(\Omega)$ and almost surely, where $t\mapsto R_{t}$ denotes the variance function of $G$, which will be supposed to be a continuous function, of bounded variations; the meaning of the different terms will be explained below. The Itô formula we provide here is, at our best 
knowledge, one of the most general one for Gaussian processes that are not semimartingales. 
Itô stochastic calculus provides a non-anticipative stochastic integral \textit{wrt} 
semimartingales. However, Itô's theory does not apply 
anymore when the Gaussian process considered is not a 
semimartingale.
Two\footnote{The enlargement of filtration technique is a 
third method to extend Itô integral for non semimartingale 
(see \cite{YorMansuy} and references therein). However we 
will not discuss it  in this paper since it is very rarely used in 
the literature.} main and parallel ways have been 
developed over the years to build a stochastic calculus with 
respect to Gaussian processes; the Itô integral \textit{wrt} 
Brownian motion
 being at the intersection of all these 
approaches. Precisely, one has the:

-\  \textit{trajectorial or pathwise extensions,}

-\  \textit{functional extensions.} 

The trajectorial approach, initiated by \cite{Yo36}, provides generalizations of the Riemann–Stieltjes integral that are: the pathwise forward-type Riemann–Stieltjes integral (introduced in \cite{Fo81}; see also \cite{SoVi13} and references therein) and 
the pathwise generalized Lebesgue–Stieltjes integrals (introduced in \cite{Z}). The reader interested in this approach, that also provides Itô formulas, will find in  \cite{FV} a very complete overview. Let us also mention the stochastic calculus via regularization (see \cite{RV07} and references therein), that is also a generalization of Itô integral. 
Since the \textit{pathwise extensions} of Itô integral require, by their very 
definition, that the stochastic integral is built $\omega$ by 
$\omega$, it will clearly appear that they are of a completely 
different nature from our definition of stochastic integral 
(that will be given in Definition \ref{oezifhherioiheroiuh} below). For this reason we will not compare, in this work, our approach to the \textit{pathwise} ones.   

Our main interest here consists in the \textit{functional} approach. The  \textit{functional extensions} are rooted in the extension of 
 Itô integral \textit{wrt} Brownian motion to anticipative integrands built by Hitsuda in \cite{Hi72,Hi78} and Skorohod in \cite{SK75}. 
  In \cite{GaTr82} it was proved that the stochastic integral 
\textit{wrt} Brownian motion and the adjoint of the 
derivative operator, on the Wiener space, coincide. This 
result led to many developments in (anticipative) stochastic calculus  
with respect to Gaussian processes, the most significant 
of which is \cite{nualart}. This latter article provides, using  Malliavin 
calculus, not only a divergence type integral with respect to continuous 
Volterra processes but also
Itô formulas. 
In fact, all the \textit{functional extensions} of Itô integral  developed to build a stochastic integral \textit{wrt} Gaussian processes so far
have been developed using the divergence type integral. One can divide these \textit{functional extensions} into two groups, depending if the 
assumptions are made on the kernel $K$ (first group) or on the 
covariance function $R$ (second group). 
The first group is composed of \cite{nualart} and \cite{MV05}, while the 
second one is composed of \cite{NuTa06,KRT07,KrRu10,LN12,SoVi14}. The stochastic calculus we 
propose in this work belongs to the first group since the set of assumptions 
we make is about the kernel $g$; however it does not use the divergence type integral.
The stochastic calculus we provide here allows us  to 
develop a White Noise Theory-based anticipative stochastic calculus with respect to all Gaussian processes that have an integral representation over a real (maybe infinite) interval. As stated in the beginning of this section, this class of Gaussian processes is very rich. Moreover, the stochastic calculus developed in the present work also allows us to get, not only Itô formulas but also Tanaka formulas, as well as occupation times formulas for local times of any $G$ in $\mathscr{G}$. While such results seem to be out of range for most of the intrinsic methods mentioned above, they are easily obtained using the White Noise Theory-based anticipative stochastic calculus we present here (note that all the results on Gaussian local times processes obtained using the present work are presented in \cite{JL17-2}). 
\textit{\bfseries Outline of the paper}

The remaining of this paper is organized 
as follows.  In Section 
\ref{qdokpsdksqopdqskopqsdkopdqs}, 
we recall some basic
facts about white noise theory and about the family of 
operators ${(M_{H})}_{H\in(0,1)}$, which is instrumental for 
our running example, which is presented at the end of the section. In Section 3 we define
the stochastic integral  \textit{wrt} any $G$ in $\mathscr{G}$. An Itô formula in $L^{2}(\Omega)$ is established  in the first part of Section \ref{Itô}. 
A complete comparison of our Itô formula
with all the Itô formulas for Gaussian processes, provided so far in the 
literature of \textit{functional extensions} of itô integral 
ends this section.
 In Section \ref{ozfjozrio}, 
 we compare our stochastic integral with respect to 
 elements of  $\mathscr{G}$, 
 to the divergence type 
 integrals, provided in \cite{nualart,MV05} and to Itô integral. 
In particular, we show therein
  how our integral fully generalizes the one built in \cite{MV05}.

\section{
Background on White noise theory \& on operators ${(M_{H})}_{H\in(0,1)}$}
\label{qdokpsdksqopdqskopqsdkopdqs}
  Introduced by T. \hspace{-0.1cm}Hida in \cite{Hi75}, White Noise Theory is, roughly speaking, the stochastic analogous of deterministic generalized functions (also known as tempered distributions). The idea is to realize nonlinear functional on a Hilbert space as functions of white noise (which is defined as being the time derivative of Brownian motion). 
White Noise theory has now  many application fields, such as 
quantum dynamics, quantum field theory, molecular biology, 
mathematical finance(\textit{e.g.} \cite{SCJLJL}), among many others (see \cite{Hi04} for 
more details).
One can find very good introductions (and more!) to White Noise Theory in 
\cite{HKPS,Kuo2,Si12} (see also references therein). One may also refer to 
\cite{HOUZ} for the study, in the white noise theory's framework, of 
stochastic differential equations as well as stochastic partial differential equations. 
%
  We recall in this section the standard set-up for classical white-noise theory. Readers interested in more details about White Noise Theory may refer  to \cite{HKPS,Kuo2} and \cite{Si12}.
%
%
\subsection{The spaces of stochastic test functions and stochastic distributions}
\label{cddcwcdsd}

Define $\N$ (resp. $\N^*$) the set of non negative integers (resp. positive integers). Let  $\sS(\R)$ be the Schwartz space endowed with its usual topology (\textit{i.e.} a family of functions ${(f_n)}_{n \in \N}$ 
of $\sS(\R)^{\N}$ is said to converge to $0$ if for all $(p,q)$ in ${\N}^2$ we have $\lim\limits_{n \to +\infty}  \ \sup \{\ | x^p \  f^{(q)}_{n}(x)|; \ x \in \R\} = 0$).
%
%
%
 Denote
$\sS'(\R)$ the space of tempered distributions, which is 
the dual space of $\sS(\R)$, and $\widehat{F}$ or $\cF(F)$ the 
Fourier transform of any element $F$ of $\sS'(\R)$. For every positive 
real $p$, denote $L^{p}(\R)$ the set of measurable functions $f$ such 
that $\int_{\R} {|f(u)|}^{p} \ du<+\infty$. When $f$ belongs to $L^1(\R)$, $
\widehat{f}$ is defined on $\R$ by setting $
\widehat{f}(\xi) := \int_{\R} e^{-ix\xi}  f(x) \ dx$.
Define the measurable space $(\Omega,\cF)$ by setting $\Omega:= {\sS'}(\R)$ and 
$\cF := \cB({\sS'}(\R))$, where $\cB$ denotes the $\sigma$-algebra of Borel 
sets.
The  Bochner-Minlos theorem ensures that there exists a unique probability measure $\mu$ on $(\Omega,\cF)$ such that, 
for every $f$ in $\sS(\R)$, the map $<.,f>:(\Omega,\cF) \rightarrow \R$ defined by $<.,f>(\omega) = <\omega,f>$ (where 
$<\omega,f>$ is by definition   $\omega(f)$, \ie the action of $\omega$ on $f$) is a centred 
Gaussian random variable with variance equal to ${\|f\|}^2_{L^2(\R)}$ under $\mu$. The map $f \mapsto <.,f>$ being an isometry from 
$(\sS(\R),{< , >}_{L^2(\R)})$ to $(L^2(\Omega,\cF,\mu), {< , >}
_{L^2(\Omega,\cF,\mu)})$, it may be extended to $L^2(\R)$. One may thus consider the centred Gaussian random variable $<.,f>$, for any $f$ in $L^2(\R)$.
In particular, let $t$ be in $\R$, the indicator function ${\i1}_{[0,t]}$ is defined by setting: $\i1_{[0,t]}(s):=1$ if $0 \leq s \leq t$, $\i1_{[0,t]}(s):=-1$ if  $t \leq s \leq 0$ and  $ \i1_{[0,t]}(s) :=0$ otherwise.
%
%
%
%
Then the process ${({\widetilde{B}}_t)}_{t \in \R}$, where ${\widetilde{B}}_t (\omega):= {\widetilde{B}}(t,\omega):= {<\omega, \i1_{[0,t]}>}$ is a standard Brownian motion with respect to $\mu$. It then admits a continuous version which will be denoted $B$.
Define, for $f$ in $L^2(\R)$, $I_1(f)(\omega) := {<\omega, f>} $. Then $I_1(f)(\omega)  =\int_{\R} f(s) \ dB_s (\omega) \hspace{0.2cm}  \mu-{\text {a.s.}}$, where  $\int_{\R} f(s)\ dB_s$ denotes the Wiener integral of $f$.
For every $n$ in $\N$, let $e_n(x):=   {(-1)}^n \  {\pi}^{-1/4} {(2^n n!)}^{-1/2} e^{x^2/2} \frac{d^n}{dx^n}( e^{-x^2})$ be the $n$-th 
Hermite function. It is well known (see $\cite{Tha}$) that ${(e_k)}_{k\in \N}$ is a family of functions of $\sS(\R)$ that forms an 
orthonormal basis of $L^2(\R,dt)$.
The following properties about the Hermite functions (the proof of which can be found in \cite{Tha}) 
will be useful.
\begin{theo}
\label{ozdicjdoisoijosidjcosjcsodijqpzeoejcvenvdsdlsiocfuvfsosfd}
There exist positive constants $C$ and $\gamma$ such that,  for every $k$ in $\N$,
\begin{equation*}
 {|e_k(x)|  } \leq C \ \big( {(k+1)}^{-1/12}  \cdot {\i1}_{\{|x| \leq 2\sqrt{k+1}\}} + e^{-\gamma x^2}\cdot {\i1}_{\{|x| > 2\sqrt{k+1}\}}\big).
\end{equation*}
%
%
\end{theo}

 Let  ${({| \ |}_p)}_{p\in\zZ}$ be the family  norms defined by ${|f|}
 ^2_p:=   \sum^{+\infty}_{k=0} {(2k + 2)}^{2p} \  {<f,e_k>}^2_{L^2(\R)}
 $, for all $(p,f)$ in $\zZ \times L^2(\R)$. The operator $A$, defined on $\sS(\R)$, by setting $A:= -\frac{d^2}{dx^2} + 
x^2 +1$, admits the sequence ${(e_n)}_{n \in \N}$ as eigenfunctions 
and the sequence $({{2n+2}})_{n \in \N}$ as  eigenvalues. 
Define, for $p$ in $\N$, the spaces $\sS_p(\R):=\{f \in L^2(\R), \  {|f|}_{p} <+\infty \}$ and  $
\sS_{\hspace{-0.15cm}-p}(\R)$ as being the completion of $L^2(\R)$ with 
respect to the norm ${{|\  \  |}_{-p}}$.
%
%
We summarize here the minimum background on White Noise Theory, written \textit{e.g.} in \cite[p. 692-693]{LLVH}. More precisely, let $(L^2)$ 
denote the space $L^2(\Omega,\cG,\mu)$, where $\cG$ is the $\sigma$-
field generated by ${(<.,f>)}_{f \in L^2(\R)}$. According to  Wiener-Itô's theorem, for every random variable $
\Phi$ in $(L^2)$ there exists a 
unique sequence ${(f_n)}_{n \in \N}$ of functions in ${\widehat{L}}
^2(\R^n)$ such that $\Phi$ can be decomposed as $\Phi =  {\sum^{+ \infty}
_{n = 0} I_n(f_n)}$, where ${\widehat{L}}^2(\R^n)$ denotes the set of all 
symmetric functions $f$ in $L^2(\R^n)$ and $I_n(f)$ denotes the $n-$th 
multiple Wiener-Itô integral of $f$ with the convention that $I_0(f_0) = f_0$ 
for constants $f_0$.  
For any $\Phi:=  {\sum^{+ 
\infty}_{n = 0}  \ I_n(f_n})$ satisfying the condition ${\sum^{+ \infty}_{n = 0}  
n!\ {|A^{\otimes n}f_n|}^2_{0} }< +\infty$, define the element $\Gamma(A)
(\Phi)$ of $(L^2)$ by $\Gamma(A)(\Phi):= {\sum^{+ \infty}_{n = 0} \ 
I_n(A^{\otimes n} f_n)}$, where $A^{\otimes n}$ denotes the $n-$th tensor 
power of the operator $A$ (see \cite[Appendix E]{Jan97} for more details 
about tensor products of operators).
The operator $\Gamma(A)$ is densely defined on $(L^2)$. It is invertible 
and its inverse  ${\Gamma(A)}^{-1}$ is bounded.
We note, for $\varphi$ in $(L^2)$, ${\|\varphi\|}^2_0:={\|\varphi\|}^2_{(L^2)}$. 
For $n$ in $\N$, let $\bD\text{om}({\Gamma(A)}^n)$ be the domain of the  
$n$-th iteration of $\Gamma(A)$. Define the  family  of norms  ${({\|\ \|}_p)}
_{p \in \zZ}$ by:
\begin{equation*}
 {\|\Phi\|}_p :=  {\|\Gamma(A)^p\Phi\|}_{0} =  {\| \Gamma(A)^p(\Phi)\|}_{(L^2)},  \hspace{1cm} \forall p \in \zZ,\hspace{0.5cm} \forall\Phi \in (L^2)\cap \bD{\text{om}}({\Gamma(A)}^p).
\end{equation*}
For $p$ in $\N$, define $({\cS}_{p}):=\{\Phi \in (L^2): \   \Gamma(A)^p(\Phi) \  \text{exists and belongs to}  \  (L^2) \}$ and define $({\cS}_{-p})$ as the completion of the space  $(L^2)$ with respect to the norm ${{\|\ \|}_{-p}}$.
As in \cite{Kuo2}, we let $(\cS)$ denote the projective limit of the sequence ${( (\cS_{p}))}_{p \in \N}$ and ${(\cS)}^*$ the inductive limit of the sequence  ${(({\cS_{-p}}))}_{p \in \N}$. This means in particular that ${(\cS)} \subset (L)^{2} \subset {(\cS)}^*$ and that ${(\cS)}^*$ is the dual space of $(\cS)$. Moreover, 
$(\cS)$ is called the space of stochastic test functions while ${(\cS)}^*$ the Hida distribution space. We will note  $<\hspace{-0.2cm}<\hspace{-0.1cm} \ ,\  \hspace{-0.1cm}>\hspace{-0.2cm}>$  the 
duality bracket between ${(\cS)}^*$ and $(\cS)$. If $\phi, \Phi$ belong to $(L^2)$, then we have the equality $<\hspace{-0.2cm}<\hspace{-0.1cm} \ \Phi , \varphi\  \hspace{-0.1cm}>\hspace{-0.2cm}> = {<\Phi,\varphi>}_{(L^2)} = \E[\Phi \  \varphi]$. Besides, denote $< , >$  the duality bracket between $\sS'(\R)$ and $\sS(\R)$ and recall that every tempered distribution $F$ can be written as $F = {\sum^{+ \infty}_{n = 0} \ <F,e_n>}\ e_n$, where the convergence holds in $\sS'(\R)$. The next proposition, that will be used extensively in the sequel, is a consequence of the definition of $(\cS)$ and  ${(\cS)}^*$.
\begin{prop}
\label{laiusdhv}
Let $F$ be in in $\sS'(\R)$. Define  $<.,F>:={\sum^{+ \infty}_{n = 0} \  
<F,e_n>\  <.,e_n>}$. Then there exists $p_0$ in $\N$ such that  that $<.,F>
$ belongs to  $({\cS_{-p_0}})$,  and hence to ${(\cS)}^*$. Moreover we have 
${\|<.,F>\|}^2_{-p_0}  = {|F|}^2_{-p_0}$.
    Conversely, define $\Phi:={\sum^{+ \infty}_{n = 0} \ b_n <.,e_n>}$, where $
{(b_n)}_{n \in \N}$ belongs to ${\R}^{\N}$. Then $\Phi$ belongs to ${(\cS)}^*
$ if and only if there exists an integer $p_0$ in $\N$ such that ${\sum^{+ 
\infty}_{n = 0}  \ b^2_n \ {(2n+2)}^{-2p_0}}  < +\infty$. In this latter case $F:= 
{\sum^{+ \infty}_{n = 0} \ b_n e_n}$   belongs to  $\sS_{\hspace{-0.15cm}-
p_0}(\R)$ and then to $\sS'(\R)$. It moreover verifies the equality ${|F|}^2_{-p_0} = 
{\sum^{+ \infty}_{n = 0}  \ {b^2_n} {(2n+2)}^{-2p_0} } = {\|\Phi\|}^2_{-p_0} $.
\end{prop}

\subsection{\oS-process, \oS-derivative and \oS-integral}

Let $(\R,\cB(\R),m)$ be a sigma-finite measure space. Through this section, $I$ denotes an element of $\cB(\R)$.
 A measurable function $\Phi:I \rightarrow $\oS$ $ is called a stochastic distribution process, or an \oS-process. An  \oS-process $\Phi$ is 
said to be differentiable at $t_0 \in I$ if {$\lim\limits_{r 
\to 0} \ r^{-1}\ (\Phi_{t_0+r}-\Phi_{t_0})$} exists in \oS. We 
note $\frac{d\Phi_{t_{0}}}{dt}$ the \oS\hspace{-0.1cm}-
derivative at $t_0$ of the stochastic distribution process $
\Phi$. $\Phi$ is said to be differentiable over $I$ if it is 
differentiable at every $t_0$ of $I$.
%
It is also possible to define an \oS-valued integral in the following way (one may refer to \cite[p.$245$-$246$]{Kuo2}  or \cite[Def. $3.7.1$ p.$77$]{HP} for more details). 
\begin{theodef}[{\bfseries integral in  \oS}]
\label{izeufheriduscheirucheridu}
Assume that $\Phi:I \rightarrow \ooS$ is weakly in $L^{1}(I,m)$, \ie assume that for all $\varphi$ in  ${(\cS)}$, the 
mapping $u \mapsto\ <\hspace{-0.2cm}<\hspace{-0.1cm} \  \Phi_{u},\ 
\varphi \hspace{-0.1cm}>\hspace{-0.2cm}>$,  from $I$ to $\R$, belongs to  
$L^{1}(I,m)$. Then there exists an unique element in \oS, noted $\int_{I} 
\Phi_{u} \ m(du)$, 
such that,  for all  $\varphi$ in $(\cS)$,
 \vspace{-1ex}
\begin{equation*}
{<\hspace{-0.2cm}<\hspace{-0.1cm} \ \int_{I} \Phi(u) \ m(du), \varphi \ \hspace{-0.1cm}>\hspace{-0.2cm}>} = 
\int_{I} <\hspace{-0.2cm}<\hspace{-0.1cm} \ \Phi_{u},
\varphi\  \hspace{-0.1cm}>\hspace{-0.2cm}> \ m(du). 
\end{equation*}
 \end{theodef}

%

We say in this case that $\Phi$ is \oS-integrable on $I$ (with respect to the measure $m$), in the {\it Pettis sense}. In the sequel, when we do not specify a name for the integral (\textit{resp.} for the measure $m$) of an \oS-integrable process $\Phi$ on $I$, we always refer to the integral in Pettis' sense (\textit{resp.} to the Lebesgue measure).

\subsection{S-transform and Wick product}
\label{ksksksks}
For $f$ in $L^2(\R)$, define the \textit{Wick exponential} of $<.,f>$, noted $:e^{<.,f>}:$, as the $(L^2)$ random variable equal to
$e^{<.,f> - \frac{1}{2}{|f|}^2_0}$. 
The $S$-transform of an element $\Phi$ of $(\cS^*)$, noted $S(\Phi)$, is defined as the function from  $\sS(\R)$ to $\R$ given by $S(\Phi)(\eta):= {\begin{g}{\Phi}{\ \hspace{-0.35cm}:e^{<.,\eta>}:\hspace{0.075cm}}  \end{g} }$ for any $\eta$ in $\sS(\R)$.
For any $(\Phi,\Psi)\in \ooS\times \ooS$, there exists a unique element of \oS, called the Wick product of  $\Phi$ and $\Psi$,  and  noted $\Phi \diamond \Psi$,  such that $S(\Phi\diamond \Psi)(\eta) = S(\Phi)(\eta) \  S(\Psi)(\eta)$ for every $\eta$ in $\sS(\R)$.  
Note that, when $\Phi $ belongs to $(L^2)$, $S\Phi(\eta)$ is nothing but $\E[\Phi :e^{<.,\eta>}: ] =  e^{ -\frac{1}{2}{|\eta|}^2_0} \  \E[\Phi \  e^{<.,\eta>} ]$. The following result
will be intensively used in the sequel.
\vspace{-1ex}
\begin{lem}{\cite[Lemma 2.3.]{JLJLV1}}
\label{dede}
For any $(p,q)$ in $\N^2$ and $(X,Y)$ in $({\cS}_{-p}) \times ({\cS}_{-q})$,
\vspace{-1ex}
\begin{equation*}
|S(X \diamond Y)(\eta)| \leq {\|X\|}_{-p} \  {\|Y\|}_{-q} \  e^{{|\eta|}^2_{\max \{p;q\}}}. 
\end{equation*}

\end{lem}
%
%




Some useful properties of S transforms are listed in the proposition below. The proof of the results stated in this proposition can be found in \cite[Chap $5$]{Kuo2}.

\begin{prop}[Some properties of S transforms]
\label{qmqmqmmmmm}
When $\Phi$ is deterministic then $\Phi \diamond 
\Psi = \Phi \ \Psi$, for all $\Psi$ in \oS.  Moreover, 
 let $\Phi=\sum^{+\infty}_{k = 0}  a_k\hspace{-0.1cm} <\hspace{-0.1cm}.,\ \hspace{-0.1cm}e_k>$ and $\Psi=\sum^{+\infty}_{n = 0}  I_{n}(f_{n})$ be in  \hspace{-0.1cm} \oS. \hspace{-0.2cm} Then their S-transform is given, for every $\eta$ in $\sS(\R)$, by $S(\Phi)(\eta) = \sum^{+\infty}_{k = 0} a_k\ {<\eta,e_k>}_{L^2(\R)}$\  and $S(\Psi)(\eta) = \sum^{+\infty}_{k = 0} \ <f_{n},\eta^{\otimes n}>$.
Finally, for every $(f,\eta,\xi)$ in $L^{2}(\R) \times \sS(\R)\times \R$, we have the equality:
\begin{equation}
\label{oeij12}
S(e^{i\xi <.,f>})(\eta) = e^{\frac{1}{2}({|\eta|}^{2}_{0} + 2 i \xi <f,\eta> -{\xi}^{2} {|f|}^{2}_{0})}.
\end{equation}
\end{prop}

One may refer to \cite[Chap.$3$ and $16$]{Jan97}  for more details about Wick product. The following results on the S-transform will be used extensively in the sequel. See \cite[p.$39$]{Kuo2}  and \cite[p.280-281]{HKPS}  for proofs. Denote $\cF(A;B)$ 
the set of $B$-valued functions defined on $A$.

\begin{lem}
\label{dkdskcsdckksdksdmksdmlkskdm} 
The $S$-transform verifies the following properties:
\bit
\iti The map $S:\Phi\mapsto S(\Phi)$, from \oS into $\cF(\sS(\R);\R)$, is injective.
\itii Let $\Phi:I\ \rightarrow \ooS$ be an \oS process. If  $\Phi$ is \oS-integrable over $I$ \textit{wrt} $m$, then one has,  for all $\eta$ in $\sS(\R)$,
 $S(\int_{I} \Phi(u) \  m(du) )(\eta) = \int_{I} S(\Phi(u)) (\eta) \  m(du)$.
\itiii  Let $\Phi:I \rightarrow \ooS$  be an \oS-process differentiable at $t\in I$. Then, for every $\eta$ in $\sS(\R)$ the map $u\mapsto [S \Phi(u)](\eta)$ is differentiable at $t$ and verifies $\displaystyle{S[\tfrac{d\Phi}{dt}(t)](\eta) = \tfrac{d}{dt}\big[ S[ \Phi(t)](\eta) \big]}$.
\eit
\end{lem}
The next theorems provide a criterion for integrability in \oS,  in term of $S$-transform. 

\begin{theo}{\cite[Theorem 13.5]{Kuo2}}
\label{peodcpdsokcpodfckposkcdpqkoq}
Let $\Phi:I \rightarrow {(\cS)}^*$ be a stochastic 
distribution such that, for all $\eta$ in $\sS(\R)$, the real-
valued map $t\mapsto S[\Phi(t)](\eta)$ is measurable and such
that there exist a natural integer $p$, a real $a$ and a  
function $L$ in $L^1(I,m)$ such that $ |S(\Phi(t))(\eta)| 
\leq L(t) \  e^{a {|\eta |}^2_{p} }$, for all $\eta$ of $\sS(\R)
$ and for almost every $t$ of $I$. Then $\Phi$ is \oS-
integrable over $I$, \textit{wrt} to $m$.
\end{theo}


We end this section with the following theorem that will be useful in the next section.

\begin{theo}{\cite[Theorem $2.17$]{ben1}} 
\label{lmsdlmsddmlslmsdlmsdl} 
For any  differentiable map $F:I\rightarrow {\sS'}(\R)$,  the element $<\hspace{-0.2cm}.,F(t)\hspace{-0.2cm}>$ is a differentiable stochastic distribution process which satisfies the equality:  

\cl{$\frac{d}{dt} <.,F(t)> \ = \  <.,\frac{d F}{dt} (t)>$.}
\end{theo}

%


\vspace{-2ex}
\subsection*{Gaussian Processes in $\mathscr{G}$ of ``reference''}
\vspace{-1ex}
To see in what extent the stochastic calculus \textit{wrt} Gaussian processes we present here generalizes the one 
provided in the 
literature so far, 
we will consider, throughout this paper, a running  example, made with elements of $\mathscr{G}$ that are Brownian motion and Brownian 
bridge, fractional and multifractional Brownian motions as 
well as $\mathscr{V}_{\gamma}$ - processes (the last three processes being defined below).  
\smallskip

\textit{\bfseries Fractional and Multifractional Brownian motions}
\smallskip

Readers interested in an exhaustive presentation of fBm or mBm may refer to \cite{Nu} for fBm and to \cite{LLVH} for mBm, as well as to the references therein. Introduced in \cite{Kol} and popularized in \cite{Mandelbrot1968}, fBm is a centered Gaussian process, the covariance function of which is denoted $R_{H}$ and is given by:
\vspace{-1ex}
\begin{equation*}
  R_{H}(t,s) :=  \frac{1}{2} ( {|t|}^{2{H}} + {|s|}^{2{H}} - {|t-s|}^{2{H}}),
\end{equation*}
where $H$ belongs to 
$(0,1)$, and is usually called the Hurst exponent. When $H=1/2$, fBm reduces to standard 
Brownian motion. Among many other properties, fBm is able to match any prescribed constant local regularity and to model phenomena that presents long range dependence. These properties made this process very popular in many fields such as mathematical finance, Internet traffic 
modeling, image analysis and synthesis, physics and more.

MBm, which is a Gaussian extension of fBm, was 
introduced in \cite{PL} and in \cite{ABSJDR} in order to 
match any prescribed non-constant deterministic local 
regularity and to decouple this property from long range 
dependence (this impossibility of doing so for fBm constitutes one of the most severe drawbacks of this process). To obtain mBm, 
the idea is to replace the constant Hurst parameter $H$ of 
fBm by a deterministic function $t\mapsto h(t)$ ranging in $
(0,1)$. Several definitions of mBm exist and the reader 
interested in the evolution of these definitions may refer to 
\cite{LLVH} and references therein. We will only give here 
the definition 
of mBm given in \cite[Definition 1.2]{LLVH}, which is not only the most 
recent but also includes all previously known ones. A mBm 
on $\R$, with functional parameter $h:\R \rightarrow (0,1)$, 
is a Gaussian process $B^h:={(B^h_t)}_{t\in \R}$ defined, 
for all real $t$, by $B^h_t:=\sifbm(t,h(t))$, where  $
\sifbm:={(\sifbm  (t,H))}_{(t,H) \in \R\times(0,1)}$ is 
fractional Brownian field on $\R\times(0,1)$ (which means 
that  $\sifbm$ is a Gaussian field, such that, for every $H$ 
in $(0,1)$, the process ${(\sifbm  (t,H))}_{t \in \R}$ is a fBm 
with Hurst parameter $H$). In other words, a mBm is 
simply a ``path'' traced on a fractional Brownian field. 
Note also that when $h$ is constant, mBm reduces to fBm. 
 The literature on Stochastic integration \textit{wrt} fBm is extensive now. The reader interested in an exhaustive overview of the subject may refer to \cite{Nu,COU07,Mis08} for divergence type integral and to \cite{ben1,bosw,ell,Nu05} for integral in the white noise theory framework.
\noindent More recent, the literature on Stochastic integration \textit{wrt} mBm is less rich. Nevertheless, one may cite \cite{dozzi} for a divergence type integral \textit{wrt} to a Volterra-type mBm and \cite{JLJLV1,JL13} for a Wick-Itô multifractional integral (\ie an integral \textit{wrt} to normalized mBm, in the White Noise theory framework). 
Note moreover that \cite{LLVH} provides a general method of integration \textit{wrt} to all classes of mBm, that does not only apply for divergence type integral and white noise theory integral but also for pathwise integral. 
 
\textit{\bfseries $\mathscr{V}_{\gamma}$ - processes}

\cite{MV05} provides a stochastic calculus, \textit{wrt} a particular class of 
Volterra processes, that we will denote $
\mathscr{V}_{\gamma}$ - processes in the sequel.
For any deterministic function $\gamma:\R_{+}\rightarrow \R$, $
\mathscr{V}_{\gamma}$ - processes 
are defined in \cite[Proposition 1]{MV05} as being the processes, denoted ${\widetilde{B}}^{\gamma}:= {({\widetilde{B}}^{\gamma}_{t})}_{t\in [0,T]}$, by setting:
\vspace{-2ex}
\begin{equation}
\label{erofkorekdeddeddeezzeep}
{\widetilde{B}}^{\gamma}_{t}:=\int^{t}_{0} \varepsilon(t-s) \ dW_{s}; \hspace{0.5cm} \forall t\in [0,T],
\end{equation}
with $\gamma:\R_{+}\rightarrow \R$ such that $
\gamma^{2}$ is of class $C^{2}$ everywhere in $\R_{+}$ 
except in $0$; and such that ${(\gamma^{2}})'$ is non increasing. The map $\varepsilon:\R^{*}_{+}\rightarrow \R$ is 
defined by setting $\varepsilon:= \sqrt{{(\gamma^{2}})'}$. Subset of $\mathscr{G}$, the set 
$\mathscr{V}_{\gamma}$ contains Gaussian processes, that can be 
more irregular than any fBm. However it does not contain fBm (nor mBm) since $\mathscr{V}_{\gamma}$ only contains processes the regularity of which
 remains constant along the time). It will be shown in Remark \ref{peorfkerpofkeprofkepr}
 that the stochastic integral built in 
\cite{MV05} is a 
particular case of the stochastic 
integral we build in this work. 
For notational simplicity we will refer to these processes as the Gaussian 
processes of ``reference''.  
\subsection[Family of operators ${(M_{H})}_{H\in(0,1)}$ and a classical example of Gaussian processes in $\mathscr{G}$]{Operators ${(M_{H})}_{H\in(0,1)}$ and a classical set of Gaussian processes in $\mathscr{G}$}
\label{eoajeorijzp}
The operator $M_{H}$ will be useful in the sequel, not only to provide one with a representation of fBm and of mBm under the form \eqref{erofkorekp}, but also to verify that Assumptions \Aa
, we will make in Section \ref{sjjjj}, hold for both fBm and mBm. 
Let $H$ belongs to $(0,1)$; following \cite{ell} and \cite[Section $2.2$]{JLJLV1}, define the $L^{2}(\R)$-valued operator $M_H$, in the Fourier domain by:
\begin{equation*}
\widehat{M_H(u)}(y) := \tfrac{\sqrt{2\pi}}{c_H} \hspace{0.1cm} |y|^{1/2-H} \ \widehat{u}(y), \hspace{0.5cm} \forall y \in \R^*,
\end{equation*}
where $c_{x}$ is defined, for every $x$ in $(0,1)$ by $c_x:=  {\big(\tfrac{2\pi}{\Gamma(2x+1) \sin(\pi x)} \big)}^{\frac{1}{2}}$.
This operator is well defined on the homogeneous Sobolev space
$ L^2_H(\R):= \{u \in {\sS'}(\R) \hspace{0.1cm}:\hspace{0.1cm} \widehat{u} = T_f; 
\hspace{0.1cm} f \in L^{1}_{loc}(\R) \hspace{0.1cm}\textnormal{and}
\hspace{0.1cm} \| u \|_H < +\infty \},$
%
%
%
where the norm $\ds \|\cdot\|_H$ derives from the inner product denoted ${\langle \cdot , \cdot \rangle}_H$, which is defined on $\ds L^2_H(\R)$ by: 
\begin{equation*}
 {\langle u,v \rangle}_H := \frac{1}{c^2_H} \int_\R |\xi|^{1-2H} 
{\widehat{u \hspace{0.1cm}}(\xi) }\ \overline{{\widehat{v \ }
(\xi) }} \ d{\xi}. 
\end{equation*}

$M_{H}$ being an isometry from $(L^2_H(\R), \ds \|\cdot\|_H)$ 
into $(L^2(\R), \ds \|\cdot\|_{L^2(\R)})$, it is clear that, for every $(H,t,s)$ in $(0,1)\times \R^{2}$, ${<M_{H}(\i1_{[0,t]}), M_{H}(\i1_{[0,s]})>}_{L^{2}(\R)}= R_{H}(t,s)$.
%
%
%
We will say that an mBm is normalized when its covariance function, denoted $R_{h}$, verifies the equality:
\vspace{-2ex}

\begin{equation}
\label{azazazaz}
R_{h}(t,s)=     \tfrac{ c^2_{h_{t,s}}}{c_{(h(t))}c_{(h(s))}}  \ \big[\tfrac{1}{2} \big( {|t|}^{2h_{t,s}} + {|s|}^{2h_{t,s}}  - {|t-s|}^{2h_{t,s}} \big)\big],
\end{equation}
where $h_{t,s} := \frac{h(t)+h(s)}{2}$ and $c_{x}$ has been above, right after $\widehat{M_H(u)}(y)$.

\begin{ex}[Gaussian Processes in $\mathscr{G}$ of ``reference'']
\label{oevieroivjedzezdzeddo}
Let $H$ be real in $(0,1)$ and $h:\R \rightarrow (0,1)$ be a deterministic measurable function. Define the processes 
\vspace{-0.5cm}

\begin{align*}
&B:=\{   <.,\i1_{[0,t]}>;\ t\in\R   \};& \
&\widehat{B}:=\{<.,\i1_{[0,t]} - t\cdot \i1_{[0,1]}>;\ t\in[0,1]   \};&\\
&B^{H}:=\{ <.,M_{H}(\i1_{[0,t]})>;\ t\in\R   \};& \
& B^{h}:=\{ <.,M_{h(t)}(\i1_{[0,t]})>;\ t\in\R   \};\ & \\
&{\widetilde{B}}^{\gamma}:=\{ <.,\i1_{[0,t)} \cdot \varepsilon(t-.))> \text{ if }  t\in\R^{*}_{+}  \ \& \  {\widetilde{B}}^{\gamma}_{0}:= 0 \}.& \
\end{align*}
%
We know, thanks to Section \ref{cddcwcdsd}, that $B$ is a 
Brownian motion on $\R$. Moreover, since for any $g_{t}$ in $L^{2}(\R)$, $<.,g_{t}>  \stackrel{\text{a.s.}}
{=} \int_{\R} \ g_{t}(u)\ dB_{u}$, it is clear, in view of the definition of ${\langle \cdot , \cdot \rangle}_H$, that ${B}^{H}$ is a fBm of Hurst index $H$, 
that  $B^{h}$ is a normalized mBm of functional parameter 
$h$, that $\widehat{B}$ is Brownian bridge on $[0,1]$ and 
that ${\widetilde{B}}^{\gamma}$ is a $\mathscr{V}_{\gamma}$ - process (defined in \eqref{erofkorekdeddeddeezzeep}).
\end{ex}

\vspace{-0.1cm}

A word on notation: $B^H_.$ or $B^{h(t)}_.$ will always 
denote an fBm with Hurst index $H$ or $h(t)$, while 
$B^h_.$ will stand for an mBm. 
Throughout this paper, unless otherwise specify, we will neither specify the 
value of $H$ in $(0,1)$, when we consider a fBm, nor the $(0,1)$-valued 
function $h$ when we consider a mBm, nor the function $\gamma$ of a $\mathscr{V}_{\gamma}$ - process. 
\section{Stochastic integral with respect to Gaussian process}
\label{sjjjj}
The first part of this section is devoted to the definition of 
the time derivative, in the Stochastic 
distribution sense, of any element $G:=(G_t)_{t\in \rR}$  
of $\mathscr{G}$.
 We then compute the $S$-transforms of processes $G$ and of its time derivative.
 The Wiener integral  \textit{wrt}  $G$ is 
presented in Subsection \ref{depzokd}, whereas the 
stochastic integral  \textit{wrt}  $G$ is built in 
Subsection \ref{ozieffjioefjoezifjioezzjioefjfeoiz}. 
We keep the same notations as in Section \ref{qdokpsdksqopdqskopqsdkopdqs}. In particular, the probability space $(\Omega,\cF,\mu)$, described in the previous section is now fixed.
Denote $G:=(G_t)_{t\in \rR}$ the process defined, for every $t$ in $\rR$, 
by $G_{t} :=  <.,g_{t}>$, where ${(g_t)}_{t\in\rR}$ is a family of functions of $L^2(\R)$. As we saw in Example 
\ref{oevieroivjedzezdzeddo}, $G$ is a Gaussian process which fulfills 
the equality $G_t \stackrel{{\text{a.s.}}}{=} \int_{\R} g_t(u) \ dB_u$. 
Denote $(t,s)\mapsto R_{t,s}$ the covariance function of $G$. We hence have 
$R_{t,s}:=\mathbf{E}[G_t\ G_s]={<g_{t},g_{s}>}_{L^{2}(\R)}$, for every $
(s,t)$ in $\rR^{2}$. We will note in the sequel $R_{t}$ instead of $R_{t,t}
$. For the sake of notational simplicity we assume that $G_0 
\stackrel{{\textit a.s.}}{=}0$. Moreover, when the Gaussian process $G$ 
will admit a continuous modification, we will systematically use it and still 
call it $G$.
%
\subsection{White Noise derivative of $G$}
\label{oeifjoerfjerij}
Define the map  $g:\rR\rightarrow {\sS'}(\R)$ by setting 
$g(t) := g_t$. When $g$ is differentiable at point $t$, one 
denotes $g '_t$ its derivative. Denote $\lambda$ the Lebesgue measure on $\rR$ and define $L^1_{l\widetilde{oc}}
(\rR):=\{f:\rR\rightarrow \R \text{ is measurable }; \ f\in L^1((a,b)), \hspace{0.1cm} \text{for all finite interval  } (a,b) \textit{ s.t. } [a,b] \subset \rR \}$. 
In this section and in the next one  
(namely in Sections \ref{sjjjj} and \ref{Itô}), we will make the following assumption:
\makeatletter\tagsleft@true\makeatother

\begin{minipage}{1\textwidth}
\vspace{-4ex}
  \begin{center}
\textcolor{white}{
\begin{subequations}
 \label{A}
 \begin{align}
  a & = b \tag*{\text{($\mathscr{A}$)}} \label{Aa}\\
  a & = b \tag*{$\mathscr{A}_{\text{(i)}}$} \label{Ai}\\
  a & = b \tag*{$\mathscr{A}_{\text{(ii)}}$} \label{Aii}\\
  a & = b \tag*{$\mathscr{A}_{\text{(iii)}}$}\label{Aiii} \\
  a & = b \tag*{$\mathscr{A}_{\text{(iv)}}$} \label{Aiv}
 \end{align}
\end{subequations}
}
\vspace{-21.5ex}
\begin{align}
\label{Aba}
\tag*{(\hspace{-0.1ex}$\mathscr{A}_{\textcolor{white}{\text{iiiiii}\hspace{-1.6ex})}}$\hspace{-1.7ex})} 
\hspace{4ex} \begin{cases} 
(\text{i})  \text{ The map $g$ is continuous on } \rR,\\
(\text{ii})  \text{  The map $g$ is differentiable } \lambda-\text{almost everywhere on } \rR,  \\
(\text{iii}) \text{   There exists }  q   \text{ in } \N^*\text{ such that } {t\mapsto {|g'_t|}_{-q} \text{ belongs to } L^1_{l\widetilde{oc}}(\rR),}\\
(\text{iv}) \text{   For every } (a,b) \text{ in } {\rR}^{2} \text{   such that } a \leq b, \text{ one has, in $\sS'(\R)$, the equality: }\\
\textcolor{white}{pfojerpfoerofjerof} 
\end{cases} 
\end{align}

\vspace{-8.5ex}
\begin{equation}
\label{zodicjzoscjoifjs}
\tag{$E_{a,b}$}
g_{b} - g_{a} = \int^{b}_{a} \ g'_{u} \ du.
\end{equation}
   \end{center}
\end{minipage}

\makeatletter\tagsleft@false\makeatother 
\vspace{1ex}

Proposition \ref{CS} below will provide an easy way to check whether Assumption \Aa holds or not. Besides, define the set $\rR_{D}$ by setting $\rR_{D}:=\{t\in\rR; \vspace{1ex} \ g \text{  is differentiable at point } t\}$. Of course  
$L^1_{l\widetilde{oc}}(\rR)$ contains in particular all measurable functions $f:
\rR_{D}\rightarrow \R$ such that  $f\in L^1((a,b)\cap\rR_{D})$ (that we will denote $L^1((a,b))$ in the sequel, by abuse of notation), for every finite interval  
$(a,b)$ \textit{ s.t. } $[a,b] \subset \rR$. For the sake of notational simplicity we will write  \ref{Ai}, 
(\text{resp.} \ref{Aii}, \ref{Aiii} or \ref{Aiv}), in the sequel, 
when one wants to refer to statement (i) (\text{resp.} to (ii), (iii) or (iv)) of Assumption \Aap.

Making Assumptions \ref{Ai}  and \ref{Aii}  seems reasonable since we want to \enquote{differentiate}, 
with respect to $t$, the Gaussian process $G$, which 
trajectories are, in general, not differentiable in the strong sense ({\textit{e.g.} the Brownian motion). The interest of Assumptions  \Aaiii  and \Aaiv will be explained when it will be needed (in Section \ref{depzokd}, right after Definition \ref{WienWick}). 
\smallskip

\vspace{-1ex}

\begin{rem}
\label{ozeufhioeruvheirezzezeedzedzedzdu}
${\text{{\bfseries{1.}}}}$\  A first consequence of Assumption \ref{Aa} is that $g$ is 
``weakly'' locally absolutely continuous on $\rR$; that is that the map $t\mapsto <g_{t},\eta>$ is absolutely continuous on every finite interval $[a,b]$ of $
\rR$, for every $\eta$ in $\sS(\R)$.

${\text{{\bfseries{2.}}}}$\ If $g$ would have been a real-valued function, 
Assumption \ref{Aa} would have been nothing but the  
local absolute continuity of $g$ on $\rR$. However, $g$ is $\sS'(\R)$-valued. Thus, and even if a notion of absolute continuity 
exists for $\sS'(\R)$-valued functions (see 
\cite[Definitions $3.6.2$ \& $3.2.4$]{HP}), the absolute continuity of $g$ on an interval $[a,b]$ of $\rR$ does not entail the differentiability of $g$ in general (see an example that illustrates this fact in \cite{HP}, right above Theorem $3.8.6$).

\end{rem}

\vspace{-1.5ex}

An easy way to see if Assumption \Aa holds is to check if the sufficient condition provided in the following proposition, and that will be used a lot in the sequel, holds.

\begin{prop}
 \label{CS}
 A sufficient condition for Assumption \Aa to be verified is that:
\vspace{-4.5ex}

\textcolor{white}{
\begin{subequations}
 \label{Ddedede}
 \begin{align}
  a & = b \tag*{\text{($\mathscr{D}$)}} \label{D}\\
  a & = b \tag*{$\mathscr{D}_{\text{(i)}}$} \label{Di}\\
  a & = b \tag*{$\mathscr{D}_{\text{(ii)}}$} \label{Dii}
   \end{align}
\end{subequations}
}

\vspace{-17ex}

\vspace{-2ex}

$$
\hspace{-3.5ex}(\mathscr{D})
\begin{cases} 
(\text{i})  \text{ The map $g$ is continuous on $\rR$ and differentiable on  every finite interval } (a,b) \textit{ s.t. } [a,b] \subset \rR,\\
(\text{ii}) \text{   There exists }  q   \text{ in } \N^*\text{ such that } {t\mapsto {|g'_t|}_{-q} \text{  belongs to } L^1_{l\widetilde{oc}}(\rR).} 
\end{cases} 
$$


\end{prop}

\begin{pr}
 
Indeed, these two conditions obviously entail that Assumptions \Aai to 
 \Aaiii hold. Moreover, these two conditions also entail Equality $(E_{c,d})$, for every $[c,d] \text{ in } \rR_{D}$. $\rR$ being a closed interval of $\R$, and in view of $(ii)$ of Assumption \Da, there can be $0$,$1$ or $2$ points, at the maximum, that belong to $\rR$ but not to $\rR_{D}$. Let us treat this latter case only and denote $a$ and $b$ these two points.  The continuity of $g$ at points $a$ and $b$, from one hand, and the 
Lebesgue dominated convergence theorem, from the other hand, give us the equality $(E_{a,b})$.
\end{pr}
As previously,  and for the sake of notational simplicity we will write  \Dai (\textit{resp.} \Daii ) when one wants to refer to (i) or to (ii) of Assumption \Dap. 
Assumption \Aaiip, Proposition \ref{laiusdhv} and Theorem \ref{lmsdlmsddmlslmsdlmsdl} 
lead to the following definition of Gaussian  white noise.
%
\begin{theodef}[{\bfseries Gaussian White Noise}]
\label{oijoifjsoidjdoijfsoijfsoijdsoifjdoisdjosdjf}
Define for every $t$ in $\rR_{D}$,
\begin{equation}
\label{oerijefoijefjfijvoidfovijdfovijdfo}
W^{(G)}_{t} := \ <.,g'_t>,
\end{equation}
where the equality holds in  ${(\cS)}^*$. Then $({W^{(G)}_{t})}_{t \in \rR_{D}}$ is a  ${(\cS)}^*$-process and is the ${(\cS)}^*$-derivative
of the process ${(G_t)}_{t \in \rR_{D}}$. We will sometimes note $\frac{dG_t}{dt}$ instead of $W^{(G)}_{t}$. 
\end{theodef}

 
Using Proposition \ref{laiusdhv} one easily sees that \eqref{oerijefoijefjfijvoidfovijdfovijdfo} also reads, for every $t$ in $\rR_{D}$:

\vspace{-1.5ex}

\begin{equation*}
W^{(G)}_{t}  = \sum^{+\infty}_{k = 0} \    {<g'_t,e_k>}  {<.,e_k>}  = \    \sum^{+\infty}_{k = 0} \  (\tfrac{d}{dt} <g_t,e_k>) \ <., e_k>. 
\end{equation*}
\begin{prop}
\label{oedfvdfoivjdfoivjdfoij}
 The map $t \mapsto {\| W^{(G)}_{t}\|}_{-p}$ is continuous if and only if $t \mapsto {|g'_t|}_{-p}$ is continuous. 
\end{prop}

{\bfseries Proof.}
Thanks to Proposition \ref{laiusdhv}, one can write ${\| W^{(G)}_{t}\|}_{-p} =  {|g'_t|}_{-p}$, $\forall\ (p,t)$ in $\N^*\times \rR_{D}.  \hfill  \square$

\vspace{0.2cm}

As the next example shows Assumptions \Da (and therefore Assumption \Aap) holds in the case of all Gaussian processes in $\mathscr{G}$ of  ``reference''.
Denote, for every $n$ in $(1/2,+\infty)$,
\vspace{-0.25cm}

\begin{equation}
\label{zoijeroif}
\cR_n:= \sum^{+\infty}_{k = 0} {(2k+2)}^{-2n}.
\end{equation}

\begin{ex}
\label{oevieroivjeo}
${\text{{\bfseries{1.}}}}$\  (Brownian motion on $\R$ \& Brownian bridge on $[0,1]$). For the Brownian motion on $\R$ (\textit{resp.} the Brownian bridge on $[0,1]$), one has $\rR=\rR_{D}=\R$, and, for every real $t$, $g'_t = \delta_t$  (\textit{resp.} $\rR=\rR_{D}=[0,1]$ and $g'_t = \delta_t - \i1_{[0,1]}$). Both maps $g$ clearly fulfills Asumption $\Di$. Moreover, for every $p$ in $\N^{*}$, the maps $t 
\mapsto {|g'_t|}_{-p}$ are continuous and bounded on $\R$, which shows that \Daii holds. Indeed, using both: the relation $e'_{k}(x)=\sqrt{\frac{k}{2}} e_{k-1}(x)- \sqrt{\frac{k+1}{2}} e_{k+1}(x)$, valid for all positive integer $k$, and
Theorem \ref{ozdicjdoisoijosidjcosjcsodijqpzeoejcvenvdsdlsiocfuvfsosfd}, we get the existence of a real $C'$, independent of $t$ and $p$, such that: $ \forall (p,t) \in \N^*\times \R, \hspace{2ex} {|
\delta'_t|}^2_{-p} =  \sum^{+\infty}_{k = 0} \ e^2_k(t) \  \ {(2k
 +2)}^{-2p}\leq  C' \cdot \cR_{p}.$
%
%

${\text{{\bfseries{2.}}}}$\   (Fractional \& Multifractional Brownian motions on $\R$) In both these cases, one has $\rR=\rR_{D}=\R$.
Thanks to  \cite[Remark 4.3 and Proposition 4.10]{JLJLV1}, we know that  
Asumption \Daip, as well as the fact that $t \mapsto {\| W^{(B^H)}_{t}\|}_{-
p}$ (\textit{resp. $t \mapsto {\| W^{(B^h)}_{t}\|}_{-p}$}) is continuous and 
bounded on any compact set of $\R$,  are verified for every $p\geq 2$ and $H
$ in $(0,1)$ (resp. $h$ differentiable with locally bounded derivative).

${\text{{\bfseries{3.}}}}$\   (The process ${\widetilde{B}}^{\gamma}$) In this case, $\rR=\R_{+}$ and $\rR_{D}=\R^{*}_{+}:=\R\backslash\{0\}$. Moreover Assumption \Da is also fulfiled, as  Theorem \ref{pzoefkpoezzokfzepokzpzefj1286564536} below shows.
\end{ex}

The following theorem, the proof of which can be found in Appendix 
\ref{appendice2}, shows that $\mathscr{V}_{\gamma}$ - processes  fulfill  Assumption \Da. 
This will be crucial to show, in Section \ref{zeofkepofkpo123879987}, that the 
stochastic integral \textit{wrt} $\mathscr{V}_{\gamma}$- processes, developed in \cite{MV05}, is a particular case of the one we build in this work. Let us first define the two maps $E:\R_{+} \rightarrow \R$ and $\mathscr{E}:\R_{+} \rightarrow \R$ by setting:
\vspace{-0.35ex}
\begin{equation*}
E(x):=
\begin{cases}
 \int^{x}_{0}\ \varepsilon(u) \ du, & {\text{ if}} \hspace{0.25cm}  x\in\R^{*}_{+},\\
0 &   {\text{ if}} \hspace{0.25cm}  x=0, \\
\end{cases}
\hspace{1cm} \& \hspace{0.75cm} 
\mathscr{E}(x):= \int^{x}_{0}\ E(u) \ du.
\end{equation*}
 For every $\Psi$ in ${\sS'}(\R)$, $\Psi'$ will denote the 
derivative of $\Psi$, in the sense of tempered distributions\footnote{One therefore has $ <(G_{t})',\varphi> = -<G_{t},
\varphi'> $, for every $(t,\varphi)$ in $\R^{*}_{+}\times \sS(\R)$.}.

\begin{theo}
\label{pzoefkpoezzokfzepokzpzefj1286564536}
The map $\Phi: \R_{+} \rightarrow \sS'(\R)$ defined by setting: 
\vspace{-2ex}
\begin{equation*}
\Phi(t):= \Phi_{t}:=
\begin{cases}
 \i1_{[0,t)} \cdot \varepsilon(t-\cdot)  &  {\text{ if }}  \hspace{0.25cm}  t\in\R^{*}_{+},\\
0 &   {\text{ if}} \hspace{0.25cm}  t=0, \\
\end{cases}
\end{equation*}
fulfill Assumption \Aap. More precisely, it is differentiable on $\R^{*}_{+}$ and, $\forall$  $t$ in  $\R^{*}_{+}$,
\vspace{-1ex}
\begin{equation}
 \label{erfjeorfjeorfjeroiferofije}
\Phi'(t):=\frac{d}{dt}[\Phi(t)] = F_{t} - (G_{t})'+ (H_{t})'',
\end{equation}
where $F_{t}, G_{t}$ and $H_{t}$ all belong to ${\sS'}(\R)$ 
and are defined by setting, $\forall\ t$ in $\R^{*}_{+}:$ 
\vspace{-1ex}
\bit
\item [\tiny$\bullet$] $F_{t}:= \frac{\varepsilon(t-\cdot)}{t} \i1_{[0,t)} + 
\big(\varepsilon(t) - \frac{E(t)}{t}\big)\ \delta_{0}$; \ \hfill
 {\tiny$\bullet$} $H_{t}:=  \ds \frac{((t-\cdot)\cdot E(t-\cdot) - \mathscr{E}(t-\cdot))   }{t} \cdot \i1_{[0,t)}$;
\item [\tiny$\bullet$] $G_{t}:=  \ds\big(E(t) - \tfrac{\mathscr{E}(t)}{t}\big)\ \delta_{0} + u\mapsto \left(\frac{u\ \varepsilon(t-u) - E(t-u)}{t}\right) \i1_{[0,t)}(u)$.
\eit

Furthermore, the map  $t\mapsto {|\Phi'(t)|}_{-q}$ belongs to $\underset{\ \ b\in\R^{*}_{+}}{\cap} L^2((0,b))$, for every integer $q\geq 3$.
\end{theo}

\subsection{Generalized functionals of $G$}
\label{eprofkerpokveprvkerpovepovk}
In order to establish easily 
that the map $t\mapsto f(G_t)$ is \oS-integrable and integrable with respect to itself, 
when $f$ is function of polynomial growth, 
we introduce here the  generalized 
functionals of $G$, using \cite[Section $7.1$]{Kuo2}. We 
identify, here and in the sequel, any function $f$ of $L^{1}
_{loc}(\R)$ with its associated tempered distribution, denoted noted $T_f$, 
when it exists. In particular, one notes 
in this case: $<f,\phi>  \ = \int_{\R} f(t) \  \phi(t) \  dt$,  for 
every $\phi$ in  $\sS(\R)$. In this latter case we say that 
the tempered distribution $T:=T_{f}$ is of function type. 
Define the sets $\Z_R:=\{t \in \rR; \ R_{t} = 0\}$ and $
\Z^c_R:=\{t \in \rR; \ R_{t} > 0\}$.


%

%

\begin{theodef}
\label{lmqsdlmflkjdg}
 Let $F$ be a tempered distribution. For every $t$ in $\Z^c_R$, define 
 \vspace{-1.5ex}
 
\begin{equation*}
 F(G_t):= \frac{1}{\sqrt{2\pi {R}_{t}}} \sum\limits^{+\infty}_{k=0} \frac{1}{\ k! \ 
{{R}^{k}_{t}}}  <F, {\xi}_{t,k}> I_k\left(   {g_t }^{\otimes k} \right),
\end{equation*}

 where, for every $(x,k)$ in $\R\times \N$, ${\xi}_{t,k}(x) := {\pi}^{1/4} {(k!)}
^{1/2}  {{R}^{k/2}_{t}}\  {\exp}{\{-\frac{x^2}{4 {{R}_{t}}}\}}  {e}_{k}{(x/(\sqrt{2 
{R}_{t}}))}$. Then for all real $t$, $F(G_t)$ is a Hida distribution, called {\it 
generalized functional of $G_t$}.
 \end{theodef}

\begin{pr}
This is an immediate consequence of \cite[p.$61$-$64$]{Kuo2}, by taking $f := g_t$.
\end{pr}

\begin{rem}
As shown in \cite{ben1}, when $F=f$ is of function type, $F(G_t)$ coincides with $f(G_t)$. 
\end{rem}

The following theorem yields an estimate of ${\|F(G_t)\|}^2_{-p}$ which will be useful in the sequel.

\begin{theo}
Let  $p$ be in $\N$. Then there is a constant $D_p$, 
such that:
\vspace{-1ex}
\begin{equation}
\label{spdsdpokpovkopsdkvspodvksopsdvpoksvdokdsvkopsvdopqbschfgsdvhjerhbtkljtyopfdhvlkdfhu}
\forall \ F \in {\sS}_{\hspace{-0.15cm}-p}(\R),\ \forall \ t \in \Z^c_R, \hspace{0.5cm} {\|F(G_t)\|}^2_{-p} \leq \ \max\{   {{R}^{-2p}_{t}}; \ {{R}^{2p}_{t}} \} \    {{R}^{-1/2}_{t}}   \  \ D_p  \ {|F|}^2_{-p}.
\end{equation}
\end{theo}

\begin{pr}
This is a simple consequence of the following more general result:  let $f$ be a nonzero function in $L^{2}(\R)$, $p \in \N$ and $F \in {\sS}_{\hspace{-0.15cm}-p}(\R)$. There exists a constant $D_p$, independent of $F$ and $f$, such that $ {\|F(<.,f>)\|}^2_{-p} \leq \ \max\{   {|f|}^{-4p}_{0}; \ {|f|}^{4p}_{0} \} \    {{|f|}^{-1}_{0}}   \  \ D_p  \ {|F|}^2_{-p}$.
%
The line of the proof is the same as in \cite[Theorem $3.3$]{ben1} by replacing  there ${t}^{2H}$ by  ${|f|}^{2}_0$. 
\end{pr}

\subsection{S-Transform of $G$ and $W^{(G)}$}
\label{eotivj}

The following theorem makes explicit the $S$-transforms of $G$, of the Gaussian white noise $W^{(G)}$ and of generalized functionals of $G$. Denote $\gamma$ the heat kernel density on $\R_+\times\R$ \textit{i.e.} 
\vspace{-0.25cm}

\begin{equation}
\label{fctiongammadzedezdz}
\gamma(t,x):= \tfrac{1}{\sqrt{2\pi t}}\exp{\{ \tfrac{-x^2}{2t}  \}} \text{ if } t \neq 0 \text{ and } 0 \text{ if } t = 0. 
\end{equation}

The results provided in Theorem \ref{tardileomalet} below will be used a lot in the sequel, and in Section \ref{depzokd}.

\begin{theo}
\label{tardileomalet}
For every $\eta$ in $\sS(\R)$ one has the following 
equalities:
\vspace{-0.25cm}

\bit
\iti  $S(G_t)(\eta) = {<g_t,\eta>}_{L^2(\R)}$, for every $t$ in $\rR$,
\vspace{-0.2cm}
\itii $S(W^{(G)}_{t})(\eta) = {<g'_t,\eta>} =\frac{d}{dt} [{<g_t,\eta>}_{L^2(\R)}] $, for every $t$ in  $\rR_{D}$;
\itiii  For $p \in \N$, $F \in \sS_{\hspace{-0.15cm}-p}(\R)$, and $t$ in $\Z^c_R$,  \  $S(F(G_t))(\eta) = \left< F, \gamma \left(R_{t}, .- {<g_t, \eta>}  \right) \right>$.
\smallskip

Furthermore, there exists a constant $D_p$, independent of $F,t$ and $\eta$, such that:
$$\forall t \in \Z^c_R, \hspace{1cm} \ {|S(F(G_t))(\eta)|}^2  \leq \ \max\{   {{R}^{-2p}_{t}}; \ {{R}^{2p}_{t}} \} \    {{R}^{-1/2}_{t}}   \  \ D_p  \ {|F|}^2_{-p} \  \exp\{{| \eta |}^2_{p}\}.$$

\eit
%
%
%
%
%
%
%
%
\end{theo}


\begin{pr}
$(i)$ Obvious in regard of Proposition \ref{qmqmqmmmmm}. Point $(ii)$ is a straightforward consequence of $(iii)$ in Lemma \ref{dkdskcsdckksdksdmksdmlkskdm}, and of  \eqref{oerijefoijefjfijvoidfovijdfovijdfo}. The equality in $(iii)$ results from \cite[Theorem $7.3$ p.$63$]{Kuo2} with $f = g_{t}$. The inequality results from \eqref{spdsdpokpovkopsdkvspodvksopsdvpoksvdokdsvkopsvdopqbschfgsdvhjerhbtkljtyopfdhvlkdfhu} as in \cite[Theorem $3.8$]{ben1}.
\end{pr}

%
%

Before giving the general result on stochastic integral \textit{wrt} 
$G$ we deal, in the next subsection, with Wiener integral  \textit{wrt}  $G$.


\subsection{Wiener integral with respect to $G$}
\label{depzokd}
%
%
%
%
%
\hspace{0.5cm}In all this subsection one denotes $I$ a Borel set of $\rR$ and $f:\rR \rightarrow \R$ a deterministic and measurable function
We want 
to define the integral of $f$, on $I$, with respect to $G$. Since the map $s \mapsto 
G_s$ is (weakly) differentiable on $I$, one may think to define formally the 
Wiener integral \textit{wrt} $G$, denoted  $\int_{I}\  f(s) \ \di G_s$, by setting:
\begin{equation}
\label{oeriheorihv}
 \int_{I} \  f(s) \ \di G_s:=\int_{I}\   f(s) \cdot \frac{dG_s}{ds} \ ds = \int_{I}\   f(s) \cdot W^{(G)}_{s}  \ ds,
\end{equation}
 assuming $s\mapsto  f(s) \cdot W^{(G)}_{s}$ is \oS-integrable on $I$. More precisely we have the following definition.

\begin{defi}{(Wiener integral with respect to $G$) \\}
\label{WienWick}
For any Borel set $I$ of $\rR$ and any deterministic measurable function $f:I\mapsto\R$ such that $s\mapsto  f(s) \ W^{(G)}_{s}$ is \oS-integrable on $I$, one says that  $\int_{I}  f(s) \ \di G_s$, defined by \eqref{oeriheorihv}, is the Wiener integral of $f$ on $I$, with respect to $G$, if $\int_{I}  f(s) \ \di G_s$ belongs to $(L^2)$.
\end{defi}

%
%
%
Even if, in practice, there will often exist an integer $q$ in $\N$ such that the map $t\mapsto {|g'_t|}_{-q}$ is bounded (as it was the case in Example \ref{oevieroivjeo}), it seems more than reasonable to expect, even if $t\mapsto {|g'_t|}_{-q}$ is not bounded, that,  for any finite interval $[a,b]$ of $\rR$,
\begin{equation}
 \label{eioroeirjeoricj}
\tag{$*$}
\int_{(a,b)}  1 \ \di G_s  \text{ is well-defined and such that:}
\int_{(a,b)}  1 \ \di G_s = G_{b} - G_{a},  \text{ in  \oS}. 
\end{equation}
Thanks to Equality \eqref{oeriheorihv} and  Theorem \ref{tardileomalet}, it is clear that \eqref{eioroeirjeoricj} entails, among other consequences, that:  
\vspace{-4ex}

\begin{align}
\label{ozirfjoerfidezdzejeoij}
& \bullet \text{The map } s\mapsto <g'_{s},\eta>  \text{ belongs to } \ L^1_{l\widetilde{oc}}(\rR),   \text{ for every } \eta \text{ in } \sS(\R),\\
\label{ozirfjoerfidezdzejedzzdoij}
& \bullet \int_{(a,b)} \ <g'_{s},\eta> \ ds =\ <g_{b},\eta> - <g_{a},\eta>,  \text{ for every } \eta \text{ in } \sS(\R)  \text{ and } (a,b) \text{ in } {\rR}^{2}.
\end{align}

Besides, using Proposition \ref{laiusdhv}, it is easy to establish that:
\begin{equation}
 \label{oiezdozeijdozeijo}
|<F,\varphi>| \leq {|F|}_{-q} \ {|\varphi|}_{q}\ ; \ \hspace{3ex} \forall (F,\varphi,q) \in  \sS'(\R) \times \sS(\R)  \times \N.
\end{equation}

In view of \eqref{oiezdozeijdozeijo}, it appears that Assumption \Aaiii is almost necessary to get \eqref{eioroeirjeoricj}, if one deals with Pettis integrals, and necessary if one deals with Bochner integral\footnote{See Appendix \ref{appendiceB} for precisions about Bochner integrals.}. Moreover, and by the very definition of the space $\sS'(\R)$ as the inductive limit of the sequence ${(\sS_{\hspace{-0.15cm}-p}(\R))}_{p \in \N}$, Assumption \Aaii entails that, for every compact set $
 \cK$ of $\rR_{D}$, there exists an integer $q$ in $\N$ such 
 that ${|g'_{t}|}_{-q} <+\infty$, for every $t$ in $\cK$. Thus \Aaiii appears 
 to be only a slight reinforcement of \Aaiip. 
Besides, it is clear that \eqref{ozirfjoerfidezdzejedzzdoij} is nothing but Assumption \Aaivp.
Thus, the simple considerations given in \eqref{eioroeirjeoricj}, as well 
as the ones given right above Remark 
\ref{ozeufhioeruvheirezzezeedzedzedzdu} (about \Aai \& \Aaiip), entail that Assumption \Aa is almost minimal (\ie necessary) to get a reasonable notion of Wiener integral. We will show further
 that these assumptions are also sufficient to provide us with a general non-anticipative stochastic integral.
Denote $\cE(\rR)$ the set of step functions on $\rR$. We have the following property.

\begin{prop}
\label{zpofeazzezazakpsokfepsofkspkod}
For any $f$ in  $\cE(\rR)$, $\int_{\rR} \hspace{-0.05cm} f(u) \  \di G_u$ is a Wiener integral with respect to $G$.
Moreover, let $[a,b]$ be a finite interval of $\rR$, then $\int^b_a \  \di G_u = G_b - G_a$ almost surely.
\end{prop}

\begin{pr}
Fix $\eta$ in $\sS(\R)$. From $(ii)$ of Theorem \ref{tardileomalet}, $t\mapsto S(f(t)  \ W^{(G)}_t)(\eta)$ is measurable on $\rR$. Moreover we have, thanks to Lemma \ref{dede},  $|S(f(t)  \ W^{(G)}_t)(\eta)| \leq   {{|g'_t|}_{-q}} \  \underset{t\in\rR}{\sup}{|f(t)|}  \ e^{{|\eta|}^2_{-q}}$, where $q$ is the integer given by Assumption \Aaiii. 
Theorem \ref{peodcpdsokcpodfckposkcdpqkoq} then applies and entails that $f$ is 
$dG$-integrable on  $(a,b)$. 
Furthermore, thanks to Lemma \ref{dkdskcsdckksdksdmksdmlkskdm}, one has the equality: $S(\int^b_a \  dG_u)(\eta) = \int^b_a \ S(W^{(G)}_u)(\eta)\ du =   [S(G_u)(\eta)]^b_a = S(G_b - G_a)(\eta)$. 
%
%
%
%
  The equality, in $({\cS}^*)$,  follows from the injectivity of the S-transform. Finally, since $G_b - G_a$ belongs to $(L^2)$, the equality $\int^b_a \  \di G_u = G_b - G_a$ holds in $(L^2)$ and hence almost surely.
\end{pr}

 The following theorem gives a sufficient condition for an integral, of the form \eqref{oeriheorihv}, to be a Wiener integral. Denote $\leadsto$ the equality in law.

\begin{theo}
\label{sdsiduhsdiuhfsdiudfhdsufhisduhfisduhfsiuhfsidhfiuhdfisudhfiu} 
Assume that there exists $q_0$ in $\N$ such that the 
map $s \mapsto  f(s)\cdot {|g'_s|}_{-q_0} $ belongs to 
$L^1(\rR)$. Then $Z:=\int_{\rR} f(s) \ \di G_{s}$ is an 
element of  ${(\cS)}^*$, which verifies,  $Z 
=  \sum^{+ \infty}_{k = 0}  ( \int_{\rR} f(s)  \ 
{<g'_s,e_k>}\ ds )\     {<.,e_k>}$ in \oS. 
Moreover $Z$ is a Gaussian random variable if and 
only if $\sum^{+ \infty}_{k = 0}  {\big( \int_{\rR} f(s) \ 
{<g'_s,e_k>}\ ds \big)}^2 < +\infty$. In this latter case, 
on has:
\begin{equation*}
Z \leadsto \cN\bigg(0, \sum^{+ \infty}_{k = 0}  {\big(\int_{\rR}  f(s)\ {<g'_s,e_k>}\ ds \big)}^2\bigg).
\end{equation*}
%
%
\end{theo}


%
%
%

\begin{pr} 
In order to show that equality $\int_{\rR}  f(s) \ \di G_s =  \sum^{+ 
\infty}_{k = 0}  \left( \int_{\rR} f(s) \ {<g'_s,e_k>} \ ds  \right)     {<.,e_k>}
$ holds in ${(\cS)}^*$, let us establish points {\bfseries a), b)} and 
{\bfseries c)} below.

{\bfseries a)} $s\mapsto f(s)\cdot  W^{(G)}_s$ is ${(\cS)}^*$-integrable over $\rR$. 
One can use Thm. \ref{peodcpdsokcpodfckposkcdpqkoq} since one has, $\forall (\eta,s)$ in $\sS(\R)\times \rR_{D}$ and using Lemma \ref{dede}, $|S(f(s) {W^{(G)}_s})(\eta)|  \leq {|f(s)| \  {\|W^{(G)}_s\|}_{-q_0}} e^{ { |\eta|}^2_{q_0}}   \leq  |f(s)|\  {|g'_s|}_{-q_0}  \  e^{ { |\eta|}^2_{q_0}}$.
%

{\bfseries b)} $\Psi_f:=\sum^{+ \infty}_{k = 0}  \left( \int_{\rR} f(s)\  {<g'_s,e_k>} \ ds  \right)  {<.,e_k>}$ belongs to $({\cS_{-p_0}})$,  as soon as $p_0 \geq q_0 +1$. 
Let $p_0$ be in $\N$ such that $p_0\geq q_0 
+1$. Recall the definition of  $\cR_n$ given in \eqref{zoijeroif}. Proposition \ref{laiusdhv} and \eqref{oiezdozeijdozeijo} entail that ${\|\Psi_f\|}^2_{-p_0} \leq  {\|s\mapsto f(s)\cdot {|g'_s|}_{-q_0} \|}^2_{L^1(\rR)}  \ \cR_{p_0-q_0} < +\infty$.
%
%
%
%
%

{\bfseries c)}  $\Phi_f:=\int_{\rR}  f(s) \ \di G_s$ is equal to  $\Psi_f$  in ${(\cS)}^*$.  
Define the ${(\cS)}^*$-process  $\tau$ and the family of  ${(\cS)}^*$-processes ${(\tau_N)}_{N\in\N}$ by setting, for every 
real $s$, $\tau(s) :=  \sum^{+ \infty}_{k = 0}  f(s) <g'_s,e_k> {<.,e_k>}$,  and $\tau_N(s) := \sum^{ N}_{k = 0}  f(s)   
<g'_s,e_k> {<.,e_k>}$. Obviously we have  $\Phi_f = \int_{\rR} \tau(s) \ ds$, $\Psi_f = {\lim\limits_{N \to +\infty}} \int_{\rR}  
\tau_N(s) \ ds $ in \oS. It 
then remains to show that  $ \Phi_f = \lim\limits_{N \to +\infty} \int_{\rR} \tau_N(s) \ ds$ in \oS. For this purpose, we 
use Theorem \ref{cc}. Let $(p_0,n)$ be a couple of integers with $p_{0} \geq q_0 +1$. It is easily seen that $
\tau_N$  and  $\tau$ are weakly measurable on $\rR$ (see Definition \ref{bb}) and, using the same upper-bound we used in {\bfseries b)}, that $
\tau_N(s)$  and  $\tau(s)$   belong to  $(\cS_{-p_0})$  for every real $s$. Moreover, using Proposition \ref{laiusdhv} and, again, the upper-bound we used in {\bfseries b)}, it is clear that  both functions $s \mapsto  {{\|\tau_N(s)\|}}_{-p_0}$ and $s \mapsto  {\|\tau(s)\|}_{-p_0}$ belong to $L^1(\rR,ds)
$  since ${\|\tau_N(s)\|}^{2}_{-p_0} \leq  {\|\tau(s)\|}^{2}_{-p_0} \leq  {f^{2}(s)}\  {{|g'_s|}^2_{-q_0}}  \ \cR_{p_0-q_0} $. We hence 
have shown that both functions $
\tau_N(.)$  and  $\tau(.)$ are Bochner integrable on $\rR$. Besides, for every $(n,m)$ in $\N^2$ with $n\geq m$, 
we have, thanks to the previous upper bound, $\int_{\rR}  {\|\tau_n(s)-\tau_m(s) \|}_{-p_0}  ds \leq   
\int_{\rR}  \|{\hspace{-0.1cm}\underset{k = m+1}
{\overset{+ \infty}{\sum}}  f(s) \ {<g'_s,e_k>} {<.,e_k>}
{{\|}_{-p_0}}}  \ ds \leq  \cR_{p_0-q_0} \cdot {\|s\mapsto f(s) 
{|g'_s|}_{-q_0} \|}_{L^1(\rR)}$.
%
%
%
%
%
It is then clear that the left hand side of the previous 
inequality tends to $0$ as $(n,m)$ tends to $(+\infty,+\infty)$. Theorem \ref{cc} (see Appendix \ref{appendiceB}) applies and establishes {\bfseries c)}.  Finally, $Z$ is the $(L^2)$-limit of a sequence of independent Gaussian variables if $ \sum^{+ \infty}_
{k = 0}  {( \int_{\rR} f(s) {<g'_s,e_k>} \ ds)}^{2}<+\infty$. The equality $\E[Z^2]=\sum^{+ \infty}_{k = 0}  {\big(\int_{\rR}  f(s)\ {<g'_s,e_k>}\ ds \big)}^2$  then becomes obvious. 
 \end{pr}

\begin{ex}
\label{oevieroivjeodzdez}
${\text{{\bfseries{1.}}}}$\  If $G$ is a Brownian motion, point 1 of Example \ref{oevieroivjeo},  Theorem \ref
{sdsiduhsdiuhfsdiudfhdsufhisduhfisduhfsiuhfsidhfiuhdfisudhfiu} as well as the equality $\sum^
{+ \infty}_{k = 0}  {\left(\int_{\R}  f(s)\ {<\delta_s,e_k>}\ ds \right)}^2 = {\|f\|}^{2}_{L^2(\R)}$, allow us to define the Wiener integral of $f$, in sense of Definition \ref{WienWick}, for 
any $f$ in $L^{2}(\R)$. This shows that our definition of Wiener integral \textit{wrt} Brownian motion and the usual one both coincide exactly. Besides,  it is clear that $I_{f}:=\int^{1}_{0} f(s)\diamond B_{1} \ ds$ is an $(L^{2})$ random variable if and only if $f$ belongs to $L^{2}([0,1])$. Therefore  Theorem \ref{peodcpdsokcpodfckposkcdpqkoq}
allows us to define the Wiener integral of $f$ \textit{wrt} Brownian bridge, in sense of Definition \ref{WienWick}, if and only if $f$ belongs to $L^{2}([0,1])$.

${\text{{\bfseries{2.}}}}$\ The case of Wiener integral \textit{wrt} fBm (\textit{resp.} wrt mBm) has been treated in \cite[Section 4]{JLJLV1} (\textit{resp. in \cite[Sections 2.3 \& 4]{JLJLV1})}.  
In view of, the previous point of this example,  Example \ref{oevieroivjeo} and  Theorem \ref{sdsiduhsdiuhfsdiudfhdsufhisduhfisduhfsiuhfsidhfiuhdfisudhfiu}, one can extend \cite[Proposition $4.31$]{JLJLV1} and claim
 that $\int_{\R}  f(s) \ \di B^{H}_s$ is the Wiener integral of $f$, \textit{wrt} $B^{H}$, for every function $f$ in $L^{1}(\R) \cap L^2_{H}(\R)$, where $ L^2_{H}(\R)$ has been defined in Subsection \ref{eoajeorijzp}. The fonctions for which one can defined a Wiener integral wrt mBm are included into $\overline{\cE(\R)}^{{<,>}_{h}}$, where ${<,>}_{h}$ denotes the inner product, defined (in \cite[Sections 2.3 \& Proposition 3.1.]{JLJLV1}) by setting ${<\i1_{[0,t]},\i1_{[0,s]} >}_{h} = R_{h}(t,s)$, and where $R_{h}$ has been defined in \eqref{azazazaz}.

${\text{{\bfseries{3.}}}}$\
 In the case of $V_{\gamma}$ - processes, one can improve Proposition 14 of \cite{MV05}. Indeed, denote $\cH$ the set of all functions for which \cite[Section 3]{MV05} define a Wiener integral. For any 
 $\eta:[0,T]\rightarrow \R_{+}$, continuous and increasing in a neighborhood of $0$ and such 
that $\lim_{0+} \eta =0$, define the set $\C^{\eta}:=\{f
\in L^{2}([0,T]), \ \underset{0\leq r<s \leq T}{\sup} |f(s)-f(r)| / \eta(s-r) <+
\infty \}$. In order to show that  $\cH$ contains $\C^{\eta}$,  \cite[Proposition 14]{MV05} has to require an additional assumption on $\eta$. No such assumption is required here. 
Using only 
Theorem 
\ref{sdsiduhsdiuhfsdiudfhdsufhisduhfisduhfsiuhfsidhfiuhdfisudhfiu}, one easily sees that, for every $T>0$, the process 
$({Z_t)}_{t \in [0,T]}$, defined by $Z_t:= \int^t_{0}  f(s) \ \di G_s$,  where $f$ belongs to $\C^{\eta}$, is a Gaussian process. One just needs to see that, 
$|f(s)|\cdot {|\Phi'(s)|}_{-q_0}   \leq  (M + |f(r_{0})|) \  (1+{|\Phi'(s)|}_{-q_0})$, for every $(f,s)$ in $\C^{\eta}\times [0,T]$, where $M:=\underset{0 \leq u \leq T} {\sup} \    {\eta(u)} \cdot  \underset{0\leq 
r<s \leq T}{\sup} |f(s)-f(r)| / \eta(s-r)$, $\Phi$ and $\Phi'$ have been 
defined in Theorem \ref{pzoefkpoezzokfzepokzpzefj1286564536} and 
$r_{0}$ is any real in $(0,T]$ such that $|f(r_{0})|<+\infty$. Using 
Theorem \ref{pzoefkpoezzokfzepokzpzefj1286564536}, one concludes 
that $s\mapsto |f(s)|\cdot {|\Phi'(s)|}_{-q_0}$ belongs to $L^{1}((0,T))$ and 
one then uses Thm. 
\ref{sdsiduhsdiuhfsdiudfhdsufhisduhfisduhfsiuhfsidhfiuhdfisudhfiu} 
\end{ex}
\begin{rem}
\label{aokaokzokokok}
In fact one can extend the notion of Wiener integral \textit{wrt} any $G$ in $\mathscr{G}$ in two ways. The first way, which is also the more general one, is given in Point  ${\text{{\bfseries{1.}}}}$ of Remark \ref{ozeidjzeoijdzeo}. The second way it is the following:  If the bilinear form ${<,>}_{R}$, defined on $\cE(\rR)\times \cE(\rR)$ by setting ${<\i1_{[0,t]},\i1_{[0,s]}>}_{R}:= R(t,s)$ is an inner product; assuming there exists an isometry, denoted $M:(\cE(\rR),{<,>}_{R})\rightarrow (L^{2}(\R),{<,>}_{L^{2}(\R)})$, such that  $g_{t}:=M(\i1_{[0,t]})$, then one can extend the notion of Wiener integral to any elements of $\overline{\cE(\rR)}^{{< >}_{R}}$. This latter space contains in general not only functions but also tempered distributions. This general method applies to fBm and mBm (see \cite[Section 3]{JLJLV1}), as well as to Volterra processes.
\end{rem}

\begin{rem}
As it is explained in \cite[Example 3.3.]{SoVi14} the Brownian bridge 
admits several representations, among which are the orthogonal one, 
the Fredholm one and the canonical one. It is clear that both the orthogonal and the canonical representations of the Brownian bridge on $[0,T]$ fulfill 
Assumption \Da on $\rR:=[0,T]$
   (since  there exists $q$ in $\N^{*}$ such that $t\mapsto {|g'_t|}_{-q}$ belongs to $L^{1}([0,b],dt)$, for every $b$ in $[0,T)$. This result can be extended to Gaussian bridges (see \cite{GSV07} for more details about this latter notion), assuming the Gaussian process $G:={(G_{t})}_{t\in[0,T]}$ fulfills Assumption \Aap.
\end{rem}

%



\subsection{The Wick-Itô integral with respect to Gaussian processes}
\label{ozieffjioefjoezifjioezzjioefjfeoiz}						

%

We still assume in this section that Assumption \Aa holds. We are now able to define, and give the main properties, of the Wick-Itô integral \textit{wrt} $G$.
We still denote $I$ a Borel set of $\rR$ and 
let $X:={(X_{t})}_{t\in \rR}$ be an \oS-valued 
process. Because the belonging to \oS\ is not stable by ordinary 
product, one can not generalize \eqref{oeriheorihv} to any 
\oS-valued process $X$, by simply  setting:
\vspace{-0.25cm}

\begin{equation*}
 \int_{I}  X_{s} \ \di G_s:=\int_{I} X_{s} \cdot \frac{dG_s}{ds} \ ds = \int_{I}  
 X_{s}  \cdot W^{(G)}_{s}  \ ds.
\end{equation*}

However, since the  belonging to \oS\  is  stable by Wick product one may extend the integral \eqref{oeriheorihv} to \oS-valued processes $X$ in the following manner.

\begin{defi}[Wick-Itô integral \textit{wrt} Gaussian process]
\label{oezifhherioiheroiuh}
Let $X:\rR \rightarrow {(\cS)}^*$ be a process such that the process 
$t \mapsto X_t \diamond W^{(G)}_t$ is  \oS-integrable on $\rR$. 
The process $X$ is then said to be $dG$-integrable on $\rR$ (or 
integrable  on $\rR$), \textit{wrt} the Gaussian process $G$. The 
Wick-Itô integral of $X$ \textit{wrt} $G$, on $\rR$, is defined by 
setting:
\vspace{-1ex}

\begin{equation}\label{eigfrretth}
\int_{\rR} X_s \ \di G_s :=  \int_{\rR} X_s \diamond W^{(G)}_s \ ds.
\end{equation}

For any Borel set $I$ of $\rR$, define $\int_{I} X_s \ \di G_s:= \int_{\rR} \ {\i1}
_{I}(s) \ X_s \ \di G_s$.
\end{defi}

 The Wick-Itô integral of an $ ({\cS}^*)$-valued process, \textit{wrt} $G$ is then an element of \oS. It is easy to see that Wick-Itô integration \textit{wrt} $G$, is linear and that Definition \ref{oezifhherioiheroiuh} is coherent with Definition \ref{WienWick}, of Wiener integral, we gave in the previous subsection. Moreover, and as it will be stated in Proposition \ref{zpofkpsokfepsofkspkod} below, one of the advantages of Definition \ref{oezifhherioiheroiuh} is that our integral \textit{wrt} $G$ is centered, assuming it belongs to $(L^2)$. 
Used a lot in the sequel of this paper, the following condition ensures the integrability, on $I$, of an \oS-valued process $X$, \textit{wrt} $G$.



\vspace{0.5ex}

Let $X:I \rightarrow {(\cS)}^*$ be an \oS-valued process. Denote the following condition:
\textcolor{white}{
\begin{subequations}
 \begin{align}
  a  = b \tag*{$\mathcal{I}$} \label{i}\\
  a  = b \tag*{$\mathcal{I}_{p,q}$} \label{ipq}
 \end{align}
\end{subequations}
}
\vspace{-11ex}
\hspace{7ex}
$$\mathcal{(I)}\hspace{0.25cm} \begin{cases} (i): \hspace{0.25cm} t \mapsto S(X_t)(\eta)  \text{ is measurable on } I, \text{ for all } \eta \text{ in } \sS(\R) \\
(ii): \hspace{0.25cm} \exists (p,q) \in \N^2\text{ such that the map } t\mapsto {\| X_t\|}_{-p}  {\| W^{(G)}_t\|}_{-q}  \text{  belongs to } L^1(I,dt). \end{cases} 
$$

%
%
%
%
%
%
%

When the processes $X$ and $G$ satisfy condition \iIp on $I$, we will say that $(X,G)$ satisfies \iIp or \iIlp, when we want to specify the value of $p$ and $q$ in $(ii)$ of Condition \iIpp. We will use the following theorem a lot in the sequel. 

\begin{theo}
\label{dazdzedzedzedzedzedze}
If $(X,G)$ satisfies condition \iIlp on $I$, then $\int_{I} X_s \ \di G_s$ is well-defined and belongs to $({\cS_{-r}})$ for every $r\geq 2+\max \{p;q\}$. Moreover there exists a real constant $C$, independent of $X$ and $G$ , such that: 
\begin{equation*}
\forall \ r\geq 2+\max \{p;q\}; \hspace{0.75cm} {\big\| \int_{I} X_t \ \di G_t \big\|}_{-r} \leq C \  \int_I {\| X_t\|}_{-p}  \ {\| W^{(G)}_t\|}_{-q}\ dt.  
\end{equation*}
\end{theo}
\begin{pr}
$\forall \eta \in \sS(\R)$, the measurability on $I$ of  $t \mapsto S(X_t\diamond W^{(G)}_t)(\eta)$ is clear since
$S(X_t\diamond W^{(G)}_t)(\eta) =  S(X_t)(\eta) \ {<g'_t,\eta>}$.
Condition \iIp being verified, we use Lemma \ref{dede} to get, for all $\eta$ in $\sS(\R)$, $|S(X_t \diamond W^{(G)}_t)(\eta)|\leq  e^{{|\eta|}^2_{\max \{p;q\}}}\   {\|X_t\|}_{-p} \  {\|W^{(G)}_t\|}_{-q}$, where
$p$ and $q$ are given by condition \iIp\hspace{-1ex} . Theorem \ref{peodcpdsokcpodfckposkcdpqkoq} then clearly applies.
The upper-bound in the theorem, as well as the existence of $r$ and $C$, results from \cite[Thm. 13.5]{Kuo2}.
 \end{pr}
 
%
%
%
%
%
%
%
%
%
%
%
%
%
 

 We can now give the first properties  of the Wick-Itô integral \textit{wrt} $G$.
%
%
 
\begin{prop}
\label{zpofkpsokfepsofkspkod}
\bit
\iti Let $I$ be a Borel subset of $\rR$ and $X:\hspace{-0.1cm}I\rightarrow  $\oS$ $ a  $dG$-integrable process over $I$. Assume that 
$\int_{I} X_s \ \di G_s$ belongs to $(L^2)$. Then $\E[\int_{I}  X_s \  \di G_s] = 0$.


\itii Let $[a,b]\subset \rR$. The $({\cS}^*)$-process $\Psi$ defined by $ \Psi(t):= \int^t_{a}  \ X_s \ \di G_s$ is  continuous on $[a,b]$, as soon as $(X,G)$ satisfies condition \iIp on $[a,b]$.
\eit
\end{prop}

{\bfseries Proof.} \textit{$(i)$} That $S(\int_{I} \ X_s \ \di G_s)(0) = 
\int_{I} S(X_s)(0) \ S(W^{(G)}_s)(0) \ ds = 0$ is clear since 
$S(W^{(G)}_s)(0) = {<g'_s,0>} = 0$. Now, it sufficient to note 
that  $\E[U] = S(U)(0)$ for every {\textit r.v.} $U$ in  $(L^2)$. 

\par \noindent \textit{$(ii)$} The integrability of $X$ 
\textit{wrt}  $G$ is proved by Theorem 
\ref{dazdzedzedzedzedzedze}. Let $t_0$ be fixed in 
$(a,b)$. In order to establish the continuity of $\Psi$ 
in $t_0$ we are going to use \cite[Theorem 8.6]
{Kuo2}. By symmetry one may assume that $t_0 \geq t
$.  \cite[Theorem 8.6]{Kuo2} applies since we clearly 
have:

\vspace{-0.25cm}

\begin{itemize}
\ita $|S(\Psi(t)-\Psi(t_0))(\eta)|\leq e^{{|\eta|}^2_{\max \{p;q\}}}\ \int^{t_0}_{t} \   {\|X_u\|}_{-p} \  {\|W^{(G)}_u\|}_{-q} \ du \underset{ t\to t_0 \textcolor{white}{ze}}{\longrightarrow 0}$;
\vspace{-0.5cm}
\itb $|S(\Psi(t))(\eta)| \leq e^{{|\eta|}^2_{\max \{p;q\}}}\ \int_{[a,b]} \   {\|X_u\|}_{-p} \  {\|W^{(G)}_u\|}_{-q} \ du $.$\hfill \square$
\eit

\begin{prop}
\label{oeijoeirjveorivjoerivjev}
Let  $(X,G)$ be a couple of processes that satisfies condition \iIp on $\rR$. Define, for every $n$ in $\N$, the process $G^{(n)}:={(G^{(n)}_t)}_{t\in\rR}$ by setting $G^{(n)}_{t} := <.,g^{(n)}_{t}>$, where $g^{(n)}_{t}$ belongs to $L^2(\R)$ and let ${(X^{(n)})}_{n\in\N}:=\{{(X^{(n)}_t)}_{t\in\rR};\ n\in\N\}$ denote a sequence of \oS-valued processes. 
Let us write the following conditions:
\bit
\itap $(X,G^{(n)})$ satisfies condition \iIp on $I$,  uniformly\footnote{\ie $\exists\ (p,q) \in \N^{2}$,  such that  $s\mapsto {\| X_s\|}_{-p}  {\| W^{(G^{(n)})}_s\|}_{-q}$    belongs to  $L^1(I,ds)$, for every $n$ in $\N$.} in $q$.
\itapp $\exists\ (r,r_{1}) \in \N\times (0,+\infty]$ such that:  ${\| {{G}^{(n)}_{.} - {G}_{.}\|}_{-r}} \underset{n\to +\infty}{\longrightarrow 0}$, \text{where the convergence holds}
 \text{both pointwisely on $I$,} and in $L^{r_1}(I)$,
\itappp  $X$ is \oS-differentiable on $I$ and there exist $(a,l)\in \N\times\R$ and a function $L\in 
L^{r_{2}}(I,dt)$ s.t.
\begin{equation*}
|\tfrac{d}{ds}[S(X_{s})(\eta)]| 
\leq L(s) \  e^{a {| \eta |}^2_{l} }, 
\end{equation*}
 for all $\eta$ of $
\sS(\R)$ and for \textit{a.e.} $s$ of $I$, where $r_{2} \in (0,+\infty]$ is such that ${r^{-1}_{1}}+{r^{-1}_{2}} =1$,
\eit
\bit
\itappq $\exists (r_{1}, r_{2}) \in {(0,+\infty]}^{2}$ with  ${r^{-1}_{1}}+{r^{-1}_{2}} =1$, such that:
\vspace{-0.25cm}
\bit
\iti $\|{s\mapsto |{g}'^{(n)}_{s} - {g}'_{s}|}_{-q}\|_{L^{r_1}(I)} \underset{n\to +\infty}{\longrightarrow 0}$ ; \hspace{0.5cm} $(ii)$
\hspace{0.1cm}  $s\mapsto {\|X_{s}\|}_{-p}$ belongs 
to $L^{r_{2}}(I)$.
\eit
\eit
If conditions $(a_{i})_{i\in\{1;2;3\}}$ or both conditions $(a_{1})$ and $(a'_{2})$ are fulfilled, then one has the convergence:
\begin{equation*}
\int_{I}  X_s \ \di G^{(n)}_s \underset{n\to+\infty}{ \longrightarrow } \int_{I}  X_s \ \di G_s \hspace{0.5cm} \text{in \oS}. 
\end{equation*}

Besides, denote the following conditions:

 \vspace{-0.2cm}

\bit

\itaps $(X^{(n)},G)$ satisfies condition \iIp on $I$, uniformly
in $p$.

\itapps $\exists\ r \in \N$ such that:  ${\| {{X}^{(n)}_{.} - {X}_{.}\|}_{-r}} \underset{n\to +\infty}{\longrightarrow 0}$ pointwise.
\itappps  Both $X$ and ${(X^{(n)})}_{n\in\N}$ are \oS-differentiable on $I$. Moreover, $\exists (l,a,r_{2})\in\N\times \R \times (0,+\infty]$ and a function $L$ in 
$L^{r_{2}}(I,dt)$ such that, for every $(n,\eta)$ in $\N\times \sS(\R)$ and \textit{a.e.} $s\in I$,
\begin{equation*}
|\frac{d}{ds}[S(X^{(n)}_{s})(\eta)]| + |\frac{d}{ds}[S(X_{s})(\eta)]| 
\leq L(s) \  e^{a {| \eta |}^2_{l} }, 
\end{equation*}
 \itapppps \bit
\iti $\exists\ r' \in \N$, \textit{s.t.} $\|s\mapsto {\|\frac{d}{ds}[X_{s}-X^{(n)}_{s}]\|}_{-r'}\|_{L^{r_2}(I)} \underset{n\to +\infty}{\longrightarrow 0}$ \hspace{0.3cm}  \text{and}  \hspace{0.25cm} $(ii)$
\hspace{0.05cm}  $s\mapsto {|g_{s}|}_{-q}$ belongs 
to $L^{r_{1}}(I)$,
where $r_{1} \in (0,+\infty]$ is such that ${r^{-1}_{1}}+{r^{-1}_{2}} =1$.
\eit

\itappqs $\exists (r_{1}, r_{2}) \in {(0,+\infty]}^{2}$ with ${r^{-1}_{1}}+{r^{-1}_{2}} =1$, such that:
\vspace{-0.25cm}
\bit
\iti $\|s\mapsto {\|X_{s}-X^{(n)}_{s}\|}_{-p}\|_{L^{r_2}(I)} \underset{n\to +\infty}{\longrightarrow 0}$
to $L^{r_{2}}(I)$\hspace{0.5cm} $(ii)$ \hspace{0.1cm} 
$s\mapsto |{g}'_{s}|_{-q} \in L^{r_1}(I)$.  
\eit

\eit

If conditions $(b_{i})_{i\in\{1;2;3;4\}}$ or both conditions $(b_{1})$ and $(b'_{2})$ are fulfilled, then one has the convergence:
\begin{equation*}
\int_{I}  X^{(n)}_s \ \di G_s \underset{n\to+\infty}{ \longrightarrow } \int_{I}  X_s \ \di G_s, \hspace{0.5cm} \text{in \oS}.
\end{equation*}

%
%
%
\end{prop}

\begin{pr}
The scheme of the proof is ``symmetric in $G^{(n)}$ and $X^{(n)}$''; we will then only show the convergence $ \int_{I}  X_s \ \di G^{(n)}_s \underset{n\to+\infty}{ \longrightarrow } \int_{I}  X_s \ \di G_s.$ Denote $A_{n}:= \int_{I}  X_s \ \di G_s - \int_{I}  X_s \ \di G^{(n)}_s$; let us show that assumptions of \cite[Theorem 8.6]{Kuo2} are fulfilled. The existence of all following integrals come form Theorem \ref{peodcpdsokcpodfckposkcdpqkoq}.
\vspace{0.1cm}

$\bullet$ Case where $a_{1}$) \& $a'_{2})$ are fulfilled:
\vspace{0.1cm}

\cl{${|S(A_{n})(\eta)|}^{r_{1} r_{2}}= {(|\int_{I}   S(X_{s})(\eta) \ S(W^{(G^{(n)})}_{s}\hspace{-0.2cm}-W^{(G)}_{s})(\eta) \ ds|)}^{r_{1} r_{2}} \hspace{-0.25cm}\leq (\int_{I} {\|X_{s}\|}^{r_{2}}_{-p} \ ds) (\int_{I} {|{g}'^{(n)}_{s} - {g}'_{s}|}^{r_{1}}_{-q} \ ds.$}
\vspace{0.2cm}

$\bullet$ Case where $(a_{i})_{i\in\{1;2;3\}}$ are fulilled :
\vspace{0.1cm}
Let us assume that $I=[0,t]$. An integration by parts yields:
\vspace{-0.5cm}

\begin{align*}
|S(A_{n})(\eta)| &= |\int_{I}\   S(X_{s})(\eta) \ S(W^{(G^{(n)})}_{s}-W^{(G)}_{s})(\eta) \ ds| = |\int_{I}\   S(X_{s})(\eta) \ <{{(g^{(n)}_{.})'}}_{s} - g'_{s},\eta> \ ds| \\
& = <g^{(n)}_{t}-g_{t},\eta> S(X_{t})(\eta) - \int_{I}\ <g^{(n)}_{s}-g_{s},\eta> \ \tfrac{d}{ds}[S(X_{s})(\eta)]   \ ds\\
& \leq e^{(1+a) {|\eta|}^{2}_{(p+r+l)}} (\|X_{t} \|_{-p} {\| {{G}^{(n)}_{t} - {G}_{t}\|}_{-r}} \ + \int_{I}  {\|{G}^{(n)}_{s} - {G}_{s}\|}_{-r}  \ L(s)\ ds).
\end{align*}
The H\"older inequality then allows one to establish the two conditions of \cite[Theorem 8.6]{Kuo2} and therefore achieves the proof.
\end{pr}
\begin{rem}
\label{apqorperzokaqepi}
 ${\text{{\bfseries{1.}}}}$\ The advantage of condition \iIp is that it allows  
 us to make assumptions on both elements of the couple $(X,G)$ 
 instead of making assumptions only on $X$ or only on $G$. 
Thus, the more informations on  the 
 ``regularity'' of $X$ (\textit{resp.} of $G$) one gets, the less informations 
 one needs on the ``regularity'' of $G$ (\textit{resp.} of $X$). 

 ${\text{{\bfseries{2.}}}}$\ It is clear, in $a_{2})$ of $(iii)$ of Proposition \ref{oeijoeirjveorivjoerivjev}, that one can also choose the pointwise convergence in $(L^{2})$ or in probability instead of convergence in $(\cS_{-r})$. 

 ${\text{{\bfseries{3.}}}}$\ When $G$ 
is an fBm (resp. an mBm), the 
Wick-Itô integral \textit{wrt} $G$ given by Definition \ref{oezifhherioiheroiuh}
is nothing but the  
fractional (resp. multifractional) Wick-Itô integral 
defined in \cite{ell,bosw,ben1} (resp. in 
\cite{JLJLV1,LLVH, JL13}).
\end{rem}
It is of interest to have also a criterion of integrability for generalized functionals of $G$. This will provide a very simple proof of the fact that both $\int^b_a \ f(G_t) \ dt$ and $\int^b_a \ f(G_t) \di G_t$ exist in \oS.
%
%
\begin{theo}
\label{dlsdqsdhuqsgdudguqsgdqudgquygdsudgqsuydgqsuydgqsuyqdgq}
Let $p$ be in $\N$, $[a,b]$ be an interval of $\Z^c_R$ and let $F$ be in $\sS_{\hspace{-0.15cm}-p}(\R)$. 
If $t\mapsto \max \{R^{-p-1/4}_{t}; R^{p-1/4}_{t} \}$ belongs to  $L^1([a,b])$ (\textit{resp.} there 
exists an integer $q$  such that the map $t\mapsto {|g'_t|}_{-q}\max \{R^{-p-1/4}
_{t}; R^{p-1/4}_{t} \}$ belongs to  $L^1([a,b])$), then the stochastic 
distribution process $F(G_t)$ is ${(\cS)^*}$-integrable  (\textit{resp.} $dG$-integrable) on $[a,b]$ (\textit{resp.} on $(a,b)$).
\end{theo}

\begin{pr}
Lemma \ref{dede} and Equality \eqref{spdsdpokpovkopsdkvspodvksopsdvpoksvdokdsvkopsvdopqbschfgsdvhjerhbtkljtyopfdhvlkdfhu} both apply and allow us to use Theorem \ref{peodcpdsokcpodfckposkcdpqkoq}.
%
%
\end{pr}

\begin{rem}
\label{iuchiushchius}
Of course conditions of Theorem \ref{dlsdqsdhuqsgdudguqsgdqudgquygdsudgqsuydgqsuydgqsuyqdgq} are obviously verified when the infimum of 
$t\mapsto R_{t}$ on $[a,b]$ is positive and when its supremum is upper-bounded on $[a,b]$. Moreover, in the particular case where these latter conditions hold, Theorem \ref{dlsdqsdhuqsgdudguqsgdqudgquygdsudgqsuydgqsuydgqsuyqdgq} entails that both quantities $\int^b_a \ f(G_t) \ dt$ and $\int^b_a \ f(G_t) \ \di G_t$ exist in \oS, as soon as $f$ is a function of polynomial growth.
\end{rem}

%


\begin{ex}[Computation of $\int^T_0 \ G_t \ \di G_t$]
\label{zeodijze}
Let $T>0$. Assume that $[0,T]\subset \rR$ and that $t\mapsto R_{t}$ is 
upper-bounded on $[0,T]$,  then  the following equality holds almost surely 
and in $(L^2)$.
\begin{equation}
\label{eoijeoriccc}
\int^T_0 \ G_t \ \di G_t = \tfrac{1}{2} \ {(G}^{2}_T - R_{T}).
\end{equation}
%
%
%
\end{ex}
This result will be obtained as a direct consequence of Itô formulas provided in Section \ref{Itô}. The direct proof is therefore left to the reader.

\begin{rem}  In the previous example, we could have replaced the assumption $t\mapsto R_{t}$ is 
upper-bounded on $[0,T]$ by $\int^{T}_{0} \   R_t  \cdot {|g'_{t}|}_{-q} \ dt <+\infty$. 
\end{rem}

To end this section, we present a simple but classical stochastic differential equation, driven by a Gaussian process. We need first to generalize the definition of the \textit{Wick exponential}, given at the beginning of Subsection \ref{ksksksks}, to 
the case where $\Phi$ belongs to $(\cS)^{*}$. For any $\Phi$ in \oS and $k$ in $\N^*$ let $\Phi^{\diamond k}$ denotes the Wick product of $\Phi$, taken $k$ times. For any $\Phi$ in \oS such that the sum $\sum^{+\infty}_{k = 0} \frac{\Phi^{\diamond k}}{k!}$ converges in $(\cS)^{*}$, define the Wick-exponential of $\Phi$,  and denote $\exp^{\diamond}\Phi$, the element of  \oS defined  by $\exp^{\diamond}\Phi:= \sum^{+\infty}_{k = 0} \frac{\Phi^{\diamond k}}{k!}$. For $f$ in $L^2(\R)$ and  $\Phi:= <.,f>$, it is easy to verify that $\exp^{\diamond} \Phi = :e^{<.,f>}:$.


%
%


%

\begin{ex}[The Gaussian Wick exponential]
Let $\rR=\R_{+}$ and let us consider the following Gaussian SDE:
\vspace{-0.5cm}
\begin{align*}
\label{}
(E):\hspace{0.1cm}     {  \left\{
\begin{aligned}
dX_t &=  \alpha(t) X_t \ dt + \beta(t) X_t \ \di G_t   \\
X_0 &\in   {(\cS)}^*, 
\end{aligned}
\right. }
\end{align*}
%

where $\alpha:\R_{+}\rightarrow\R$  and  $\beta:\R_{+}\rightarrow\R$ are two deterministic continuous functions. Of course $(E)$ is a shorthand notation for $X_t = X_0 + \int^t_0 \  \alpha(s) \ X_s \ ds + \int^t_0 \  \beta(s) \ X_s \ \di G_s$.
%
%
%
As in \cite{HOUZ}, it is easy to guess the solution.
Let us define the process $Z$ by setting:


\begin{equation}
\label{msmslslqmls}
Z_t := X_0 \diamond \exp^\diamond { \bigg(  \int^t_0 \alpha(s) \ ds +   \int^t_0 \beta(s) \ \di G_s \bigg)} ,\hspace{0.5cm} t \in \R_+,
\end{equation}

\end{ex}

\begin{theo}
\label{eprfjerpoifjeroi}
The process $Z$ defined by \eqref{msmslslqmls} is the unique solution, in \oS, of (E).
\end{theo}

\begin{pr}
This is a straightforward application of \cite[Theorem $3.1.2$]{HOUZ}.
\end{pr}


%
%


%
%
%

\section{Itô Formula}
\label{Itô}
The main result of this section is Theorem 
\ref{ergoigfjerfijoijgfoifjgoidfjgofd}, which provides an Itô Formula in $
(L^2)$, for $\cC^{1,2}$ functions, with sub-exponential growth. This 
latter result is given in Subsection \ref{pokprofke}, while the end of this 
section is devoted to a complete comparison between our Itô formula and all the Itô formulas for Gaussian processes provided so far in the 
literature of \textit{functional extensions} of itô integral that are: \cite[Thms $1$  \& $2$]{nualart}, \cite[Thm $31$]{MV05}, \cite[Thm $1$]{NuTa06}, \cite[Cor. $8.13$]{KRT07}, \cite[Prop. 11.7]{KrRu10}, \cite[Thm. $3.2$]{LN12}. 
It will, in particular, show the generality of the Itô formula for Gaussian processes of the form \eqref{erofkorekp} we establish here.
Let us first recall a few basic facts about \textit{Lebesgue-Stieljes} \&  \textit{Riemann-Stieljes} integrals, that will be used extensively in the remaining part of this work. 
Let $[a,b]$ be an interval of $\R$ and $j:[a,b]
\rightarrow \R$ be a function of bounded variation. Denote $\alpha_{j}$ the signed measure such that $j(t)=\alpha_{j}([a,t])$, for every $t$ in $[a,b]$. For any function $f:[a,b]\rightarrow \R$,  denote $ \int^{b}_{a} f(s) \ dj(s)$ or $ \int^{b}_{a} f(s) \ d\alpha_{j}(s)$ the \textit{Lebesgue-Stieljes integral} of $f$ with respect to $j$, assuming it exists. In this latter case, we will write  that $f\in L^{1}(I,dj(t))$ or $L^{1}(I,\alpha_{j})$.
%
%
%
In the particular case where the function $f$ is continuous on $[a,b]$, the \textit{Lebesgue-Stieljes integral} of $f$ exists and is also equal to the \textit{Riemann-Stieljes integral} of $f$, which is denoted and defined by:
\begin{equation}
\label{pzodkk}
\text{(R.S.) }  \int^{b}_{a}\ f(s)\  d j(s):= \underset{\pi \to 0}{\lim}  \sum^{n}_{i=1}\ f(\xi^{(n)}_{i}) \ (j(x_{i})-j(x_{i-1})),
\end{equation}
where the convergence holds uniformly on all finite partitions ${\mathscr {P}}^{(n)}_{\pi}:=\{a:=x_{0} \leq x_{1} \leq \cdots \leq x_{n}:=b\}$ of $[a,b]$ such that $\underset{1\leq i \leq n}{\max} (x_{i}-x_{i-1}) \leq \pi$ and such that $\xi^{(n)}_{i}$ belongs to $[x_{i-1},x_{i}]$. 
The following result, will be used extensively in the sequel of this section.

\begin{lem}
\label{zpodkzpeodkzep}
Let $[a,b]$ be a finite interval of $\R$, $I$ (\textit{resp.} $J$) an interval of  $\R_{+}$ (\textit{resp.} of $\R$) and let $L:[a,b]\times I\times J$ be a $C^{1}$-function. Let $f:[a,b]\rightarrow I$ and $j:[a,b]\rightarrow J$ be two continuous functions of bounded variation on $[a,b]$. Then one has the following equality:
\vspace{-1ex}
\begin{align}
 \label{zpdozepkddzzdzdzdpk}
 L(b,f(b),j(b)) -  L(a,f(a),j(a)) &= \int^{b}_{a}\    \frac{\partial L}{\partial u_{1}}(s,f(s),j(s))  \ ds +  \int^{b}_{a}\    \frac{\partial L}{\partial u_{2}}(s,f(s),j(s))  \ df(s) \notag \\
&+ \int^{b}_{a}\    \frac{\partial L}{\partial u_{3}}(s,f(s),j(s))  \ dj(s).
\end{align}
\end{lem}

\begin{pr}
All the integrands in the right hand side of \eqref{zpdozepkddzzdzdzdpk} are continuous. Thus the \textit{Lebesgue-Stieljes} integrals in the right hand side of \eqref{zpdozepkddzzdzdzdpk} are also \textit{Riemann-Stieljes} integrals. It is then 
easy to deduce \eqref{zpdozepkddzzdzdzdpk}, using \eqref{pzodkk}.
\end{pr}
In view of Theorem-Definition \ref{izeufheriduscheirucheridu}, it is clear that we can extend the notion of integral in \oS to the case where $m$ is a signed measure (the notation remaining the same). We will therefore keep the same notations for this integral, whatever the measure $m$ is (signed or positive).
In the remaining of this paper, and unless otherwise specify,
the measure $m$ denote a measure, that may be $\sigma$-finite or signed.
%
\subsection{Itô Formula in $(L^2)$ for $\cC^{1,2}$ functions with sub-exponential growth}
\label{pokprofke}
Let us begin with the following lemma, the proof of which is an immediate consequence of \cite[Theorems 1,2 p.88-89]{Wid}.

%

\begin{lem}
\label{podiuhyiuuhiudokdfzdokfdepsoksdpokfsdpokfspokdpfoksdpfskfodkof}
Let $T>0$ and $v:[0,T]\times \R \rightarrow \R$ be a continuous function 
such that there exists  a couple $(C_T,\lambda_T)$ of
$\R\times\R^*_+$ such that $\underset{t \in [0,T] }{\max}{|v(t,y)|} \leq C_T \
e^{\lambda_T y^2}$ for all real $y$. Define; for every $a>\lambda_T$, the map  $J_v:\R_+ \times (0,1/4a) \times \R \rightarrow \R$ by 
setting: 
\begin{equation}
\label{osicosicjsdoicj}
J_v(t,u_{1},u_{2}):= \int_{\R} \  v(t,x)\cdot \gamma(u_{1},x-u_{2}) \  dx.  
\end{equation}
Then $J_v$ is 
well defined. Moreover \hspace{-0.5cm} $\lim\limits_{(t,u_{1},u_{2}) \to (t_{0},0^+,l_0)}  
\hspace{-0.15cm} J_{v}(t,u_{1},u_{2}) = v(t_{0},l_{0})$, $\forall\  (t_{0},l_0)$ in $[0,T]\times \R$.
\end{lem}

%

%
%

It is easy to extend \cite[Thm. $2.8$]{ben2} to the case of a Borel measure $m$ instead of the Lebesgue measure. The next result, which constitutes this extension, is more suitable that Thm. \ref{peodcpdsokcpodfckposkcdpqkoq}, when one deals with $L^2$-valued integrands. The proof being obvious, is then left to the reader.
\begin{theo}
\label{ozeijdoezjdozeijdoizjdoijezoijzdijz}
Let $m$ be a positive measure on $(\rR,\cB(\rR))$ and $X:\rR \rightarrow (L^2)$ be such that the function $t\mapsto S(X_t)(\eta)$ is measurable, for all $\eta$ in $\sS(\R)$, and such that $t \mapsto {\|X_t\|}_0$ belongs to $L^1(\rR,m)$. Then $X$ is \oS-integrable over $\rR$ and verifies:
\vspace{-2ex}
\begin{equation*}
{\big\| \int_{\rR} \ X_t \ m(dt) \ \big\|}_0~\leq~\int_{\rR} \ {\|X_t\|}_0  \ m(dt). 
\end{equation*}
\end{theo}

%

Through this subsection, we assume that $T>0$ and define $\rR:=[0,T]$. We can now give the main result of this section. Denote $
\cC^{1,2}([0,T]\times\R,\R)$ the set of functions of 
two variables which belongs to $\cC^{1}([0,T],\R)$ as 
function of their first variable and to $\cC^{2}(\R,\R)$ 
as function of their second variable. The main result of this section is the following.

\begin{theo}
\label{ergoigfjerfijoijgfoifjgoidfjgofd}
Let $T >0$. Let $f$ be a $\cC^{1,2}([0,T]\times\R,\R)$ function. 
Furthermore, assume that there are constants $C \geq 0$ and $\lambda 
< {(4 {\underset{t \in [0,T] }{\max}{R_{t}}})}^{-1}$ such that for all $(t,x)$ 
in $[0,T]\times \R$,
\smallskip

\begin{equation}
\label{orijv}
{\underset{t \in [0,T] }{\max}{\big\{ \big|f(t,x)\big|,\big|\tfrac{\partial f}{\partial t}(t,x)\big|,\big|\tfrac{\partial f}{\partial x}(t,x)\big|,\big|\tfrac{\partial^2 f}{\partial x^2}(t,x)\big| \big\}} \leq C e^{\lambda x^2}}.
\end{equation}

Assume moreover that Assumption \Aa holds and that the map $t\mapsto R_{t}$ is both continuous and of bounded variations on $[0,T]$.
Then the following equality holds in $(L^2)$:
\vspace{-1ex}
\begin{equation}
\label{ojzdojzdodjz}
f(T,G_{T}) = f(0,0) + \int^T_0 \ \tfrac{\partial f}{\partial t}(t,G_{t}) \ dt + \int^T_0 \ \tfrac{\partial f}{\partial x}(t,G_{t}) \ \di G_{t} 
				+ \tfrac{1}{2} \ \int^T_0 \ \tfrac{\partial^2 f}{\partial x^2}(t,G_{t}) \ dR_{t}.
\end{equation}
\end{theo}
\begin{pr} 
The general technique of proof of the
Itô formula via the S-transform can be traced back to \cite{Kub}. The general structure of this proof is similar to the proof of  \cite[Theorem 5.3]{ben2}. However, one can not follow this latter completely since one does not assume (as it is the case for fBm) that $\z=\{0\}$, where we set $\z:= \Z_{R} \cap [0,T]$. Equality \eqref{ojzdojzdodjz} may be rewritten as:
\vspace{-1ex}
\begin{equation}
\label{jwhcsdihcsdiuchisuh}
\int^T_0 \ \tfrac{\partial f}{\partial x}(t,G_{t}) \ \di G_{t} = f(T,G_{T}) - f(0,0) - \int^T_0 \tfrac{\partial f}{\partial t}(t,G_{t}) \ dt -\tfrac{1}{2}  \int^T_0   \tfrac{\partial^2 f}{\partial x^2}(t,G_{t}) \ dR_{t}.
\end{equation}
 Thanks to \eqref{orijv} we may write, for every $K$ in $\big\{f, \tfrac{\partial f}{\partial t}, \tfrac{\partial f}{\partial x}, \tfrac{\partial^2 f}{\partial x^2}  \big\}$ and $t$ in $[0,T]$, that $\E\big[{K(t,G_{t})}^2\big] \leq M^{2}$, where we set $M^{2}:=C^2 \ {(1-4 \lambda  \overline{R}    )}^{-1/2}$ and  $\overline{R}:=\sup\{R_{t};\ t \in [0,T]\}$. Moreover, $t 
 \mapsto {\|K(t,G_{t})\|}_0$ belongs to $L^1([0,T], dt)$ while  $t\mapsto  {\| 
 \frac{\partial^2 f}{\partial x^2}(t,G_{t}) \|}_0$ belongs to $L^1([0,T], dR_{t})$. 
 The measurability of the maps $t\mapsto S(K(t,G_{t})(\eta)$ will become 
 clear thanks to \eqref{hiuiuiuiuiudd}. A simple application of Theorem \ref
 {ozeijdoezjdozeijdoizjdoijezoijzdijz} then yields that all members on the 
 right hand side of \eqref{jwhcsdihcsdiuchisuh} exist and are in $(L^2)$. Moreover, Lemma \ref{dede} provides the upper-bound $|S(\tfrac{\partial f}{\partial x}(t,G_{t})  \diamond W^{(G)}_t)(\eta)|	\leq   M \ {|g'_t|}_{-q} \ e^{{|\eta |}^2_{q}}$,
for all $(\eta,t)$ in $\sS(\R)\times [0,T]$, where $q$ is given by Assumption \Aaiii.  A straightforward application of Theorem 
\ref{peodcpdsokcpodfckposkcdpqkoq} then shows that $\int^T_0 \ 
\frac{\partial f}{\partial x}(t,G_{t}) \ \di G_{t}$ belongs to \oS. 
In order to prove the theorem, it then just remains to show that the $S$-transform of both sides of \eqref{jwhcsdihcsdiuchisuh} are equal.
For this purpose, we first give an integral representation of the S-transform of $K(t,G_{t})$.
Since  $\E[:e^{<.,\eta>}:]=1$, for every  $\eta$ in $\sS(\R)$, one can define a probability measure $\Q_{\eta}$ on the space $(\Omega,\cF)$ by setting   $\frac{d\Q_{\eta}}{d\mu} \stackrel{\text{def}}{=} :e^{<.,\eta>}:$, where  $\frac{d\Q_{\eta}}{d\mu} $ denotes the Radon-Nikodym derivative of $\Q_{\eta}$ with respect to $\mu$.  To make computations easier we use the following obvious fact: $\cL^{\mu}_{X+\sS(X)(\eta)} = \cL^{\Q_{\eta}}_{X}$, for every centered Gaussian random variable $X$ and $\eta$ in $\sS(\R)$, and where $\cL^{\rho}_{Y}$ denotes the law of a random variable $Y$ under the probability measure $\rho$.
In view of this fact, it is clear that $G_t$ is a Gaussian variable with mean ${<g_t,\eta>}$ and variance $R_{t}$, under the probability measure $\Q_{\eta}$. One then gets, for every $t$ in $[0,T]$ and $\eta$ in $\sS(\R)$:
\vspace{-1ex}
\begin{align}
\label{hiuiuiuiuiudd}
S(K(t,G_{t}))(\eta) &= \E_{\Q_{\eta}}[K(t,G_{t})] = \int_{\R} \ K\big(t,u \ {R^{1/2}_{t}} + {<g_t,\eta>} \big)  \ \tfrac{1}{\sqrt{2\pi}} \ e^{-u^2/2} \ du.
\end{align}



 Denote $\w:=[0,T]\backslash  \z$ and let $\eta$ be in $\sS(\R)$. In view of \eqref{hiuiuiuiuiudd} we get:
%
\vspace{-1ex}
\begin{empheq}[left={S(K(t,G_{t}))(\eta)=\empheqlbrace}]{alignat=2}
& \int_{\R} \ K\big(t,v) \ \gamma (R_{t}, v-<g_t,\eta>)\ dv, && \hspace{0.2cm}\forall\ t\in \w,  \\
\label{eoijecjaaa}
& K(t,0), && \hspace{0.2cm} \forall\ t\in \z.
\end{empheq}


 Let $a$ be a real in $(\lambda, {(4\overline{R})}^{-1})$. Thanks to Lemma \ref{podiuhyiuuhiudokdfzdokfdepsoksdpokfsdpokfspokdpfoksdpfskfodkof}, we know that the map $J_{K}$ is well defined on $\Sigma_{a}:=[0,T]\times (0,1/4a)\times\R$ and we clearly have:
\vspace{-1.5ex}

\begin{equation}
\label{hiuiuiuiuiuddpzofkeprofkeprofkerpfkpok}
S(K(t,G_{t}))(\eta) = J_{K}(t,R_{t},<g_{t},\eta>), \hspace{0.2cm}  \forall\ t\in \w.
\end{equation}

Moreover, it is clear that $J_{f}$ is a $C^{1}$-function on  $\Sigma_{a}$. Denote, for every $\eta$ in $\sS(\R)$,  $j_{\eta}:[0,T]\rightarrow \R$ the map defined by $j_{\eta}(t):= <g_{t},\eta>$. According to Point $1$ of Remark \ref{ozeufhioeruvheirezzezeedzedzedzdu}, $j_{\eta}$ is absolutely continuous on $[0,T]$.
 We first have the following result.
 
\begin{lem}
\label{epokep}
Let $a$ be a real in $(\lambda, {(4\overline{R})}^{-1})$ and $J$ be the map, 
defined on $\Sigma_{a}$ by \eqref{osicosicjsdoicj}. For every $(t,\eta)$ in $\w\times \sS(\R)$, 
one has the following equalities:
\begin{equation*}
J_{\hspace{-0.1cm}{\scriptstyle\tfrac{{\partial^2}\hspace{-0.1cm}f}{\partial x^2}}}\hspace{-0.1cm}(t,R_{t},j_{\eta}(t)) \hspace{-0.1cm}= \hspace{-0.1cm} 2 \ \tfrac{{\partial} J_{f} }{\partial u_{1}}(t,R_{t},j_{\eta}(t)) ; \hspace{0.1cm}  J_{\hspace{-0.1cm}{\scriptstyle\tfrac{{\partial}\hspace{-0.05cm}f}{\partial x}}}\hspace{-0.1cm}(t,R_{t},j_{\eta}(t)) \hspace{-0.1cm}=  \hspace{-0.1cm} \tfrac{{\partial} J_{f} }{\partial u_{2}}(t,R_{t},j_{\eta}(t)) \hspace{-0.1cm}; J_{{\scriptstyle\tfrac{{\partial}\hspace{-0.05cm}f}{\partial t}}}(t,R_{t},j_{\eta}(t)) = \hspace{-0.1cm}\tfrac{{\partial} J_{f} }{\partial t}(t,R_{t},j_{\eta}(t)).
\end{equation*}
%
%
%
\end{lem}
\vspace{-0.15cm}

{\bfseries Proof of Lemma \ref{epokep}.}
Using from one hand the equality $\tfrac{\partial \gamma}{\partial t} =  \tfrac{1}{2} \tfrac{\partial^2 \gamma}{\partial x^2}$, valid on $\R^{*}_{+}\times \R$, and, form the other hand the theorem of differentiation under the integral sign, in a neighborhood of every $(t,u_{1},u_{2})$ in $\Sigma_{a}$, provide equalities stated in Lemma \ref{epokep} on each $(t,u_{1},u_{2})$ of $\Sigma_{a}$ and then allows us to conclude.$\hfill \square$
\smallskip

Using Lemma \ref{epokep}, one gets, for every $\eta$ in $\sS(\R)$, 


\vspace{-0.5cm}

\begin{align}
\label{lopkio1}
&\cI^{\text{\tiny{(1)}}}_{\eta}:=\hspace{-0.1cm}\int^{T}_{0}   S(\tfrac{\partial^2 f}{\partial x^2}(t,G_{t}))(\eta) \ dR_{t} = \int_{\w}   J_{\tfrac{\partial^2 f}{\partial x^2}}(t,R_{t},j_{\eta}(t)) \ dR_{t} = 2 \int_{\w}   \tfrac{\partial J_{f}}{\partial u_{1}}(t,R_{t},j_{\eta}(t)) \ dR_{t}, \\
\label{lopkio2}
&\cI^{\text{\tiny{(2)}}}_{\eta}:=\int^{T}_{0}\   S(\tfrac{\partial f}{\partial x}(t,G_{t}))(\eta) \ S(W^{(G)}_{t})(\eta)\ dt = \int_{\w}\   \tfrac{\partial J_{f}}{\partial u_{2}}(t,R_{t},j_{\eta}(t)) \ dj_{\eta}(t), \\
\label{lopkio3}
&\cI^{\text{\tiny{(3)}}}_{\eta}:=\int^{T}_{0}\   S(\tfrac{\partial f}{\partial t}(t,G_{t}))(\eta)\ dt = \int_{\w}\   \tfrac{\partial J_{f}}{\partial t}(t,R_{t},j_{\eta}(t)) \ dt + \int_{\z}\   \tfrac{\partial f}{\partial t}(t,0) \ dt.
\end{align}

%
%
%
%
%
%
%
 
Thus, in order to end the proof, one just has to establish the following equality:
\vspace{-0.5cm}

\begin{align}
\label{ojzdojepfoerpozdodjz}
S(f(T,G_{T}))(\eta) - S(f(0,0))(\eta)  = \int_{\w}\   &\tfrac{\partial J_{f}}{\partial t}(t,R_{t},j_{\eta}(t)) \ dt + \int_{\w}\   \tfrac{\partial J_{f}}{\partial u_{1}}(t,R_{t},j_{\eta}(t)) \ dR_{t} \notag \\
&+ \int_{\w}\   \tfrac{\partial J_{f}}{\partial u_{2}}(t,R_{t},j_{\eta}(t)) \ dj_{\eta}(t) + \int_{\z}\   \tfrac{\partial f}{\partial t}(t,0) \ dt.
\end{align}


Since $\w$ is an open set of $[0,T]$, that does not contain $0$, it can be written under the form
\begin{equation}
\label{ocdedcdspcokdspockzpcozdodjz}
\w=\underset{i\in \N}{\bigsqcup} (a_{i}, b_{i})\ {\sqcup}\ (b,T],    
\end{equation} 
where all the intervals in \eqref{ocdedcdspcokdspockzpcozdodjz} are disjoint and where, by convention, $(x,y)=(x,y]=\emptyset$, for every reals $x$ and $y$ such that $x\geq y$. Note moreover that every element of $\{a_{i},b_{i},\ i\in\N\}$ (as well as $b$, if $(b,T]\neq\emptyset$) belongs to $\z$. We need to distinguish between two cases:
 \smallskip
  \smallskip

{\bfseries First case:} $\exists\ (a',b') \in {(0,T)}^{2}$ with $a'<b'$ \textit{s.t.} $(0,a')$ and $(b',T]$ are both subsets of $\Gamma_{R}$.
\vspace{0.1cm}

Define $a:=\sup\{a'\in [0,T], \ \text{s.t.} \ (0,a')\subset \w \}$ and $b:=\inf\{b' \in [0,T], \ \text{s.t.} \ (b',T] \subset \w \}$. Even if one has to consider a 
subset $I$ of $\N$, one may assume, and we will in the sequel, that $
(a_{i},b_{i})\neq \emptyset$, for every $i$ in $\N$. 
One can find $\rho$ in $\R^{*}_{+}$ such that $(\rho, a-\rho)\neq 
\emptyset$, $(b+\rho,T-\rho)\neq \emptyset$. Moreover, for every $i$ in $\N$, one can find $\rho_{i}$ in $\R^{*}_{+}$ such that $(a_{i}+\rho_{i},b_{i}-\rho_{i})\neq \emptyset$. 
Since all these intervals belong to $\w$, one can apply Lemma 
\ref{zpodkzpeodkzep}, with $L=J_{f}$, on each one of them. 
We then get, for 
every interval $(x,y)$ in the set of intervals $\Upsilon:= \{(\rho, a-\rho), (b+\rho,T-\rho), (a_{i}+\rho_{i},b_{i}-\rho_{i}); \ i \in \N\}$,
\vspace{-1ex}
\begin{align}
\label{paokpzok}
J_{f}&(y,R_{y}, j_{\eta}(y))-J_{f}(x,R_{x}, j_{\eta}(x)) \notag \\
\tag{$\Lambda_{x,y}$}
&= \int^{y}_{x}\   \tfrac{\partial J_{f}}{\partial t}(t,R_{t},j_{\eta}(t)) \ dt +  \int^{y}_{x}\  \tfrac{\partial J_{f}}{\partial u_{1}}(t,R_{t},j_{\eta}(t)) \ dR_{t} + \int^{y}_{x}\  \tfrac{\partial J_{f}}{\partial u_{2}}(t,R_{t},j_{\eta}(t)) \ dj_{\eta}(t).
 \end{align}

 For any interval $(x',y')$ which belongs to $\{(0,a), (b,T), (a_{i},b_{i});\  i \in \N\}$, there exists a sequence of elements ${(x_{n},y_{n})}_{n\in\N}$ in $
 \Upsilon^{\N}$ such that $ (x_{n},y_{n}) \to (x',y')$, as $n\to +\infty$ and such that $[x_n,y_n] \subset (x',y')$. Lemma 
 \ref{podiuhyiuuhiudokdfzdokfdepsoksdpokfsdpokfspokdpfoksdpfskfodkof} 
 then provides the convergence of the left hand side of $
 (\Lambda_{x_{n},y_{n}})$ to $f(y',0)-f(x',0)$, if $(x',y')$ 
 belongs to $\{(0,a), (a_{i},b_{i}); i \in \N\}$, and to  $S(f(T,G_{T}))(\eta)-f(b,
 0)$, if $(x',y') = (b,T)$. 
Besides, the Lebesgue's dominated convergence theorem applies to each integrand of  the right hand side of $(\Lambda_{x_{n},y_{n}})$,  since they are all continuous. This provides the convergence of the right hand side of $(\Lambda_{x_{n},y_{n}})$ to $ \int^{y'}_{x'}\   \tfrac{\partial J_{f}}{\partial t}(t,R_{t},j_{\eta}(t)) \ dt +  \int^{y'}_{x'}\  \tfrac{\partial J_{f}}{\partial u_{1}}(t,R_{t},j_{\eta}(t)) \ dR_{t} + \int^{y'}_{x'}\  \tfrac{\partial J_{f}}{\partial u_{2}}(t,R_{t},j_{\eta}(t)) \ dj_{\eta}(t)$, for any $(x',y')$ in $\{(0,a), (b,T), (a_{i},b_{i});\  i \in \N\}$.
In view of \eqref{eoijecjaaa} and \eqref{hiuiuiuiuiuddpzofkeprofkeprofkerpfkpok}, and making the summation of $(\Lambda_{a_{i},b_{i}})$, over all $i\in \N$, we then get:  
\vspace{-0.5cm}
  
\begin{align}
\label{zdozeodkzepdokezzpezdedokzede}
S(f&(T,G_{T}))(\eta) - S(f(0,0))(\eta) -( f(b,0)- f(a,0) -\sum_{i\in \N} (f(b_{i},0)-f(a_{i},0)))  \notag \\
&= \int_{\w}\   \tfrac{\partial J_{f}}{\partial t}(t,R_{t},j_{\eta}(t)) \ dt + \int_{\w}\   \tfrac{\partial J_{f}}{\partial u_{1}}(t,R_{t},j_{\eta}(t)) \ dR_{t} + \int_{\w}\  \tfrac{\partial J_{f}}{\partial u_{2}}(t,R_{t},j_{\eta}(t)) \ dj_{\eta}(t).
\end{align}

Denote $\Delta:= \int_{\z}\   \tfrac{\partial f}{\partial t}(t,0) \ dt$,  one has the equality:

\vspace{-3ex}

\begin{align}
\label{zdozeodkzepdokezzpokzede}
\Delta&=  \int_{[0,T]}\   \tfrac{\partial f}{\partial t}(t,0) \ dt
 - \int_{\w}\   \tfrac{\partial f}{\partial t}(t,0) \ dt 
=  \int_{[0,T] \backslash (0,a]\sqcup (b,T]}\    \tfrac{\partial f}{\partial t}(t,0) \ dt  - \sum_{i\in \N} \int^{b_{i}}_{a_{i}}\    \tfrac{\partial f}{\partial t}(t,0) \ dt \notag \\
&= f(b,0)-f(a,0) - \sum_{i\in \N} (f(b_{i},0)-f(a_{i},0)).
\end{align}
Using \eqref{zdozeodkzepdokezzpokzede}, Equality \eqref{zdozeodkzepdokezzpezdedokzede} then reads:
%
%
%
%
\vspace{-2ex}
\begin{align*}
S(f&(T,G_{T}))(\eta) - S(f(0,0))(\eta) - \int_{\z}\   \tfrac{\partial f}{\partial t}(t,0) \ dt \notag \\
&= \int_{\w}\   \tfrac{\partial J_{f}}{\partial t}(t,R_{t},j_{\eta}(t)) \ dt + \int_{\w}\   \tfrac{\partial J_{f}}{\partial u_{1}}(t,R_{t},j_{\eta}(t)) \ dR_{t} + \int_{\w}\  \tfrac{\partial J_{f}}{\partial u_{2}}(t,R_{t},j_{\eta}(t)) \ dj_{\eta}(t),
\end{align*}

which is nothing but \eqref{ojzdojepfoerpozdodjz} and therfore ends the proof in this case.
\smallskip
 
 {\bfseries Second case:} There is no $(a',b')$ in ${(0,T)}^{2}$ \hspace{-0.1cm} with $a'\hspace{-0.1cm}<\hspace{-0.1cm}b'$ 
\textit{s.t.} both $(0,a')$ \& $(b',T]$ are subsets of $\Gamma_{R}$.

  Since the cases of $0$ and $T$ can be treated in the same manner, we 
  only treat here the case of $T$. We then assume that there is no $b'$ in 
  $(0,T)$ such that $(b',T]\subset \w$. 
 We need to distinguish between two cases. If there exists $b'$ in $(0,T)$ 
such that $(b',T] \subset \z$ then the problem can be reduced to establish \eqref{ojzdojzdodjz} between $0$ and $\hat{b}$, where  $\hat{b}:=\inf\{b\in  (0,T); \ [b',T]\subset \z \}$. Otherwise, one can find an increasing sequence ${(T_{n})}_{n\in\N}$ of $\w^{\N}$, which converge to $T$. For every $n$ in $\N$, denote $(a^{(n)}_{i},b^{(n)}_{i})$ the interval $(a_{i}, b_{i})$ of $\w$ which contains $T_{n}$. For every integer $n$, $(a^{(n)}_{i},T_{n}]$ is a non empty subset of $\w$. Therefore, one can use first case to establish \eqref{ojzdojzdodjz} between $0$ and $T_{n}$. To establish the equality of S-transform of both sides of \eqref{ojzdojzdodjz}, (between $0$ and $T_{n}$), it then remains to apply, from one hand Lemma \ref{podiuhyiuuhiudokdfzdokfdepsoksdpokfsdpokfspokdpfoksdpfskfodkof} to $J_{f}(T_{n},R_{T_{n}}, j_{\eta}(T_{n}))$ and, form the other hand, Lebesgue's dominated convergence theorem to the following integrals:

 \vspace{-3ex}

\begin{align*}
&\cI^{\text{\tiny{(1)}}}_{\eta,n}:=\hspace{-0.1cm}\int^{T_{n}}_{0}   S(\tfrac{\partial^2 f}{\partial x^2}(t,G_{t}))(\eta) \ dR_{t} = 2 \int_{\w^{(n)}}   \tfrac{\partial J_{f}}{\partial u_{1}}(t,R_{t},j_{\eta}(t)) \ dR_{t}, \\
&\cI^{\text{\tiny{(2)}}}_{\eta,n}:=\int^{T_{n}}_{0}\   S(\tfrac{\partial f}{\partial x}(t,G_{t}))(\eta) \ S(W^{(G)}_{t})(\eta)\ dt = \int_{\w^{(n)}}\   \tfrac{\partial J_{f}}{\partial u_{2}}(t,R_{t},j_{\eta}(t)) \ dj_{\eta}(t), \\
&\cI^{\text{\tiny{(3)}}}_{\eta,n}:=\int^{T_{n}}_{0}\   S(\tfrac{\partial f}{\partial t}(t,G_{t}))(\eta)\ dt = \int_{\w^{(n)}}\   \tfrac{\partial J_{f}}{\partial t}(t,R_{t},j_{\eta}(t)) \ dt + \int_{ {\mathcal{Z}^{\scriptscriptstyle \hspace{-0.01cm} T,(n)}_{\hspace{-0.05cm}\scriptscriptstyle R}}}\   \tfrac{\partial f}{\partial t}(t,0) \ dt,
\end{align*}

where $\w^{(n)}$ denotes  $ \w \cap [0,T_{n}]$ and
${\mathcal{Z}^{\scriptscriptstyle \hspace{-0.01cm} T,(n)}_{\hspace{-0.05cm}\scriptscriptstyle R}}$ denotes $ \z \cap [0,T_{n}]$. This result and the fact that $J_{f}$ is a $C^{1}$ function on $\Sigma_{a}$ allows us to apply \cite[Theorem $8.6$]{Kuo2} and thus to conclude.
\end{pr}

\newpage
\subsection{Comparison with other Itô formulas for Gaussian processes}
\label{pokprofkdeezdzedezdedee}


Since \cite{nualart}, many Itô formula for Gaussian 
processes have been established. If one excepts Itô formula for 
Gaussian semimartingales, that are well known, all the Itô formulas provided, for Gaussian processes in general, in the 
literature of \textit{functional extensions} so far, namely: \cite[Theorems $1$  \& $2$]{nualart}, \cite[Theorem $31$]{MV05}, \cite[Theorem $1$]{NuTa06}, \cite[Corollary $8.13$]{KRT07}, 
\cite[Proposition 11.7]{KrRu10} and \cite[Theorem $3.2$]{LN12}, 
are established using the divergence type integral.  
A requirement of all these previous references
is that the variance function $t\mapsto R_{t}$ is, at least, continuous and with bounded variations on $[0,T]$. Assuming the continuity of $R$ seems reasonable. Indeed, otherwise, as the anonymous referee noticed, Equality \eqref{ojzdojzdodjz} may fail
for very simple functions $f$, such as $f(x):=x^{4}$.
If one excepts \cite[Corollary $8.13$]{KRT07}, another requirements of theses references above
 is that the function $f$ is of class 
$C^{2}$ and, together with all its derivatives,  with sub-exponential 
growth (\ie fulfills \eqref{orijv}). 
In view of this fact, it appears that the assumptions 
made in Theorem \ref{ergoigfjerfijoijgfoifjgoidfjgofd} are minimal. However,  
to see to what extent Theorem \ref{ergoigfjerfijoijgfoifjgoidfjgofd} 
generalizes Itô formulas for Gaussian processes that already exist, 
let us make a detailed comparison.

\textit{\bfseries Comparison with the conditions on function $f$}

The function $f$ (\ie $f(t,x):=f(x)$) is assumed to be of class $C^{\infty}$ in \cite[Theorem $31$]{MV05} and in \cite[Proposition 
11.7]{KrRu10}, and of class $C^{7}$ in\cite[Theorem $1$]{NuTa06}. In \cite[Corollary $8.13$]{KRT07} and $f$ is assumed to be of class 
$C^{2}$ but not with sub-exponential growth; the second derivative of $f
$ therein is assumed to be bounded. However, since the stochastic calculus for 
Gaussian processes, developed in \cite{KRT07} requires that the 
covariance function has a planar bounded variation, which corresponds 
to ``regular'' processes (such as fBm for $H>1/2$), one easily sees that 
the price to pay for relaxing the assumption on the growth of $f$ is that 
one can not deal with irregular Gaussian processes
(that is precisely to 
overcome this limitation on the regularity of $G$ that \cite{KrRu10} has been written. However, as we stated above, this latter reference requires much more than the growth condition we make on $f$).
\smallskip

\textit{\bfseries Comparison with the assumptions made on $G$}
\vspace{-1ex}
\bit
\item Comparison with the assumptions made on $R$
\eit
\vspace{-1ex}
The variance function $t\mapsto R_{t}$ is assumed to:
be continuous and of bounded variations 
on $[0,T]$ in \cite[Corollary $8.13$]{KRT07}, 
 be of class $C^{2}$ 
 on $\R^{*}
_{+}$ in \cite[Theorem $31$]{MV05}, fulfill Assumptions $(3)$ et $(4)$ in 
\cite[Theorem $1$]{NuTa06}, fulfill Assumptions
 (A), (B) and (C)  in \cite[Proposition 11.7]
 {KrRu10} and, in \cite[Theorem $3.2$]{LN12}, to verify the two following conditions:
\vspace{-1ex}
 \bit
 \iti for every $s$ in $[0,T]$, the map $t \mapsto R(t,s)$ is absolutely continuous on $[0,T]$;
 \itii there exists $\alpha>1$ such that: $ \underset{t\in[0,T]}{\sup}\ \int^{T}_{0} \ {\big|\frac{\partial R}{\partial s} (s,t)\big|}^{\alpha} \ ds <+\infty$.
 \eit
\vspace{-2ex}
  \vspace{-1ex}
\bit
\item Other assumptions 
\eit
\vspace{-1ex}
In \cite[Theorems $1$  \& $2$]{nualart}, in addition to the assumptions made in Theorem \ref{ergoigfjerfijoijgfoifjgoidfjgofd},  the kernel $K$ has to fulfill 
Assumptions (K1) to (K3), in the singular case, and (K1) to (K4), in the 
regular case. Other assumptions on the process $G$ are difficult to 
compare in general. As it is stated in Remark 
\ref{ozeufhioeruvheirezzezeedzedzedzdu}, a key propoerty in our 
construction of integral is that the maps $t\mapsto \E[G_{t} \int_{\R} \eta (s) \ dB_{s}]$ are absolutely continuous w.r.t. the Lebesgue measure for 
every $\eta$ in $\sS(\R)$.   A related assumption in other papers is that the functions $t\mapsto \E[G_{t} G_{s}]$ or $t\mapsto \E[G_{t} B_{s}]$
  are absolutely continuous or of
bounded variations for every $s$, see \cite{LN12,nualart}. However These
assumptions do not, in general, imply each other. In an another class of papers, namely \cite{KRT07}, some quadratic variation type conditions are imposed on $G$, e.g. in \cite{KRT07} or in \cite{NuTa06}. These assumptions are even more
difficult to compare with the present setting and would lead us too far from the goal of this present work. We therefore postpone a more detailed comparison to a future work.
\smallskip

Note also that he Itô formula provided in \cite{LN12} is extended in this work, while the other results presented in \cite{LN12}  are extended in \cite{JL17-2}.

In view of the arguments developed above, it appears that the Itô formula we present here offers improvements on the ones presented in \cite{nualart,cn05,MV05,NuTa06,KrRu10, LN12}, by allowing one to 
have less restrictive hypotheses. 
 Of course all the Gaussian processes 
in\hspace{0.05cm} $\mathscr{G}$ of ``reference'' fulfill assumptions of 
Theorem \ref{ergoigfjerfijoijgfoifjgoidfjgofd}. In the case of mBm one needs to assume that $h
$ is a $C^{1}$ function with its derivative 
bounded on $\rR_{D}$. Note moreover that, applying Theorem 
\ref{ergoigfjerfijoijgfoifjgoidfjgofd} when $G$ is a fBm \textit{(resp. a mBm)} 
allows one to recover \cite[Theorem 4.1 \& Rk. 4.6]{ben1} \textit{(resp.  
\cite[Theorem $5.5$]{JLJLV1})}. When $G$ is a $\mathscr{V}_{\gamma}$ - 
process one recovers and  extends, as we showed above, \cite[Theorem 
31]{MV05}.

\vspace{-2ex}

\section{Comparison with other stochastic integrals}
\label{ozfjozrio}

\textit{\bfseries Forewords}

%
In order to define 
the divergence integral of a continuous Gaussian process $G$ in the way of \cite{nualart} and then of \cite{MV05},  it is essential to first know a  
representation of $G$ on a compact set of the form $[0,T]$. In 
general, \cite[Theorem $4.1$]{HidHit} ensures that any Gaussian 
process may be written as a sum of two terms; one of them being $
\sum_{i=1}^{N} \int^t_0 K_i(t,u)\ dW_i (u)$, where $N$ is a positive integer (possibly 
infinite) and $W$ is a Brownian motion. However it is not an easy task 
to obtain such a decomposition for a given process $G$. For instance, 
although a kernel is known for fBm, this is not the case of bifractional 
motion \cite{HV03}.
Likewise, writing the moving average and harmonizable multifractional 
Brownian motion under this form remains an open problem (see 
\cite[Section 5]{LLVH} for more details). Moreover, Gaussian bridges 
in general are an example of Gaussian processes which do not admit 
``proper'' Volterra representation, \ie that can not be written under the form 
\eqref{erofkoredefefekdeddedep} (see \cite[Ex. 3.3]{SoVi14}). 
Thus, there is no hope to use \cite{nualart} nor \cite{MV05} in order to build a stochastic integral \textit{wrt} Gaussian processes of the form \eqref{erofkorekp}, for which  one does not know any integral representation on a compact set included in $[0,T]$. To overcome this deficiency 
one then might consider \cite{NuTa06,KRT07,KrRu10,LN12,SoVi14}.
As we stated above, in these latter references 
one needs that the covariance function fulfills some requirements. 
 However, it happens sometimes that one has to deal with Gaussian 
processes, given under the form \eqref{erofkorekp}, for which one does 
not know how to compute the covariance function, such as the one where 
$g_{t}$ is defined by setting $g_{t}(u) = \i1_{[0,t]}(u)\ K_{h(t)}(t,u)$, where 
$h:[0,T]\rightarrow (0,1)$ is a continuous deterministic function and where  
the family of Kernel ${(K_{H})}_{H\in(0,1)}$ is the one defined in 
\cite[(5.8)]{Nu} in the case where $H\in(1/2,1)$ and in \cite[Proposition 
5.1.3]{Nu} 
in the case where  $H\in(0,1/2)$. 
As a consequence our stochastic calculus it is the only one
 available when the 
Gaussian process $G$ can be written under the form \eqref{erofkorekp} but not under 
any of the form \eqref{erofkoredefefekdeddedep} nor  \eqref{ooijoierfkfpopfokerfpoefoierfu};  or when the stochastic calculus provided in  \cite{KRT07,KrRu10, LN12} does not apply. 
The work  provided in \cite{KRT07} offers an alternative to the previous methods 
to build a stochastic integral \textit{wrt} continuous Gaussian
processes, for which one 
knows the covariance function. Introducing the concept of covariance 
measure structure, the authors built and developed a stochastic calculus 
\textit{wrt} ``regular'' processes (such as fBm with $H\geq 1/2$). This work has been 
extended to the ``singular'' case in \cite{KrRu10}. This approach is particularly suitable when 
the kernel is not explicitly known, under any of the representations 
 \eqref{erofkorekp} to \eqref{ooijoierfkfpopfokerfpoefoierfu}, (like in the case 
of bifractional motion \cite{HV03}). However, 
the Itô formula in \cite{KrRu10} is quite restrictive.
The conditions required in \cite{KRT07} are not so 
restrictive but they do not allow one to deal with ``irregular'' Gaussian processes, by the very 
essence of covariance measure structure.

In this section we make first, in Subsection \ref{orefij}, a comparison of the Wick-Itô 
stochastic integral we developed above with the
\textit{functional extensions} of stochastic integrals developed
 in \cite{nualart,MV05} and then,  in Subection \ref{oervoerveiorj},  with the Itô 
integral. 

\vspace{-2ex}

\subsection{Comparison with Malliavin Calculus or divergence type integrals}
\label{orefij}

We start by making the comparison between our Wick-Itô integral and 
the divergence type integral developed in \cite{nualart}. We will 
then show, in Section \ref{zeofkepofkpo123879987},  
that the Wick-Itô integral fully generalizes the (extended) Skorohod 
integral developed in \cite{MV05}. Let $T>0$ be fixed and let us take $
\rR=[0,T]$. Let $G:={(G_t)}_{t\in [0,T]}$ be a Volterra process.

\vspace{-1ex}

\subsubsection{Comparison with divergence type integral of \cite{nualart}}
\label{zeofkepofkpo123879}

The goal of this section is to compare the Wick-Itô integral \textit{wrt} $G$  to the divergence integral \textit{wrt} $G$, defined in \cite{nualart} and in \cite{MV05} and studied in \cite{nualart,Nu05} and in \cite{MV05}. In  \cite{nualart} $G$ is a assumed to be a continuous process while it is not assumed to be continuous in \cite{MV05}. One therefore will assume (in Subsection \ref{zeofkepofkpo123879}
 only) that $G$ is continuous on $[0,T]$. $G$ being a Volterra process, it can be written, for any real $t$ in $[0,T]$, 
\begin{equation*}
G_t = \int_{0}^t K(t,s)\ dW_s,
\end{equation*}
where the kernel $K(t,s)$, defined on ${[0,T]}^{2}$, is such that $K(t,s)=0$ 
on the set ${[0,T]}^{2}\backslash\{(u,v)\in (0,T] \times [0,T]: \ v< u \}
$   and verifies for any $t \geq 0$,  $\widehat{K_{t}}:= \int_0^t K(t,s)^2 ds < \infty$.

Denote $L^2(\Omega,L^2([0,T]))$ the set of 
random process $u$ such that $ \left\|  u \right\|
^2_{L^2(\Omega,L^2([0,T]))} \hspace{-0.1cm}:=  \E[{\textstyle{\int^T_0}} u^2_t\ dt]< 
+\infty$. The main result of this section is Theorem \ref{oerferoifjeoifjefi},  
which 
states that every process $u$ which belongs to $L^2(\Omega,L^2([0,T]))$ and that belongs
to the domain of the divergence of $G$  is also Wick-Itô integrable \textit{wrt} $G$, on $[0,T]$. Moreover, one has the equality  $\int^T_0 u_s\ \delta G_{s} = \int^T_0 u_s \ dG_{s} $, where $\int^T_0 u_s\ \delta G_{s}$ \hspace{-0.1cm} denotes \hspace{-0.1cm} the divergence integral \hspace{-0.1cm} on $[0,T]$, associated to $G$, that has been defined in\cite{nualart}.
 In order to state rigorously this result we briefly recall some elements and notations of stochastic calculus of variations {\it wrt} $G$ (for a presentation of Malliavin
calculus, see {\it e.g.} \cite{bally,Nu}), as well as  the approach of \cite{nualart} and \cite{Nu05} for the construction of a stochastic integral \textit{wrt} to Volterra processes. The real $T>0$ being fixed, one still note $G$ the process ${(G_t)}_{t\in[0,T]}$ since there is no risk of confusion. $G$ being a centered Gaussian process, denote $\cH_T$ the reproducing kernel Hilbert space (R.K.H.S.) 
defined as the closure of the set $\cE_T:= {\text{span}}\{\i1_{[0,t]}, \ t\in [0,T]\}$, with respect to  the inner product ${<,>}_{\cH_T}$, that has 
been defined by setting ${<\i1_{[0,t]},\i1_{[0,s]}>}_{\cH_T} := R_{t,s}$. Denote $H_1$ the first Wiener 
chaos of $G$ and $G(\varphi)$ the image in $H_1$ of an element $\varphi$ of $\cH_T$ by the isometry, between  $\cH_T$ and $H_1$, that associates $\i1_{[0,t]}$ to $G_t$.

\begin{rem}
It is not always true that the bilinear form ${< , >}_{\cH_T}$ defined by  ${<
\i1_{[0,t]},\i1_{[0,s]}>}_{\cH_T} := R_{t,s}$ is an inner product. For 
example, for the Brownian bridge  $\widehat{B}:={(\widehat{B}_t)}_{t
\in[0,1]}$  on $[0,1]$, one has ${\|\i1_{[0,1]}\|}_{\cH_T}=0$. 
For this reason we will assume in the sequel that ${< , >}_{\cH_T}$, 
defined above, is really an inner product. The reader interested in details 
on Reproducing Kernels Hilbert Spaces may refer to \cite[Chap.8]
{Jan97} as well as to \cite[Appendix B]{JLJLV1} in the case of mBm.
\end{rem}

Define $\mathcal S := \left\{
V=f\left(G(\varphi_1),G(\varphi_2), \ldots, G(\varphi_n)\right),
f \in C^{\infty}_b(\R^n), \varphi_i \in \cH_T, i=1,\ldots, n \right\}$.
%
 For an element $V$ of $\mathcal S$, one defines the derivative operator $D^G$ as:
\vspace{-0.25cm}
 
$$D^G V := \sum_{i=1}^n \tfrac{\partial f}{\partial x_i}  \left(G(\varphi_1),G(\varphi_2), \ldots,
G(\varphi_n)\right) \varphi_n.$$

The derivative operator $D^G$ is a closable unbounded operator from $L^2(\Omega)$ into $L^2(\Omega; \cH_T)$. We note $\bD_G$ the closure of $\cS$ with respect to the norm defined by $ \left\|  V \right\|_{G,1,2} := (\E[V^2] + \E[\|D^G V\|_{L^2(\Omega;\cH_T)}^2]) ^\frac{1}{2}$.
%
We denote by $\delta_G$, and call divergence integral with respect to $G$, the adjoint of the derivative operator $D^G$. The domain of $\delta^G$ in $L^2$,  denoted  $\Dom(\delta^G)$, is the set of the elements $u$ in $L^2(\Omega; \cH_T)$ such that there exists a constant $c$ verifying, for all $V$ in $\cS$, $|\E(<D^GV,u>_{\cH_T})| \leq c \ \|V\|_2$,
%
%
%
 where $\| \ \|_2$ denotes the norm in $L^2(\Omega)$. If $u$ belongs to $\Dom(\delta^G)$, $\delta^G(u)$ is the element of $L^2(\Omega)$ defined by the duality relationship: $\E(V\delta^G(u)) = \E(<D^GV,u>_{\cH_T}), {\text{  for all }} V  {\text{  in }} \bD_G$.
%
We will simply denote, in the sequel, $D, \bD, \delta$ and $\| \ \|_{1,2}$ when $G$ is a Brownian motion.
%
%
%
%
%
%
%
%
%
%
%
%
%
%
Define now  the linear operator $K_{-}:\cE_T\rightarrow  L^2([0,T])$ by $K_{-}(\i1_{{[0,t]}}) :=K(t,.)$ and denote ${\|\ \|}_{\cH_T}$ the norm on $\cH_T$ which derives from the inner product ${<, >}_{\cH_T}$. Since ${\|\varphi \|}_{\cH_T} = {\|K_{-}(\varphi) \|}_{L^2([0,T])}$, for every $\varphi$ in $\cE_T$, it is clear that the operator $K_{-}$ can be extended to a linear isometry, still denoted $K_{-}$, between  $(\cH_{T}, {\|\ \|}_{\cH_T})$ and a closed subset of $L^2([0,T])$. Besides, one can show, \cite[(12)]{nualart}, that $\Dom(\delta_G) = (K_{-})^{-1} (\Dom(\delta))$.
%
%
%
%
%
%
Moreover, for a process $v$ in $\Dom(\delta_G)$ one has:
\begin{equation}
\label{oeijrfoeirfje}
\delta_G(v)= \int_0^T (K_{-} v)(s) \ \delta W_s. 
\end{equation}
In other words, $\delta_G(v)$, the divergence integral of $v$ \textit{wrt}  $G$, also noted $\int_0^T v(s)\ \delta G_s$, verifies the equality  $\int_0^T v(s)\ \delta G_s = \int_0^T (K_{-} v)(s) \ \delta W_s$. 
In order to prove Theorem \ref{oerferoifjeoifjefi} below, one needs to define the adjoint of the operator $K_{-}$, that we will denote  $K_{+}$, not only on the set $\cE_T$ but also on $\sS(\R)$. For this reason we recall the two following hypotheses, given in \cite{Nu05} for fBm, that we will make in the sequel on the kernel $K(t,s)$.
\vspace{-1ex}
\begin{itemize}
\ia  $K(t,s)$ is continuously differentiable on $\{0<s<t\leq T\}$ and its partial derivative verify the following integrability condition:
\vspace{-1ex}
$$ \underset{\varepsilon \leq t \leq T}{\sup} \int^T_{t}  |\tfrac{\partial K}{\partial r}(r,t)|  (r-t) \ dr +  \int^t_{0}  |\tfrac{\partial K}{\partial t}(t,s)|  (t-s) \ ds <\infty,$$
\vspace{-1ex}

for any $\varepsilon$ in $(0,T)$. Moreover, $t\mapsto \int^{t\wedge b}_{0}  \tfrac{\partial K}{\partial t}(t,s)  (t\wedge b-s\vee a)_{+} \ ds$ is continuous on $(0,T]$, for all $0\leq a \leq b$. 

\iaa  The function $k(t):= \int^{t}_{0} K(t,s) \ ds$ is continuously differentiable on $(0,T]$.
\end{itemize}
\vspace{-1ex}
We present here the arguments given in \cite[Section 2]{Nu05} for fBm about the operator $K_{+}$, but in a slightly different manner. Denote $C^{1}_{b}(\R)$ the set of differentiable functions which are bounded together with its derivatives. Hypotheses H1) and H2) allow us to define the operator $K_{+}$ on $ \cE_T\cup C^{1}_{b}(\R)$ by setting, for every $t$ in $[0,T]$, $(K_{+}\varphi)(t) := k'(t) \ \varphi(t)  + \int_0^t  \tfrac{\partial K}{\partial t}(t,r)  \ (\varphi(r)-\varphi(t)) \ dr$.

In view of \cite[p.116]{Nu05}, it is easy to check that we have, for any $(\psi, \varphi)$ in $\cE_T\times\cE_T$, the equality 
\begin{equation}
\label{oefjeroifjeorfjoeri}
{<K_{+}(\varphi),\psi>}_{L^2([0,T])} \ =  \  {<\varphi,K_{-}(\psi)>}_{L^2([0,T])}.
\end{equation}
It is clear that one has, in this section, $g_t:=K_{-}(\i1_{[0,t]})$, for every $t$ in $[0,T]$. It is established in \cite[Propostion 2]{Nu05} that $g'_t$ exists and that $g'_t=K_{+}(.)(t)$ for every $t$ in $(0,T]$. However it is not possible to establish that $t\mapsto W^{(G)}_t$ is $(\cS^*)$-integrable on $[0,T]$ without any additional assumption. Moreover one needs to be able to establish that $\int^T_0 u_s \ dG_{s}$ exist for a reasonable class of processes $u$.
Thus, following  \cite[Proposition 7]{Nu05}, we will assume in the sequel the following condition:
\begin{itemize}
\label{peokepork}
\iaaa The function $C:t\mapsto |k'(t)| +  \int_0^t  |\tfrac{\partial K}{\partial t}(t,r)| \ (t-r)  \ dr$ belongs to $L^2([0,T])$.
\end{itemize}

\begin{rem}
It is clear that $H_1)$, $H_2)$ and $H_3)$ entail that Assumptions \Dai   and \Daii hold. We will show, in the next subsection (Remark \ref{peorfkerpofkeprofkepr}), that they are not always necessary. 
\end{rem}

%
%

The following result will be useful in the proof of Theorem \ref{oerferoifjeoifjefi} below.

\begin{lem}
\label{oeirvjoerivjoeivjeo}
If Assumptions $H_1)$, $H_2)$ and $H_3)$ hold, any process $u$ in $L^2(\Omega,L^2([0,T]))$ is Wick-Itô integrable with respect to $G$.
\end{lem}

\begin{pr}
The proof of this lemma, which consists on verifying that Condition \iIp is verified with $p=q=2$, can be found in \cite[Proposition 7]{Nu05}.
\end{pr}

Since one has: $|(K_{+}\eta)(t)|  \leq   |k'(t)| {\|\eta\|}_{\infty} + {\|\eta'\|}_{\infty}   \int_0^t |\tfrac{\partial K}{\partial t}(t,r)| \ (t-r) \ dr$, for every $(\eta,t)$ in $\sS(\R)\times [0,T]$,
%
%
Hypothesis H3) implies in particular that  $K_{+}(\eta)$ belongs to $L^2([0,T])$. Note moreover that, for every $\eta$ in $\sS(\R)$,
\vspace{-2ex}

\begin{equation}
\label{oziedzeoidhzeoihzeoih}
K_{+}(\eta)(t) = {<\delta_t, K_{+}(\eta)>} = \tfrac{d}{dt} {<K_{-}(\i1_{[0,t]}),\eta>}_{L^2([0,T])}.
\end{equation}


%
%
%
%
%
%
%
%
%

%
%
%
%
%
%
%
%
%
%
%
%

%
%
 The following result, which is a 
 consequence of results given in \cite[Section 2 and Proposition 7]{Nu05}, will be essential in order to prove Theorem \ref{oerferoifjeoifjefi} below. Denote  $\overline{\mathcal{E}_{T}}^{{\|\ \|}_{T}}$
 the closure of the set $\cE_T$ with respect to the norm ${\|\ \|}_{T}:= {\|\ \|}_{\cH_T}+{\|\ \|}_{L^{2}([0,T])}$. Note that $\overline{\mathcal{E}_{T}}^{{\|\ \|}_{T}} \subset \cH_T \cap L^2([0,T])$.
 \begin{lem}
\label{edzqsqsqdezdzedz}
 For any function $\psi$ in $\overline{\mathcal{E}_{T}}^{{\|\ \|}_{T}}$ and $\eta$ in $\sS(\R)$, one has:
\begin{equation}
\label{peofpoferpoferpofepfoj}
 {<K_{+}(\eta),\psi>}_{L^2([0,T])} \ =  \  {<\eta,K_{-}(\psi)>}_{L^2([0,T])}.
\end{equation}
\end{lem}

\begin{pr}

It is easy to check \eqref{peofpoferpoferpofepfoj} directly in the case where $\psi$ 
is in $\cE_T$ and  $\eta$ in $\sS(\R)$, using \eqref{oziedzeoidhzeoihzeoih}. The 
fact that, for every $\psi\in\cE_T$, ${<\psi,K_{+}(.)>}_{L^2([0,T])}$ belongs to $
\sS_{\hspace{-0.15cm}-p}(\R)$, for every $p$ in $\bN^{*}$ is also clear. Thus, for 
every $p$ in $\bN^{*}$, one easily sees that the map $\Psi_p:\psi\mapsto  {<\psi,K_{+}
(.)>}_{L^2([0,T])}$ is uniformly continuous from $(\cE_T,{\|\  \  \|}_{\cH_T})$ to  $
(\sS_{\hspace{-0.15cm}-p}(\R), {|\  \  |}_{-p})$ and can then be extended uniquely to 
$\overline{\mathcal{E}_{T}}^{{\|\ \|}_{T}}$ (we will denote $\Psi_p:=<\psi,K_{+}
(.)>$ this extension). The same 
argument can be applied to the map $\Phi_p:\psi\mapsto  {< . ,K_{-}(\psi)>}
_{L^2([0,T])}$. The equality of $\Psi_p$ and $\Phi_p$ on $\overline{\mathcal{E}_{T}}^{{\|\ \|}_{T}}$ from one hand, and the fact that ${<\psi,K_{+}(.)>} = {<\psi,K_{+}(.)>}_{L^2([0,T])}$ for any $\psi$ in $\overline{\mathcal{E}_{T}}^{{\|\ \|}_{T}}$ from the other hand allow us to 
conclude. \end{pr}

The main result of this section is the following.
\begin{theo}
\label{oerferoifjeoifjefi}
Assume that $H_1)$, $H_2)$ and $H_3)$ hold. Let $u$ be a process in $L^2(\Omega,\overline{\mathcal{E}_{T}}^{{\|\ \|}_{T}})$, then 
$u$ belongs to the domain of the divergence of $G$, and $u$ is Wick-Itô integrable on $[0,T]$  \textit{wrt}  $G$. Moreover one has the equality
\begin{equation}
\label{dzepidzeàidjze}
 \int^T_0 u_s\ \delta G_{s} = \int^T_0 u_s \ \di G_{s}.
\end{equation}
\end{theo}

%
%

\begin{pr}
 The proof we give here is a generalization, to Volterra processes, of the proof provided, in the particular case of fBm, in \cite[Proposition 8]{Nu05}. We however write it down here for reader's convenience. The fact that  $\int_0^T u_s\ \di G_s$ is well-defined has been established in Lemma \ref{oeirvjoerivjoeivjeo}. Besides, for every fixed $\eta$ in $\sS(\R)$, one  has: $ \mathcal{L}_1:= S(\int^T_{0}  u_s\ \delta G_{s})(\eta) =  S(\int_0^T (K_{-} u)(s) \ \delta W_s)(\eta) = \int_0^T S[(K_{-} u)(s)](\eta)\  \eta(s) \ ds$.
Note that the last equality results from the fact that the Wick-Itô integral  \textit{wrt}  Brownian motion generalizes the Hitsuda-Skorohod
integral (see for example \cite[Theorem 2.5.9]{HOUZ} or \cite[(13.8)]{Kuo2}). Using the previous equality, Fubini's theorem and Lemma \ref{edzqsqsqdezdzedz}, one gets:
\vspace{-0.5cm}

\begin{align*}
 \mathcal{L}_1 &=  \int_0^T  \E[(K_{-} u)(s) :e^{<.,\eta>}: ] \ \eta(s)\ ds =  \E[   :e^{<.,\eta>}:  \  <K_{-}u,\eta>_{L^2([0,T])}] \notag \\
&=  \E[   :e^{<.,\eta>}:   \ <u,K_+(\eta)>_{L^2([0,T])}] =  \int_0^T  \E[u_s\ :e^{<.,\eta>}: ] \ K_+(\eta)(s)\ ds.
\end{align*}


It then remains to use $(ii)$ of Theorem \ref{tardileomalet} as well as \eqref{oziedzeoidhzeoihzeoih} to obtain:
\vspace{-0.15cm}

\begin{equation*} 
 \mathcal{L}_1 = \int_0^T  S(u_s)(\eta)\ S(W^{(G)}_s)(\eta)\ ds =  S(\int_0^T u_s\ \diamond W^{(G)}_s\ ds)(\eta).
\end{equation*}

We hence have shown, for every $\eta$ in $\sS(\R)$, the equality $S(\int^T_{0}  u_s\ \delta G_{s})(\eta) = S(\int_0^T u_s\ \di G_s)(\eta)$.


The injectivity of $S$-transform (see $(i)$ of Lemma \ref{dkdskcsdckksdksdmksdmlkskdm}) allows us to conclude.
\end{pr}


\begin{ex}[The case of fBm]
\label{zopdijzeoidzeoidjzeoidjzeoidjz}
Let $T>0$. For any $H$ in $(0,1)$, define $B^H_t:= \int_{0}^t K_{H}
(t,s)\ dW_s$, where the kernel $K_H$ is defined in \cite[$(13)$]
{Nu05}. $B^H$ is an fBm of Hurst index $H$. 
Moreover the process $B^H$ fulfills
 H1), H2) and H3). 
This implies in particular that the Wick-Itô integral $\int_{0}^T B^H_s \ 
\di B^H_s$ exists, for any $H$ in $(0,1)$. We know moreover, thanks to 
Example \ref{zeodijze}, that the equality $\int_{0}^T B^H_s \ \di B^H_s 
\stackrel{\text a.s.}{=} \tfrac{1}{2} \ ({(B^H_T)}^{2} - T^{2H})$ is true for 
any $H$ in $(0,1)$.
On the other hand, the divergence integral \textit{wrt}  $B^H$ is only 
defined and developed in \cite[Section 8]{nualart} or in \cite[Section 6]{Nu} for $H>1/4$, as we mentioned in the introduction. One moreover knows  that $B^{H}$ does not belong to $\Dom(\delta_{B^{H}})$ when $H<1/4$ and one has to use the extended divergence integral \textit{wrt} fBm developed in \cite{cn05}.

\end{ex}

\begin{coro}
The set  ${\{ {\textstyle{\int^T_{0}}} f(s) \  \delta G_s, \   f  \in \cH_T \}}$ of Wiener divergence integral  \textit{wrt}  $G$, coincide with the set $\Theta_G:= \overline{  \{ {\textstyle{\int^T_{0}}} f(s) \  \di G_s, \   f  \in {\cE_T}\}}^{(L^2)}$ of Wick-Itô Wiener integrals  \textit{wrt}  $G$.
%
\end{coro}

%
%

\begin{pr}
The equality  $\{ \int^T_{0} f(s) \  \di G_s, \   f  \in \ \cE_T \} = \{ \int^T_{0} f(s) \  \delta G_s, \   f  \in {\cE_T}\}$  is obvious, in 
view of \eqref{dzepidzeàidjze}. Besides, the equality 
$\{ \int^T_{0} f(s) \  \delta G_s, \   f  \in \ \cH_T \} = \overline{\{  \int^T_{0}   f(s) \ \delta G_s, \ f  \in {\cE_T} \}}^{(L^2)}$ results from Meyer inequalities (see \cite[(5)]{nualart} for example).
\end{pr}

 
\begin{rem}
\label{ozeidjzeoijdzeo}
 ${\text{{\bfseries{1.}}}}$\ In many cases (such as for fBm) the equality $\overline{\mathcal{E}_{T}}^{{\|\ \|}_{T}}= \cH_T \cap L^2([0,T])$ is clear. In these situations Theorem \ref{oerferoifjeoifjefi} in a way that makes clear how the wick Itô integral \textit{wrt} $G$ generalizes the divergence one. Indeed, when equality  $\overline{\mathcal{E}_{T}}^{{\|\ \|}_{T}}= \cH_T \cap L^2([0,T])$ holds one just has to assume that $H_1)$, $H_2)$ and $H_3)$ hold. Hence, for any process $u$ in $L^2(\Omega, L^2([0,T]))$, if $u$ belongs to the domain of the divergence of $G$, then $u$ is Wick-Itô integrable on $[0,T]$  \textit{wrt}  $G$. Moreover Equality \eqref{dzepidzeàidjze} holds.
 
${\text{{\bfseries{2.}}}}$\  In view of the previous corollary, we see 
that one just has to extend the notion of Wiener integral given in 
Definition \ref{WienWick}, and call Wiener integral \textit{wrt}  
$G$ in $\mathscr{G}$, any element of $\Theta_G$, if one wants that
our set of Wiener integrals is the same that the one of \cite{nualart}. 

 ${\text{{\bfseries{3.}}}}$\   If Theorem \ref{oerferoifjeoifjefi} clearly states 
 that the Wick-Itô integral has a bigger set of integrands than the 
 divergence type integral developed in \cite{nualart}, assuming they both 
 belong to $L^2(\Omega,\overline{\mathcal{E}_{T}}^{{\|\ \|}_{T}})$, one may wonder if this fact 
 remains true outside  $L^2(\Omega,\overline{\mathcal{E}_{T}}^{{\|\ \|}_{T}})$. While this remains an open problem, here is what we can still say about it.
%
  The set $\cH_T$ may contains generalized functions (for example, one can see \cite[p.280]{Nu} or \cite[Proposition 2.11]{JLJLV1} in the case where $G$ is a fBm). When this happens ({\textit{i.e.}} when, for almost every $\omega$ in $\Omega$,  $u(\omega)$ is a generalized function which belongs to $\cH_T$ and  which is not a function), $\int^T_0 u_s\ \delta G_{s}$ has still a meaning and belongs to $L^2(\Omega)$.
 On the contrary, $\int^T_0 u_s \ \di G_{s}$ can only  be defined if $s\mapsto u_s$ is a function (an 
$(\cS^*)$-valued function but still a function). Define the space
 \vspace{-1ex}
 
 $$\Lambda:=\big\{ u \in  L^2(\Omega;\cH_T); \ u {\text{ is Wick-Itô integrable \textit{wrt} }} G  {\text{ and such that  }} {{\int^T_0}} u_s \ \di G_{s} \in L^2(\Omega)  \big\}.$$

A consequence of what we stated above is that the inclusion $\Dom(\delta_G) \subset \Lambda$ is not true. 
Note that the inclusion $\Lambda \subset \Dom(\delta_G)$ does not hold 
either. Indeed, if one considers again, as process $G$, the fBm $B^H$, 
as we did  in Example \ref{zopdijzeoidzeoidjzeoidjzeoidjz}, we know that 
$B^H$ belongs to $\Lambda$ for every $H$ in $(0,1)$, while $B^H$ 
does not belong to $\Dom(B^H) $ if $H$ is in $(0,1/4)$. Finally, the only 
thing one can say in general is that we have the dense inclusion 
$L^2(\Omega,L^2([0,T]))\cap \Dom(\delta_G) \subset \Lambda$.
\end{rem}

\subsubsection{Comparison with the divergence type integral of \cite{MV05}}
\label{zeofkepofkpo123879987}

The comparison between Wick-Itô stochastic integral and the one defined in \cite{MV05} is easier, in view of Theorems \ref{pzoefkpoezzokfzepokzpzefj1286564536} and \ref{oerferoifjeoifjefi}. Indeed, one has the following result.
\begin{theo}
\label{oerferodezdezdifjeoifjefi}
For any process $u$ such that the (extended) Skorohod integral  
 \textit{wrt}  $G$, on $[0,T]$, defined in \cite{MV05}, exists, then $u
$ is also Wick-Itô integrable \textit{wrt}  $G$, on $[0,T]$. Moreover one has the equality
\begin{equation}
\label{dzepiddezdzzeàidjze}
\cite{MV05} - \int^T_0 u_s\ \delta G_{s} = \int^T_0 u_s \ \di G_{s},
\end{equation}
where \cite{MV05} - $ \int^T_0 u_s\ \delta G_{s}$ denotes the 
(extended) Skorohod integral of $u$ \textit{wrt} $G$, defined in \cite{MV05}.
\end{theo}
\begin{rem}
\label{peorfkerpofkeprofkepr}  Note that, in this case, one does not have to make any additional assumptions  (such as H1, H2 or H3) nor that the equality $\overline{\mathcal{E}_{T}}^{{\|\ \|}_{T}}= \cH_T \cap L^2([0,T])$ holds. Moreover, and  as we stated in the introduction, this theorem as well as Theorem \ref{pzoefkpoezzokfzepokzpzefj1286564536} show that the stochastic integral \textit{wrt} to $\mathscr{V}_{\gamma}$ - processes,  built in 
\cite{MV05}, is a 
particular case of the Wick-Itô stochastic 
integral we provide here. This means that for every $\mathscr{V}_{\gamma}$ - process 
${\widetilde{B}}^{\gamma}$, and every stochastic process 
$X$, such that the integral of $X$  \textit{wrt} ${\widetilde{B}}^{\gamma}$ exists in the 
sense defined in \cite{MV05}, the Wick-Ito stochastic integral of $X$  \textit{wrt} 
${\widetilde{B}}^{\gamma}$ exists. Moreover they are equal. Finally, this also allows us to deal with non continuous Gaussian processes, as it is the case in \cite{MV05}. 

\end{rem}

\begin{pr}
 Using notations of Theorem \ref{pzoefkpoezzokfzepokzpzefj1286564536}, Lemma \ref{dede} and Cauchy-Schwarz inequality one gets, for every $\eta$ in $\sS(\R)$ and every integer $q\geq 3$:
\vspace{-1.5ex}

\begin{equation*}
\int_0^T | S(u_s)(\eta)\ S(W^{(G)}_s)(\eta)|\  \ ds  \leq {\bigg(\int^{T}_{0}\ {\|u_{s}\|}^{2}_{0} \  \ ds\bigg)}^{1/2} \ {\bigg(\int^{T}_{0}\ {|\Phi'(s)|}^{2}_{-q} \  \ ds\bigg)}^{1/2} \  e^{{|\eta|}^2_{q}}.
\end{equation*}

Since  both quantities $\left\|  u \right\|
^2_{L^2(\Omega,L^2([0,T]))}$ and $ \int^{T}_{0}\ {|\Phi'(s)|}^{2}_{-q} \  \ ds$ are finite (by assumption for the first one and as a consequence of Theorem \ref{pzoefkpoezzokfzepokzpzefj1286564536} ofr the second one), Theorem \ref{peodcpdsokcpodfckposkcdpqkoq} applies and establishes the existence of  $\int^T_0 u_s \ \di G_{s}$.
Besides, since in this case the extended domain of the  \cite{MV05} - Skorohod integral 
is, by its very definition (see \cite[Def. 27]{MV05}) a subset of 
$L^2(\Omega,L^2([0,T])$, one can use the exact same proof 
as the one of Theorem 
\ref{oerferoifjeoifjefi}; one just has to change therein $K_{-}$ (\textit{resp.} 
$K_{+}$) by $K^{*}_{\gamma}$ (\textit{resp.} $K^{*,
\text{a}}_{\gamma}$) and note that the equality given in \cite[Remark 12]{MV05} has now the role played by  
Equality \eqref{peofpoferpoferpofepfoj}, in the proof of Theorem \ref{oerferoifjeoifjefi}. The only thing which remains to be shown 
is that $\sS(\R) \subset \cH'$, where $\cH':=\{f\in L^{2}([0,T]),\ 
K^{*,\text{a}}_{\gamma} f \in L^{2}([0,T])    \}$. This latter inclusion 
results from \cite[Proposition 15]{MV05} (one just has to take therein $
\eta(s):=s^{\alpha}$ on $\R^{*}_{+}$, with $\alpha\in (1/2,1)$ and $
\eta(0):=0$ and then show that $\sS(\R) \subset \cC^{\eta}$, where $
\cC^{\eta}$ has been defined in Example \ref{oevieroivjeodzdez}. 
\end{pr}
Note that the results provided in both Thm. \ref{oerferodezdezdifjeoifjefi} and Thm.  
\ref{eprfjerpoifjeroi} allow us to think that one could solve some linear 
stochastic evolution equations driven by infinite dimensional Gaussian 
processes. 
%
%
\smallskip

In \cite{MV05} the set of Gaussian processes is smaller than $\mathscr{G}$. In \cite{NuTa06} the class of Gaussian processes considered 
is a little bit restrictive\footnote{ Moreover, while the Wick product is used to define a stochastic integral in \cite{NuTa06}, the space of stochastic distributions (which is the natural set on which one can use Wick product) is not used at all. This latter is crucial to derive occupation time formulas for local times, as we will show in  \cite{JL17-2}}, (see \cite[(2), (3) \& (4)]{NuTa06}). Moreover, since our stochastic calculus is carried out within the framework of the White Noise Theory, our stochastic 
integral does not have to be extended\footnote{Besides, if \cite{cn05} provides a method (that has been used in \cite{MV05}) to extend the divergence 
type integral \textit{wrt} fBm, this leads to require
much more regularity on the function 
$f$, to provide an Itô formula\footnote{In both Itô formulas provided in these papers (\cite[Lemma 4.3]
{cn05} and \cite[Theorem $31$]{MV05}), $f$ is assumed to be of class 
$C^{\infty}$, and such that all its derivatives have a sub-exponential 
growth.}.}, once it has been built, in order that the set of integrands is not empty or not too small, as it is the case for divergence type integral (see \cite[Remark 25 \& p.407]{MV05} and \cite{cn05}). 
Indeed, it happens that the Gaussian 
process is not even itself integrable \ie that $\int  G_{s} \ \delta G_{s}$ does not exist, (\textit{e.g.} in \cite{nualart} when $G$ is a fBm with $H\leq1/4$ or when $G$ is the process considered in \cite{MV05}). Note that the same phenomenon happens also\footnote{In this latter case, the extended domain and the initial one are not comparable (see \cite[p. 383]{LN12}). } in \cite{LN12}. A general way to 
extend the divergence integral for 
Volterra processes, assuming it exists, has been provided 
in \cite{LN05}. However\footnote{and if one excepts \cite{MV05}, the results of which we fully generalize in this paper.}, an Itô formula for extended divergence integral has 
not been provided in the same time for general Volterra processes.  
Finally, our stochastic calculus is an extension to general Gaussian 
processes of the stochastic calculus  built,  \textit{wrt} fBm in 
\cite{ell,bosw,ben1} and \textit{wrt}  mBm in \cite{JLJLV1,JL13,LLVH}. 
\subsection{Comparison with Itô Integral}
\label{oervoerveiorj}
%
%
%
The goal of this section is to compare the Wick-Itô 
integral \textit{wrt} $G$ to the Itô integral \textit{wrt} $G$, when $G$ is a  
(Gaussian) semimartingale. In this subsection we still assume that $
\rR=[0,T]$. Since the line of reasoning we are 
following would be similar if $t\mapsto G_{t}$ would not be continuous, we 
will assume, in this subsection, that $G$ is continuous. Denote, for every $t$ in $[0,T]$, ${\cU}_{t}$ the complete\footnote{If ${\cU}
_{t}$ is not complete, we complete it and still denote it ${\cU}_{t}$.} $\sigma
$-field defined by ${\cU}_{t}:=\sigma(\{G_{s};\ 0\leq s\leq t \})$ and denote $\cU$ the filtration ${(\cU_{t})}_{t\in[0,T]}$.
%
%
%
%
%
%
In this subsection one then assumes that $G={(G_{t})}_{t\in[0,T]}$ is a 
continuous centered Gaussian $\cU$-semimartingale of the form \eqref{erofkorekp}, which fulfills Assumption \Aa. Let us recall first the following result, that 
describes the structure of Gaussian semimartingales.

\begin{prop}{\cite[Prop. $2$ \& Thm $1$]{Sti83}}
\label{pzodkzpeodkzob}
The Gaussian $\cU$-semimartingale $G$ is a special $\cU$-semimartingale: \ie for almost every $(t,\omega)$ in $[0,T]\times\Omega$, one can write:
\begin{equation}
\label{ozfjeroijero}
 G_{t}= M_{t} + A_{t},
\end{equation}
where $M:={(M_{t})}_{t\in[0,T]}$ is a centered $\cU$-martingale and $A:={(A_{t})}_{t\in[0,T]}$ is a centered $\cU$-predictable  process of bounded variations. Moreover, $M$ and $A$ both belong to the same Gaussian Hilbert space as $G$.
In addition, the function of quadratic variation of $G$, denoted $t\mapsto {<\hspace{-0.1cm}G\hspace{-0.1cm}>}_{t}$, is deterministic and $M$ is bounded in $L^{p}$, for every positive real $p$.
\end{prop}

Denote $\cT:={({\cT}_{t})}_{t\in[0,T]}$ the filtration, defined by ${\cT}_{t}:=
\sigma(\{B_{s};\ 0\leq s\leq t \})$, which we suppose complete (if it is not the 
case we complete it and still denote it ${\cT}_{t}$). Through this subsection, we will denote   $\ds \mathscr{I}_{G}(X):= \int^T_0\ X_s \ d G_{s}$  the Itô  (\textit{resp.}  $\ds
\mathscr{J}_{G}(X):= \int^T_0\ X_s \ \di G_{s}$ the Wick-Itô) integral of $X$ \textit{wrt} $G$ 
on $[0,T]$, when it exists.
For any continuous 
martingale $M:={(M_{t})}_{t\in[0,T]}$, bounded in $(L^{2})$ and such that 
$M_{0}=0$, denote $L^{2}(M):=L^{2}([0,T]\times \Omega, \cP, 
d\mu\ d{<M>}_{s})$ the space of progressively measurable processes $K$ such that:
\begin{equation*}
 {\|K\|}^{2}_{L^{2}(M)}:=\E[\int^{T}_{0}\  K^{2}_{s} \ d{<\hspace{-0.12cm}M\hspace{-0.15cm}>}_{s}],
\end{equation*}
where $\cP$ denote the progressive $\sigma$-field with respect to $
\cU$. The following result will be used in order to establish Point $(i)$ of  Proposition \ref{zpfok11} below.
\begin{lem}
\label{zpfodks}
Let $G$ be a Gaussian martingale that fulfills Assumption \Aa and let $f$ be a $\cC^{1,2}([0,T]\times\R,\R)$ function. Denote $F(t,x):=\tfrac{\partial f}{\partial x}(t,x)$ and define $X_{t}:=\tfrac{\partial f}{\partial x}(t,G_{t})$. If $R$ and $f$ both fulfill conditions of Theorem \ref{ergoigfjerfijoijgfoifjgoidfjgofd}, then one has the following equality:
\begin{equation*}
\text{a.s.}\hspace{0.25cm}  \int^T_0 X_s \ d G_{s} = \int^T_0 X_s \ \di G_{s}.
\end{equation*}
In particular one has the equality:
\begin{equation}
\label{dzpzpdkzpeodk12} 
\text{a.s.}\hspace{0.25cm}  \int^T_0 \ F(s,G_s) \ d G_{s} = \int^T_0 \ F(s,G_s) \ \di G_{s}.
\end{equation}
\end{lem}
{\bfseries Proof of Lemma \ref{zpfodks}.}
$G$ being a Gaussian martingale, one gets $<\hspace{-0.1cm}G\hspace{-0.1cm}>_{t} \ = R_{t}$ almost surely, for every $t$ in $[0,T]$. Since $R$ and $f$ both fulfill conditions of Theorem \ref{ergoigfjerfijoijgfoifjgoidfjgofd}, one gets, using both Itô formulas  \eqref{ojzdojzdodjz} and \cite[Theorem $3.3$]{RY}): $\Delta_{G}(X) :=  \mathscr{I}_{G}(X)- \mathscr{J}_{G}(X) =0$,
%
%
which in particular implies that: $  \int^T_0 \ F(s,G_s) \ d G_{s} = \int^T_0 \ F(s,G_s) \ \di G_{s}$.
%
$\hfill \square$

The main result of this subsection is the following.
\begin{prop}
\label{zpfok11}
{\bfseries 1.} Assume that $G$ is Gaussian martingale, adapted to the filtration $\cT$. Let $X$ in $L^{2}(G)$ such that $X_{t}\in  
(L^{2})$, for every $t$ in $[0,T]$, and such that $(X,G)$ satisfies condition $
\iI$ (given in Section \ref{ozieffjioefjoezifjioezzjioefjfeoiz}). 
Define, for every $s$ in $[0,T]$, $U_{s}:=\inf\{t;\ R_{t}>s\}$. 
Let us write the following conditions:
\bit
\ita the map $t\mapsto R_{t}$ is strictly increasing and continuous on $[0,T]$ and such that $U$ is absolutely continuous on $[0,T]$.
\itb $G_{t_{2}} - G_{t_{1}}$ is independent of $\cT_{t_{1}}$, for every $0\leq 
t_{1} < t_{2} \leq T$.
\eit
 If one the two conditions a) or b) is fulfilled then the map $t\mapsto X_{t}$ is both $dG$-integrable and Itô-integrable, on 
$[0,T]$. Moreover we have the equality:
\begin{equation}
\label{eofijero}
 \int^T_0 X_s \ d G_{s} = \int^T_0 X_s \ \di G_{s}. 
\end{equation}
{\bfseries 2.} If the semimartingale $G$ is not a Gaussian martingale, then Equality \eqref{eofijero} does not hold in general, assuming both integral $\ds \int^T_0 X_s \ d G_{s}$ and $\ds \int^T_0 X_s \ \di G_{s}$ do exist.
\end{prop}

%
%
Note that the condition $X_{t}\in  (L^{2})$, for every $t$ in $[0,T]$ is only a slight  reinforcement of the assumption $X \in L^{2}(G)$.
%
\begin{pr}{\bfseries 1.} Let $G$ and $X$ be of the form described in point {\bfseries 1.} above. We still denote   $ \mathscr{I}_{G}(X):= \int^T_0\ X_s \ d G_{s}$  the Itô  (\textit{resp.}  $
\mathscr{J}_{G}(X):= \int^T_0\ X_s \ \di G_{s}$ the Wick-Itô) integral of $X$ \textit{wrt} $G$ 
on $[0,T]$, when it exists. The existence of $ \mathscr{I}_{G}(X)$ is clear and the existence of $\mathscr{J}_{G}(X)$ is obvious, in view of Theorem \ref{dazdzedzedzedzedzedze}.
%
%
%
%
%

{\bfseries 1.a)} $G$ being a martingale, it can be written, according to Dubins-Schwarz theorem as $G_{u}=B_{R_{u}}$, for some Brownian motion $B$ and for every $u$ in $[0,T]$. Of course we also have the equalities $B_{t}=G_{U_{t}}$ and $U_{R_{t}}= t$, for every $t$ in $[0,T]$ since $R$ is strictly increasing. It is then clear that:
\begin{equation*}
\mathscr{I}_{G}(X)= \int^T_0\ X_s \ d G_{s} = \int^T_0\ X_s \ d B_{R_{s}} = \int^T_0\ X_{U_{R_{s}}} \ d B_{R_{s}} =  \int^{R_{T}}_0 \ X_{U_{t}} \ d B_{t} = \int^{R_{T}}_0 \ X_{U_{t}} \ d^{\diamond} B_{t}.
\end{equation*}
For every $\eta$ in $\sS(\R)$, we can then write $S(\mathscr{I}_{G}(X))(\eta) =\int^{R_{T}}_{0}\ S(X_{U_{t}})(\eta)  <\delta_{t}, \eta> dt$.
Besides, we have the equality $<\i1_{[0,t]},\eta>_{L^{2}(\R)}= \E[B_{t} <.,\eta>] = \E[G_{U_{t}} <.,\eta>]=<g_{U_{t}},\eta>_{L^{2}(\R)}$. Since $U$ is absolutely continuous on $[0,T]$ so is $s\mapsto <g_{U_{s}},\eta>$. Thus this yields to:
\begin{equation*}
S(\mathscr{I}_{G}(X))(\eta)=\int^{R_{T}}_{0}\ S(X_{U_{t}})(\eta) \ d<g_{U_{t}},\eta>.
\end{equation*}
On the other hand, we have: $S(\mathscr{J}_{G}(X))(\eta)= \int^{T}_{0}\ S(X_s)(\eta) \   d<g_{s}, \eta>=: \int^{T}_{0}\ f_{\eta}(s) \   dA^{(\eta)}_{s},
$
where we have set $A^{(\eta)}_{s}:=<g_{s}, \eta>$ and $f_{\eta}=S(X_.)(\eta)$. Applying \cite[Proposition 4.10]{RY} to the positive and then the negative part of $ f_{\eta}$ we get $\int^{T}_{0}\ f_{\eta}(s) \   dA^{(\eta)}_{s} = \int^{R_{T}}_{0}\ f_{\eta}(U_{s}) \   dA^{(\eta)}_{U_{s}}$ which entails that:
\begin{equation*}
S(\mathscr{J}_{G}(X))(\eta)=\int^{R_{T}}_{0}\ f_{\eta}(U_{s}) \   dA^{(\eta)}_{U_{s}} = \int^{R_{T}}_{0}\ S(X_{U_{t}})(\eta) \ d<g_{U_{t}},\eta>
\end{equation*}
and ends the proof in this case.

{\bfseries 1.b)}
 The proof of the equality $\ds \mathscr{I}_{G}(X) = \mathscr{J}_{G}(X)$ is 
obtained by following exactly the same three steps as in the proof of 
\cite[Theorem $13.12$]{Kuo2}, in which the equality $\ds \mathscr{I}_{B}(X) = \mathscr{J}_{B}(X)$ is established ($B$ being a Brownian motion). One then just has to substitute in there the process $\varphi$ by $X$, to replace $\i1_{[t_{1},t_{2})}$ by $g_{t_{2}} - g_{t_{1}}$ in 
the first step, and noticing that one can find, for any process $X$ in $L^{2}
(G)$, a sequence ${(X_{n})}_{n\in\N}$ of simple processes such that $\lim_{n
\to+\infty}\mathscr{I}_{G}(X_{n}) =\mathscr{I}_{G}(X)$, where the 
convergence holds in $(L^{2})$.  We will only write down here the first of these three steps; in order, first to make clear the differences with the case where $G=B$, and second, to translate the proof of 
\cite[Theorem $13.12$]{Kuo2} in the notations we use in this paper.
Let $(t_{1},t_{2})$ be in ${[0,T]}^{2}$ such that $0\leq t_{1} < t_{2} \leq T$ and assume that $X_{t}:=X_{t_{1}} \ \i1_{(t_{1},t_{2}]}$, where $X_{t_{1}}$ is   ${\cU}_{1}$-mesurable. Let $X_{t_{1}}=\sum^{+ \infty}_{n = 0}\  
I_{n}({f}^{(n)}_{t_{1}})$ be the chaos decomposition of $X_{t_{1}}$. By definition of Itô integral, and using the identity:$I_{n}(k) \ I_{1}(l) \  = I_{n+1}(k\otimes l) + n \ I_{n-1}({<k,l>}_{L^{2}(\R)})$,
where $\widehat{\otimes}$ denotes the symmetric tensor product and 
which is valid for every: $n$ in $\N^{*}$, symmetric function $k$ in $L^{2}(\R^{n})$ and $l$ in $L^{2}(\R)$, one gets:
\vspace{-5ex}
\begin{multline*}
\mathscr{I}_{G}(X) = X_{t_{1}} (G_{t_{2}} - G_{t_{1}})=\sum^{+ \infty}_{n = 0}\  
I_{n}({f}^{(n)}_{t_{1}})\ I_{1}(g_{t_{2}}- g_{t_{1}}) \\
 =\sum^{+ \infty}_{n = 0}\  \big( I_{n+1}({f}^{(n)}_{t_{1}} \widehat{\otimes} (g_{t_{2}}- 
g_{t_{1}}))  +n\  I_{n-1}({<{f}^{(n)}_{t_{1}},g_{t_{2}}-g_{t_{1}}>}_{L^{2}(\R)}) \ \big).    
\end{multline*}
Lemma $3.11$ of \cite{Kuo2} applies here since\ $\cU$ is included in $\cT$. One then knows that, for every $n$ in $\N^{*}$, ${f}^{(n)}_{t_{1}}$ is equal to $0$ almost everywhere on ${[0,T]}^{n}\backslash {[0,t_{1}]}^{n}$.
Moreover, since $G_{t_{2}} - G_{t_{1}}$ is independent of $\cT_{t_{1}}$, it is clear that ${<\hspace{-0.1cm}{f}^{(n)}_{t_{1}},g_{t_{2}}-g_{t_{1}}\hspace{-0.1cm}>}_{L^{2}(\R)} =0$, for every $n$ in $\N^{*}$.
Using Proposition \ref{qmqmqmmmmm}, we get, for any $\eta$ in $\sS(\R)$, 
\vspace{-0.7cm}

\begin{align*}
S(\mathscr{I}_{G}(X) )(\eta) &= S(\sum^{+ \infty}_{n = 0}\  \big( I_{n+1}({f}^{(n)}_{t_{1}} \widehat{\otimes} (g_{t_{2}}- 
g_{t_{1}})))(\eta) = \sum^{+ \infty}_{n = 0}\   <{f}^{(n)}_{t_{1}} \widehat{\otimes} (g_{t_{2}}- g_{t_{1}}), \eta^{\otimes (n+1)}> \notag\\
&= \sum^{+ \infty}_{n = 0}\ <{f}^{(n)}_{t_{1}}, \eta^{\otimes n}>\ <g_{t_{2}}- g_{t_{1}},\eta> =   S(X_{t_{1}})(\eta) \ S(G_{t_{2}}-G_{t_{1}})(\eta) =S(\mathscr{J}_{G}(X))(\eta).
\end{align*}

The injectivity of S-transform then allows us to write $\mathscr{I}_{G}(X)  = \mathscr{J}_{G}(X)$.

%
%
%
%
%

%
%
\smallskip

%
%
{\bfseries $\textit{2}$}.
In view of Proposition \ref{pzodkzpeodkzob}, three cases are possible for the structure of the semimartingale $G$. The case where $G$ is a martingale has been treated in Point {\bfseries 1} below. Our goal here is to exhibit, when $G$ is not a Gaussian martingale, some general and simple examples for which $\ds \mathscr{I}_{G}(X)$ and $\ds \mathscr{J}_{G}(X)$
both exist and are different.
Let $f$ be a $\cC^{1,2}([0,T]\times\R,\R)$ function. Assume that both $R$ 
and $f$ fulfill conditions of Theorem \ref{ergoigfjerfijoijgfoifjgoidfjgofd}. 
Denote $Y\equiv 0$ when the process $Y:={(Y_{t})}_{t\in[0,T]}$ is such that 
\vspace{-0.25cm}

\begin{equation}
\label{ofrjeoijreo}
Y_{t}(\omega) = 0, \hspace{0.25cm} \forall (\omega,t)\in\Omega'\times [0,T], 
\end{equation}

where $\Omega'$ is measurable subset such that $\mu(\Omega')=1$. We will denote $Y\not\equiv 0$ when \eqref{ofrjeoijreo} is not satisfied.


%
The case where $M\equiv 0$ being easier we will assume that 
 $A\not\equiv0$  and $M\not\equiv0$  in \eqref{ozfjeroijero}.
Assume that both $M$ and $A$ are continuous and that there exists a map $t \mapsto g^{(1)}_{t}$ from $[0,T]$ into $
L^{2}(\R)$ such that $g^{(1)}$ fulfills Assumption \Aa
and such that $M_{t}= <.,g^{(1)}_{t}>$ almost surely, for every $t$ in $[0,T]$. 
Let us compare $\mathscr{I}_{G}(G)$ and  $\mathscr{J}_{G}(G)$. The existence of $\mathscr{I}_{G}(G)$ is clear. Moreover, using classical Itô formula, one gets:
\begin{equation*}
\mathscr{I}_{G}(G)= \int^{T}_{0} \ M_{s}\ dM_{s} + \int^{T}_{0} \ A_{s}\ dA_{s}  + A_{T}\ M_{T}.
\end{equation*}
The existence of $\mathscr{J}_{G}(G)$ is clear in view of Example 
\ref{zeodijze}. Moreover, using again Example \ref{zeodijze}, the fact that $M$ is bounded in $(L^{2})$ as well as an 
integration by parts, one gets:
\begin{equation*}
\mathscr{J}_{G}(G)= \int^{T}_{0} \ M_{s}\ \di M_{s} + \int^{T}_{0} \ A_{s}\ \di A_{s}  + A_{T}\diamond M_{T}.
\end{equation*}
%
Classical Itô formula, Example \ref{zeodijze}, 
\eqref{dzpzpdkzpeodk12} and, finally,  Propsotion \ref{pzodkzpeodkzob} and 
 Proposition \ref{qmqmqmmmmm} yields:
\vspace{-2ex}
\begin{align}
\Theta_{G}(G)  &:= \mathscr{I}_{G}(G) - \mathscr{J}_{G}(G) =\int^{T}_{0} \ 
A_{s}\ dA_{s} - \int^{T}_{0} \ A_{s}\ \di A_{s}  + A_{T}\ M_{T} -  A_{T}
\diamond M_{T} \notag \\
& = 2^{-1} \E[A^{2}_{T}] + A_{T}\ M_{T} -  A_{T}\diamond M_{T} =  2^{-1} 
\E[A^{2}_{T}] - \E[A_{T} M_{T}].
\end{align}
It is then easy to find a finite variation processes $A$, as well as a positive real $T$, and choose the map $g^{(1)}$, that defines the Gaussian martingale $M$, such that:
$2^{-1} \E[A^{2}_{T}] - \E[A_{T} M_{T}]\neq 0$.
\end{pr}

\begin{rem}

${\text{{\bfseries{1.}}}}$ In the particular case where there exists a function $f:\R\rightarrow \R$, which belongs to $L^{2}(\R)$, such that $\ds G_{t}:=\int^{t}_{0}\ f(u) \ dB_{u}$, for every $t$ in $[0,T]$, \textit{a.s.}, then all the assumptions of Propoistion \ref{zpfok11} are reduced to $X$ belongs to $L^{2}(G)$. Note also that one recovers in particular, the result  of \cite[Theorem 13.12]{Kuo2},  that is $\ds  \int^1_0 X_s \ \di B_{s} = \int^1_0 X_s \  dB_{s}$, for every $X$ in $L^{2}(B)$. 

${\text{{\bfseries{2.}}}}$ 
%
One may also remark, from what we stated in the previous sections, that the existence of  $ \int^1_0 X_s \ \di G_{s}$  does not imply the existence  
$ \int^1_0 X_s \  dG_{s}$. Conversely the existence of $ \int^1_0 X_s \  dG_{s}$ does not imply the existence $ \int^1_0 X_s \ \di G_{s}$. Three natural questions  then arise in this framework;
\vspace{-1ex}

\bit
\iti If $ \int^T_0 X_s \  dG_{s}$ exist, on which conditions on $X$ the integral $ \int^T_0 X_s \ \di G_{s}$ will exist?
\itii If $ \int^T_0 X_s \  \di G_{s}$ exist, on which conditions on $X$ the integral $ \int^T_0 X_s \ dG_{s}$ will exist?
\itiii When both the integrals $ \int^T_0 X_s \  \di G_{s}$ and  $ \int^T_0 X_s \  dG_{s}$ do exist, what is the exact link bet ween them?
\eit

%
\vspace{-2ex}
In order to answer properly to these three questions one needs to use the 
operators $D_{g_t}$ and $D^*_{g_t}$, defined in \cite[Chap $9$]
{Kuo2}, and express both our Wick-Itô integral and the Itô integral using these operators. Since this would lead us too far from the goal of this present work, we will therefore give the answer to these questions in a future work.

${\text{{\bfseries{3.}}}}$ 
In  view of Lemma \ref{zpfodks} 
, it seems that Equality \eqref{eofijero} remains true under weaker assumptions than the one proposed in Proposition \ref{zpfok11}. However, extend \eqref{eofijero} under weaker assumptions is an open problem.

${\text{{\bfseries{4.}}}}$ 
%
 Of course one can limit our definition of Wick-Itô integral to Gaussian martingales only. Then, and as it is the case for Itô integral, one can extend the definition of Wick-Itô integral \textit{wrt} $G$ to the case where $G$ is a Gaussian semimartingale, by simply setting:
 \vspace{-1ex}
\begin{equation}
\label{efopj123}
\int^{T}_{0}\ X_{s}\ d^{*} G_{s}:= \int^{T}_{0}\ X_{s}\ \di M_{s} + \int^{T}_{0}\ X_{s}\  dA_{s}, 
\end{equation}
 \vspace{-2ex}
 
where $M$ (\textit{resp.} $A$) denotes the martingale  (\textit{resp.} the bounded variation process) given by \eqref{ozfjeroijero}, and where $M$ is assumed to be of the form \eqref{erofkorekp} and fulfills Assumption \Aap.
 $\int^{T}_{0}\ X_{s}\ d^{*} G_{s}$ will then be defined as soon as each member of the right hand side of \eqref{efopj123} will exist. The Itô integral $\int^{T}_{0}\ X_{s}\  dA_{s}$, in the right hand side of \eqref{efopj123}, offers also the advantage, on  $\int^{T}_{0}\ X_{s}\  \di A_{s}$ of being defined $\omega$ by $\omega$ since it is Stieljes integral.
 \end{rem}

\vspace{-4ex}

\begin{small}
\section*{Acknowledgments}
I want to express my deep gratitude to Jacques Lévy Véhel for his advices
and for the very stimulating discussions we had about this work. I also want to thank Professor T. Hida for his warm welcome at the University of Nagoya, where a part of this paper was written, as well as Professor L. Chen and the Institute for Mathematical Sciences of Singapore (NUS), where another part of this paper was written. I also thanks the Associate Editor as well as the anonymous referee for his remarks that greatly improve the quality of this paper and especially Section \ref{pokprofkdeezdzedezdedee}.

This work is dedicated to the memory of Professor Marc Yor. 
\end{small}

\vspace{10ex}


\centerline{\Large \bfseries \textsc{Appendix}}
\vspace{-2ex}
\makeatletter
\renewcommand\theequation{\thesection.\arabic{equation}}
\@addtoreset{equation}{section}
\makeatother

\begin{appendix}

\makeatletter
\renewcommand\theequation{\thesection.\arabic{equation}}
\@addtoreset{equation}{section}
\makeatother

\vspace{-10ex}
\textcolor{white}{
\section{Appendix}
\subsection[\textcolor{black}{Bochner Integral}]{\textcolor{white}{Bochner Integral}}
\label{appendiceB}
}
\vspace{-3ex}

\section*{ A.  Bochner Integral}
\vspace{-1ex}
\begin{small}
The following notions about Bochner integral come from \cite[p.$72$, $80$ and $82$]{HP} and \cite[p.$247$]{Kuo2}.

\begin{defi}{Bochner integral \cite[p.$247$]{Kuo2}}
\label{bb}
Let $I$ be a Borelian subset of $\R$ endowed with the Lebesgue measure. One says that  $\Phi:I \rightarrow \ooS$ is Bochner integrable on $I$ if it satisfies the two following conditions:

$1$  $\Phi$ is weakly measurable on $I$ \textit{i.e} $u\mapsto <\hspace{-0.2cm}<\hspace{-0.1cm} \ \Phi_{u},\varphi\  \hspace{-0.1cm}>\hspace{-0.2cm}>$  is measurable on $I$ for every $\varphi$ in $(\cS)$.
\medskip

$2$ $\exists\ p \in \bN$ such that $\Phi_{u} \in ({\cS}_{-p})$ for almost every $u \in I$ and $u \mapsto {\|\Phi_{u}\|}_{-p}$ belongs to $L^1(I)$.
\medskip

\ \ The Bochner-integral of $\Phi$ on $I$ is denoted $\int_{I}\ \Phi_{s} \ ds $ .
\end{defi}

\begin{prop}
If $\Phi:I \rightarrow \ooS$ is Bochner-integrable on $I$  then there exists an integer $p$ such that ${\left\| \int_{I} \Phi_{s}  \ ds \right\|}_{-p} \leq \int_{I} {\left\|\Phi_{s} \right\|}_{-p} \ ds$. 
Moreover $\Phi$ is also Pettis-integrable on I and both integrals coincide on I.
\end{prop}

\begin{rem}
The previous proposition shows that there is no risk of confusion by using the same notation for both
 Bochner and Pettis integrals.
\end{rem}

\begin{theo}
\label{cc}
Let $p \in \bN$ and ${(\Phi^{(n)})}_{n \in \bN}$ be a sequence of processes from $I$ to \oS such that  $\Phi^{(n)}_{u} \in ({\cS}_{-p})$ for almost every $u \in I$ and for every $n$. Assume moreover that $\Phi^{(n)}$ is  Bochner-integrable on $I$, for every $n$, and that $\lim\limits_{(n,m) \to (+\infty,+\infty)}  \int_{I}  {\big\| {\Phi}^{(m)}_{s} - {\Phi}^{(n)}_{s}\big\|}_{-p} \  ds = 0$. Then there exists an \oS-process (almost surely $({\cS}_{-p})$-valued), denoted $\Phi$, defined and Bochner-integrable on $I$, such that 

\vspace{-0.3cm}

\begin{equation}
\label{iuerhferuhfrufhu}
\lim\limits_{n \to +\infty}  \int_{I}\  {\|\Phi_{s} - {\Phi}^{(n)}_{s}\|}_{-p} \ ds = 0
\end{equation}

Furthermore, if there exists an \oS-process, denoted $\Psi$, which verifies $\eqref{iuerhferuhfrufhu}$, then $\Psi_{s} = \Phi_{s}$ for {\textit a.e.} $s$ in $I$. Finally one has  $\lim\limits_{n \to +\infty}  \int_{I}\ \Phi^{(n)}_s \ ds =  \int_{I}\ \Phi_{s} \ ds$, where the equality and the limit both hold in \oS.   
\end{theo}
\end{small}

\textcolor{white}{
\section{Appendix}
\subsection[\textcolor{black}{Proof of Theorem \ref{pzoefkpoezzokfzepokzpzefj1286564536}}]{\textcolor{white}{Bezozeijdezoijl}}
\label{appendice2}
}

\vspace{-18ex}

\section*{B.  Proof of Theorem \ref{pzoefkpoezzokfzepokzpzefj1286564536}}
\vspace{-1ex}

\begin{small}

{\bfseries Proof:}
\label{eorfijerofijeroij}
In view of Proposition \ref{CS},  it is sufficient to show that Assumption \Da holds. Besides, it is clear that $\Phi$ is well defined on $\R_{+}$ since one has, for every $t$ in $\R_{+}$, the equality:
\vspace{-1ex}
\begin{equation}
\label{eoifjoerijoefjoij}
 \Phi_{t}= E(t) \cdot \delta_{0} - (\i1_{[0,t)} \cdot E(t-\cdot))'.
\end{equation}
It is clear that function $E$ (\textit{resp.} $\mathscr{E}$) is increasing, differentiable on $\R^{*}_{+}$ and continuous on $\R_{+}$ (\textit{resp.} increasing and of class $C^{1}$ on $\R_{+}$). Equality \eqref{eoifjoerijoefjoij} together with the properties of $E$
 entail that $\Phi$ is continuous at $t=0$. Let us now establish Equality \eqref{erfjeorfjeorfjeroiferofije}. For every $t$ in $\R^{*}_{+}$, $\varphi$ in ${\sS}(\R)$, and $r>0$, 
denote $I_{r}:= <\tfrac{\Phi_{t+r}-\Phi_{t}}{r},\varphi>$. Using the 
change of variable formula, an easy computation gives us: 
\vspace{-1.25ex}
\begin{align*}
 \hspace{-2ex}I_{r} &= \tfrac{1}{r} \left( \int^{t+r}_{0} \ \varphi(u) \  \varepsilon(t+r-u) \ du - \int^{t}_{0} \ \varphi(u) \  \varepsilon(t-u) \ du  \right)\\
& \hspace{-2ex}= \int^{1}_{0}  \hspace{-1ex}\tfrac{t}{r}  \left[ \varepsilon\big((t+r)(1-v)\big)\ \varphi(v(t+r)) - \varepsilon(t(1-v))\ \varphi(vt) \right] \hspace{-0.5ex}dv + \hspace{-1ex} \int^{1}_{0}  \hspace{-1ex}\varepsilon\big((t+r)(1-v)\big)\ \varphi(v(t+r))  dv=:I^{(1)}_{r} + I^{(2)}_{r}.
\end{align*}
\vspace{-3ex}

For every $r$ in $(0,1)$, one has $I^{(2)}_{r} = \int^{t+1}_{0} 
\i1_{(0,t+r)}(u) \ \frac{\varepsilon(u)}{t+r}\ \varphi(t+r-u) \ du$. Since 
$t$ and $\varepsilon$ are positive, Lebesgue's dominated convergence 
theorem applies and allows one to write that $ \lim_{r\to 0} I^{(2)}_{r} = \int^{t}_{0} 
 \frac{\varepsilon(u)}{t}\ \varphi(t-u) \ du$ and thus\footnote{Note that one could also have used \cite[Remark $3$]{MV05} and assume that $(\gamma^{2})'$ (and hence $\varepsilon$) is non-increasing.} that $\lim_{r\to 0} I^{(2)}_{r} = \frac{1}{t}\int^{t}_{0} \ \varphi(u)\ \varepsilon(t-u) \ du$.
%
 Besides, $I^{(1)}_{r}$ can be written under the following form:
\vspace{-3ex}

\begin{equation*}
I^{(1)}_{r} =  \int^{1}_{0} t  \varepsilon\left((t+r)(1-v)\right) 
\big(\tfrac{\varphi(v(t+r))-\varphi(vt)}{r}\big) dv +  \int^{1}_{0} \frac{t}{r}\ \varphi(vt)  \left(\varepsilon\big((t+r)(1-v)\big) - \varepsilon\left(t(1-v)\right) \right) 
 dv =:J^{(1)}_{r} + J^{(2)}_{r}.
\end{equation*}
The exact same method as the one used to compute $\lim_{r\to 0} I^{(2)}_{r}$ applies and allows one to write:
\begin{equation}
\label{perofkerpdededeok}
 \lim_{r\to 0} J^{(1)}_{r} =  \int^{1}_{0} \ tv\
 \varphi'(vt)\ \varepsilon(t(1-v)) \ dv =  \frac{1}{t}\int^{t}_{0} \ u\
 \varphi'(u)\ \varepsilon(t-u) \ du.
\end{equation}

Having in mind that $\varepsilon^{2}(r) = (\gamma^{2})'(r)$, an integration by parts in $J^{(2)}_{r}$ yields:
\vspace{-3ex}

%

\begin{align*}
J^{(2)}_{r} &=   \varphi(0)  \left( \frac{t}{t+r} \frac{(E(t+r) - E(t))}{r} - \frac{E(t)}{t+r} \right) -\frac{t}{t+r} \ \int^{1}_{0}\ \varphi'(vt)\ E((t+r)(1-v))  \ dv\\
&\hspace{2ex}+\underbrace{\frac{t}{r} \int^{1}_{0} \ \varphi'(vt)\  \left(E((t+r)(1-v)) - E(t(1-v))\right) \ dv}_{=:K_{r}}.
\end{align*}

It is then clear that:
\vspace{-3ex}

\begin{equation}
\label{perofkerpdededeok2}
 \lim_{r\to 0} J^{(2)}_{r} 
 =  \varphi(0)  \left(\varepsilon(t) - \frac{E(t)}{t}\right) - \frac{1}{t} \int^{t}_{0} \
 \varphi'(u)\ E(t-u) \ du + \lim_{r\to 0} K_{r}.
\end{equation}

Thus, it only remains to determine $\ds\lim_{r\to 0} K_{r}$. An  integration by parts in $K_{r}$ yields:
\vspace{-2ex}

\begin{equation*}
 K_{r} 
 = \ \varphi'(0)  \bigg(  \tfrac{t}{t+r}\tfrac{\mathscr{E}(t+r)-\mathscr{E}(t)}{r} - \tfrac{\mathscr{E}(t)}{t+r} \bigg) -\tfrac{t}{t+r}  \int^{1}_{0} \varphi''(vt) \  \mathscr{E}(t(1-v))  dv + \tfrac{t^{2}}{t+r}  \int^{1}_{0} \varphi''(vt) \left(  \tfrac{\mathscr{E}((t+r)(1-v)) - \mathscr{E}(t(1-v))}{r}  \right)   dv.
\end{equation*}


One then gets, after a change of varaible:
\begin{equation}
\label{eofijeorijeroifjeoij}
 \hspace{-2ex} \lim_{r\to 0} K_{r} 
= \varphi'(0)  \big( E(t) - \tfrac{\mathscr{E}(t)}{t} \big) - \frac{1}{t}\int^{t}_{0} \ \varphi''(u) \  \mathscr{E}(t-u) \ du + \int^{t}_{0} \ \varphi''(u) \ (1-\frac{u}{t})\ E(t-u) \ du.
\end{equation}

Finally, gathering $\lim_{r\to 0} I^{(2)}_{r}$ and Equalities \eqref{perofkerpdededeok} to \eqref{eofijeorijeroifjeoij} yields:
\vspace{-3ex}

\begin{multline*}
 \lim_{r\to 0} I_{r} =  \frac{1}{t}\int^{t}_{0} \
 \varphi(u)\ \varepsilon(t-u) \ du + \frac{1}{t} \int^{t}_{0} \
 \varphi'(u)\ (u\ \varepsilon (t-u)- E(t-u)) \ du \\
 \hspace{2ex}+ \frac{1}{t} \int^{t}_{0} \
 \varphi''(u)\ ((t-u)\ E(t-u)- \mathscr{E}(t-u)) \ du 
 \hspace{2ex}+\varphi(0)  \left(\varepsilon(t) - \frac{E(t)}{t}\right) +\varphi'(0)  \bigg( E(t) - \frac{\mathscr{E}(t)}{t} \bigg).
 \end{multline*}
%
This is nothing but \eqref{erfjeorfjeorfjeroiferofije}. Let us now show that $t\mapsto {|\Phi'(t)|}_{-q}\in \underset{\ \ b\in\R^{*}_{+}}{\cap} L^2((0,b))$, $\forall q\geq 3$. 
Let $b$ be a positive real and $q$ be an integer such that $q\geq 3$.
 It is sufficient to show that the map  $t\mapsto {|\Phi'(t)|}
^{2}_{-q}$ belongs to $L^1((0,b))$. Using 
\eqref{erfjeorfjeorfjeroiferofije}, one gets, for every integer $k\geq 
2$,
\begin{equation}
\label{eofijdze12eorijdedezeroifjrreoij}
{|<\Phi'(t),e_{k}>|} \leq \alpha_{k} \big( \varepsilon(t) + \tfrac{E(t)}{t} + M\big),
 \end{equation}
where $\alpha_{k}:= {\|e_{k}\|}_{\infty} + {\|e'_{k}\|}_{\infty}+ {\|e''_{k}\|}_{\infty}$ and where $M:=(t+3)\underset{s\in[0,t]}{\sup} E(s) + \underset{s\in[0,t]}{\sup} \mathscr{E}(s)+ \underset{s\in(0,t]}{\sup} \frac{\mathscr{E}(s)}{s}$. 
Using the relation $e'_{k}(x)=\sqrt{\frac{k}{2}} e_{k-1}(x)- \sqrt{\frac{k+1}{2}} e_{k+1}(x)$ (see \cite[p.$354$]{Kuo2}) as well as Theorem \ref{ozdicjdoisoijosidjcosjcsodijqpzeoejcvenvdsdlsiocfuvfsosfd}, one easily obtains that  $\alpha_{k}\leq 48 \ {(k+1)}^{2} \sum^{2}_{l=-2} \ {\|e_{k+l}\|}^{2}_{\infty} \leq 250 \ {(k+1)}^{2}$, for every  integer $k\geq 2$. Having in mind the definition of $\cR_{n}$ given in \eqref{zoijeroif} it is then clear that there exists $C>0$ which does not depend on $q$ nor $t$ such that:
\vspace{-1ex}
\begin{equation}
{|\Phi'(t)|}^{2}_{-q}=\sum^{+\infty}_{k=0}{|<\Phi'(t),e_{k}>|}^{2}\ {(2k+2)}^{-2q} \leq \ell(t)\cdot  \sum^{+\infty}_{k=0}\ \alpha^{2}_{k}\ {(2k+2)}^{-2q} \leq C \cdot \ell(t)\cdot \cR_{2-q},
\end{equation}
\vspace{-3ex}

where we have set $\ell(t):= {(\varepsilon(t) + \tfrac{E(t)}{t} + M)}^{2}$.
It then remains to show that both $\varepsilon^{2}$ and $t\mapsto  
\frac{{E(t)}^{2}}{t^{2}}$ belong to $L^1((0,b))$. The first part is 
clear since $\varepsilon^{2} = (\gamma^{2})'$. Moreover, in a 
neighborhood of $0$ one has\footnote{See \cite[p.396]{MV05} for a proof of this fact.} 
$E(t)\leq 2\ t\ \varepsilon(t)$, for $t\neq 0$. One therefore has: $\int^{b}_{0} {\big(\tfrac{E(t)}{t}\big)}^{2} dt \leq 4 \int^{b}_{0} \varepsilon^{2}(t)  dt = 4\gamma^{2}(b) <+\infty$,  which ends the proof.
$ \hfill \square$
\end{small}

\end{appendix}

\begin{scriptsize}

\end{scriptsize}

\end{document}